\documentclass[12pt, reqno]{amsart}
\usepackage[T1]{fontenc}
\usepackage{amsmath,amsthm,amscd,amssymb,amsfonts, amsbsy}

\makeatletter
\newcommand{\sectionnotoc}[1]{%
  \begingroup
  \let\@tocwrite\@gobbletwo
  \section*{#1}
  \endgroup
}
\makeatother

\usepackage{mathabx}
\usepackage{url}
\usepackage{latexsym, color, enumerate}
\usepackage{pxfonts}
\usepackage{mathrsfs}
\usepackage{marginnote}
\usepackage{todonotes}
\usepackage[colorlinks=true, pdfstartview=FitV, linkcolor=blue, citecolor=blue, urlcolor=blue]{hyperref}
\usepackage{mathtools}
\mathtoolsset{showonlyrefs}
\usepackage[margin=1 in]{geometry}
\usepackage{comment}

\allowdisplaybreaks[4]
\usepackage{lipsum}

\newcommand\blfootnote[1]{%
  \begingroup
  \renewcommand\thefootnote{}\footnote{#1}%
  \addtocounter{footnote}{-1}%
  \endgroup
}
\usepackage[numbers,sort]{natbib}

\usepackage{appendix}

\newtheorem{innercustomgeneric}{\customgenericname}
\providecommand{\customgenericname}{}
\newcommand{\newcustomtheorem}[2]{%
  \newenvironment{#1}[1]
  {%
   \renewcommand\customgenericname{#2}%
   \renewcommand\theinnercustomgeneric{##1}%
   \innercustomgeneric
  }
  {\endinnercustomgeneric}
}

\newcustomtheorem{customthm}{Theorem}

\usepackage{tikz}
\usetikzlibrary{decorations.pathmorphing}
\usepackage{multicol}
\usepackage{cancel}

\numberwithin{equation}{section}

\theoremstyle{plain}
\newtheorem{theorem}{Theorem}[section]
\newtheorem{lemma}[theorem]{Lemma}

\newtheorem{proposition}[theorem]{Proposition}

\theoremstyle{definition}
\newtheorem{definition}[theorem]{Definition}

\newtheorem{remark}[theorem]{Remark}

\newtheorem{hp}{}

\newtheorem*{theorem*}{Theorem}

\newcommand{\p}{\partial}
\newcommand{\ds}{{\frac{d}{ds}}}

\newcommand{\re}{\mathrm{Re}}

\newcommand{\Op}{\mathrm{Op}}

\newcommand{\dt}{{\frac{d}{dt}}}
\newcommand{\la}{\langle}
\newcommand{\ra}{\rangle}

\DeclareMathOperator\supp{supp}
\newcommand{\abs}[1]{\lvert#1\rvert}

\def\XXint#1#2#3{{\setbox0=\hbox{$#1{#2#3}{\int}$}
		\vcenter{\hbox{$#2#3$}}\kern-.5\wd0}}

\newcommand{\R}{\mathbb{R}}

\begin{document}

\title[]{Microlocal and global smoothing effect of variable coefficient ZK-type operators}

\author[Federico and Tramontana]{Serena Federico and Davide Tramontana}

\address{Serena Federico
\newline \indent Dipartimento di Matematica, Universit\`a di Bologna
\newline \indent Piazza di Porta San Donato 5, 40126, Bologna, Italy}
\email{serena.federico2@unibo.it}

\address{Davide Tramontana
\newline \indent Dipartimento di Matematica, Universit\`a di Bologna
\newline \indent Piazza di Porta San Donato 5, 40126, Bologna, Italy}
\email{davide.tramontana4@unibo.it}

\blfootnote{S. Federico and D. Tramontana are members of GNAMPA.}

\begin{abstract}
In this work, we study the smoothing effect of a class
of Zakharov--Kuznetsov (ZK) equations, which generalize the
Korteweg--de Vries (KdV) equation to higher dimensions. Specifically, we consider linear evolution equations generated by ZK-type operators with variable coefficients, and we study both their microlocal and global smoothing effects under suitable conditions on the coefficients of the operator. Our analysis requires a detailed study of the geometric properties of the Hamiltonian flow associated with the operator under consideration. These properties allow us to prove a dispersive and a strong smoothing effect associated with the Cauchy problem in question: under suitable moment conditions on the initial datum, the solution becomes smooth in space for positive times in a finite time interval.
\end{abstract}

\maketitle 
\tableofcontents

\section{Introduction}
The aim of this work is to study the microlocal smoothing effect associated with a linear initial value problem of the form 
\begin{equation}\label{mainprobintro}
\begin{cases}
    D_tu-Au=0 \quad \text{in} \ \ (0,T]\times \R^n, \\
    u(0,\cdot)=u_0 \in L^2(\R^n),
\end{cases}
\end{equation}
where $T>0$ and $A$ is a differential operator of order three of the form
\begin{align}
    A(x,D)=D_1A_2(x,D),\quad A_2(x,D)=\sum_{i,j}D_ia_{ij}(x)D_j,
\end{align}
with $A_2$ being a variable coefficients differential operator with principal symbol
\[
a_2(x,\xi):=\sum_{i,j=1}^n a_{ij}(x)\xi_i\xi_j,
\qquad (x,\xi)\in\mathbb{R}^{2n}.
\]
We assume that $a_2$ is uniformly elliptic and that
$a_{ij}(x)-\delta_{ij}$, together with its derivatives,
satisfies suitable decay estimates as $|x|\to\infty$. The latter condition encodes the fact that $A_2$ has symbol which is asymptotic to the symbol of the positive Laplacian at spatial
infinity (see Section \ref{sec.ZK}).

The existence of solutions to Cauchy problems of the form
\eqref{mainprobintro} has been established in \cite{Ch}
(see Section \ref{sec.smoothing} for a more detailed discussion).
We may therefore focus on the study of propagation phenomena
for this class of equations.

The operators considered in this work are (multidimensional)
variable-coefficient generalizations  of the linear
Zakharov--Kuznetsov (ZK) operator
\[
D_t-D_{x_1}\Delta,
\qquad
\Delta=\sum_{j=1}^n D_{x_j}^2,
\]
where $\Delta$ denotes the positive Laplacian.
This operator describes the dispersive part of the ZK equation,
a multidimensional generalization of the Korteweg--de Vries (KdV) equation.
The physical motivation behind this work comes from the study of ion-acoustic
waves in magnetized plasmas, whose propagation is described,
in suitable regimes, by the ZK equation (see \cite{MP,Se,KPSV}).
The operator $D_{x_1}\Delta$ represents its dispersive part:
the anisotropy in equation  in the $x_1$ direction reflects that the plasma waves propagate more
freely along magnetic field lines. This is the core of the linear dispersive dynamics
underlying this model.
Studying the smoothing properties of the corresponding evolution operator clarifies
how these dynamics affect the regularity of solutions.
In the present work, we investigate this mechanism for a class
of variable-coefficient ZK-type operators, focusing on the
relation between spatial decay of the initial datum and
microlocal regularity at positive times. In this framework, replacing the Laplacian with a second order variable-coefficient (in space) operator
allows us to model spatial variations in the dispersive properties
of the medium, extending the constant-coefficient description
to a spatially inhomogeneous setting.

More precisely, we consider operators of the form
\begin{equation}\label{eq.operatorsintro}
D_t-A=D_t-D_{x_1}A_2,
\qquad
A_2=\sum_{i,j=1}^n D_{x_i}a_{ij}(x)D_{x_j},
\end{equation}
where the coefficients are smooth and real-valued, and the matrix
$(a_{ij}(x))_{i,j=1}^n$ satisfies suitable hypotheses (see \ref{hp1}-\ref{hp3}).

\vspace{2mm}

We investigate the local-in-time smoothing effect for Cauchy problems of the form \eqref{mainprobintro}, from both microlocal and global perspectives. 
Our results are of two different type: of microlocal dispersive type and of purely microlocal type. 

From the microlocal dispersive point of view, we show that if the initial datum $u_0$ belongs to $L^2(\R^n)$
then the homogeneous solution $u$ gains microlocally one derivative, averaged in time and localized in space, with respect to the initial datum. This is more precisely expressed through the property
\begin{equation}\label{intro.dispersive.smooth}
   \re \int_0^T (u(t,\cdot),C u(t,\cdot))_{L^2}dt<\infty
    \end{equation}
with $C:=\Op(c)$ such that $c\in S(2,-1;0,1)$ is microlocally elliptic in a conic neighborhood at high frequency of $(x_0,\xi_0)$ (see the notation in Section \ref{sec.pre}). The fact that $c\in S(2,-1;0,1)$ makes the property look like a gain of regularity in weighted Sobolev spaces averaged in time. This result, derived in the local framework in \cite{C}, is also obtained in higher generality under the moment condition. More precisely we prove that, if $u_0$ satisfies a microlocal moment condition of order $K\geq 0$ integer at a point $(x_0,\xi_0)$ (see Definition \ref{moment.condition}), then there exist operators $C_\ell=\Op(c_\ell)$, with $c_\ell\in S(\ell,K-\ell-1;\rho,\delta)$, $\ell=1,\ldots,K$ microlocally elliptic in a conic neighborhood at high frequency of $(x_0,\xi_0)$, such that \eqref{intro.dispersive.smooth} holds with $C=C_\ell$ at positive times (see Theorem \ref{thm.mainest}). This can be viewed as a microlocal gain of weighted Sobolev regularity of type $H^s(\mathbb{R}^N;\langle x\ra^{K-\ell-1}dx)$. We stress that the requirement for $u_0$ is a moment condition, which differs from the usual Sobolev regularity condition required in classical local dispersive estimates. Moreover, the more we gain in Sobolev regularity the more we lose in terms of moment property. 

From the purely microlocal side, we prove that if the initial datum $u_0$ belongs to $L^2(\R^n)$ and satisfies suitable microlocal moment conditions at a point $(x_0,\xi_0)$ (see Definition \ref{def.finitemoments}), then the solution $u(t,\cdot)$ gains Sobolev regularity microlocally at $(x_0,\xi_0)$ for every $t\in(0,T]$. The amount of regularity gained depends on the order of the moment conditions imposed on the initial datum at that point. As a consequence of this microlocal result, we also establish a global smoothing effect: if the moment conditions are satisfied globally, then the corresponding gain in Sobolev regularity holds globally for every $t\in(0,T]$. In particular, in the global setting, we have that if the inital datum $u_0 \in L^2(\R^n)$ satisfies the moment condition at each order, then the solution $u$ of \eqref{mainprobintro}, with initial datum $u_0$, satisfies $u(t,\cdot) \in C^\infty(\R^n)$, for each $t \in (0,T]$.

\vspace{2mm}

The smoothing effect of operators of the form $D_t-A$ is closely related with the nontrapping property of the bicharacteristic curves associated with (the principal symbol of) the operators $A$ under consideration (see \cite{F2}). Under the nontrapping condition, high-frequency waves leave compact regions quickly enough to gain local or microlocal regularity when averaged over time.
The relationship between smoothing and nontrapping for
Schr\"odinger-type equations has also been investigated on
Riemannian manifolds (see \cite{D1,D2}), including settings with boundary (see \cite{B}). For the Laplace-Beltrami operator, whose principal symbol is given in terms of the metric, this relationship is expressed in terms of the nontrapping property of the geodesic flow (since the configuration
projections of the bicharacteristics are geodesics, up to a
constant rescaling of the flow parameter).

The strength of our result relies on the fact that for operators of the form \eqref{eq.operatorsintro}, the nontrapping property holds globally (see Lemma \ref{lemma.nontrapping}).

In the global setting, smoothing effects for variable-coefficient Schr\"odinger-type operators have been extensively investigated using a variety of approaches and techniques (see, e.g.,  \cite{D0,KPRV05,FR,FS,F1,AACG,GT}). By contrast, considerably fewer results are available for operators of order higher than two (see \cite{FT,Ch, AAC1, AAC2, AAC3}), due in part to the more involved symbolic properties of the operator and of its bicharacteristics.

A pioneering work on the microlocal smoothing effect for
Schr\"odinger-type equations is \cite{CKS}. In this work, the
authors relate the microlocal regularity of solutions at positive
times to moment conditions on the initial datum along backward
nontrapped bicharacteristics. This approach highlights the
connection between the geometry of the Hamiltonian flow,
the behavior of the initial datum at spatial infinity, and
the regularity of the solution.
Building on this work, the microlocal smoothing effect was subsequently studied in \cite{W} for Schr\"odinger type equations in the framework of scattering manifolds. The main motivation behind this geometric framework is that it encompasses asymptotically Euclidean manifolds (see \cite{W} for further details) and thus provides the natural generalization, in non-Euclidean settings, of the variable-coefficient Schr\"odinger operators studied in \cite{CKS}, \cite{KPRV05} and \cite{FT}, which were themselves obtained as asymptotically small perturbations of constant-coefficient differential operators.

Further developments on scattering manifolds were obtained
in \cite{IN}, where the microlocal smoothing was studied using
classical scattering maps and radially homogeneous wave front
sets. In the analytic setting, an analytic theory of the
quadratic scattering wave front set was developed and applied
to the Schr\"odinger equation in \cite{RZ}, while microlocal analytic smoothing estimates for Schr\"odinger
operators with analytic coefficients that are long-range
perturbations of the Laplacian were established in \cite{MNS}. In the analytic setting, microlocal smoothing results have been obtained also for operators of order $m\geq 3$ in \cite{Ta}. To the best of our knowledge, the latter is the only result concerning the microlocal smoothing effect for operators of order greater than two. 

\vspace{2mm}

In order to study the smoothing effect we used the scattering calculus, which is a pseudodifferential calculus
that allows for possibly different growth rates in the configuration
variable $x$ and the frequency variable $\xi$. In its standard form,
it is the Weyl calculus associated with the Weyl--H\"ormander metric
(see \cite{HoV3}, Section 18.5)
\[
g_{SC}
=
\frac{\abs{dx}^2}{1+\abs{x}^2}
+
\frac{\abs{d\xi}^2}{1+\abs{\xi}^2}.
\]
This calculus has been studied both in the Euclidean setting
(see \cite{P,C,C0}) and on manifolds
(see, e.g., \cite{M,MZ,Sc}).
Here, we use the scattering calculus developed in \cite{CKS},
which generalizes the standard scattering calculus on $\R^{2n}$
by allowing the parameters $\rho$ and $\delta$ to differ from
$0$ and $1$, respectively.

\vspace{2mm}

The main microlocal smoothing results for Cauchy problems of the form \eqref{mainprobintro}  we obtain in this work, under a suitable finite moment condition on the initial datum (see Definition \ref{def.finitemoments}), read as follows (recall Subsection \ref{sub.micro} for the notation).  Note that Theorem \ref{thm.mainest} is stated in great generality in terms of a nonnegative real number $K$.

\begin{theorem}\label{thm.mainest}
    Let $u_0 \in L^2(\R^n)$, with $S_0=\Op(s_0)$, $s_0=\sqrt{b_0}$, with $b_0$ as in \eqref{b0c0}, and assume 
    \begin{equation}\label{m.moment.condition}
        (u_0,S_0^\ast S_0u_0)_{L^2}<+\infty. 
    \end{equation}
    Let also $K\geq 0$ and denote by  $\kappa=K-\lfloor K \rfloor$, where $\lfloor K \rfloor$ is the integer part of $K$. 
    Then, for all $\ell=0,\dots,\lfloor K \rfloor$ there exist
    symbol pairs $(b_\ell,c_\ell)$, with $0 \leq b_\ell \in S(\ell+\kappa,K-\ell-\kappa ;\rho,\delta)$ $0 \leq c_\ell \in S(\ell+\kappa+2,K-\ell-\kappa-1;\rho,\delta)$, $0\leq \rho<\delta\leq 1$, and there exists a nested sequence of subneighborhoods 
    \[
    \mathscr{D}_{+}^{\lfloor K \rfloor}(x_0,\xi_0) \subseteq \mathscr{D}_+^{\lfloor K \rfloor-1}(x_0,\xi_0)\subseteq \dots \subseteq \mathscr{E}_+^{0}(x_0,\xi_0), 
    \]
    defined by 
    \[
    \mathscr{D}_{+}^\ell(x_0,\xi_0)=\lbrace (x,\xi); \ c_\ell(x,\xi) \geq \langle x \rangle^{K-\ell-\kappa-1}\langle \xi \rangle^{\ell+\kappa+2} \rbrace , \quad \ell=1,\dots,\lfloor K \rfloor -1,
    \]
    and for $\varepsilon>0$ small,
    \[
    \mathscr{D}_{+}^{\lfloor K \rfloor}=\lbrace (x,\xi); \ c_{\lfloor K \rfloor}(x,\xi) \geq \langle x \rangle^{-1-\varepsilon}\langle \xi \rangle^{K} \rbrace, 
    \]
    such that, for each $\ell\geq 1$, 
    \[
     \mathrm{supp}(b_\ell)\subseteq \mathscr{D}_{+}^{\ell-1}(x_0,\xi_0).
    \] 
    Furthermore, for all $\ell=1,\dots, \lfloor K \rfloor$, for all $T>0$, on denoting $B_\ell=\Op(b_\ell)$, $C_\ell=\Op(c_\ell)$ and $u(t,\cdot)$ the solution of \eqref{mainprobintro} we have 
    \begin{align}
        (i) \ \sup_{ 0 \leq t \leq T}\Bigl(t^{\mu_\ell}\re(u(t,\cdot),B_\ell u(t,\cdot))_{L^2}\Bigr)<+\infty \\
        (ii) \ \re\int_0^T t^{\mu_\ell}(u(t,\cdot),C_\ell u(t,\cdot))_{L^2}dt<+\infty,
    \end{align}
    where $\mu_\ell=\ell+\kappa(1+\varepsilon')$, with $\varepsilon'>0$.
    
    If $\ell=0$ then $\mathscr{D}_+^0(x_0,\xi_0)\equiv \mathscr{E}_+^J(x_0,\xi_0)$, with $\mathscr{E}_+^J(x_0,\xi_0)$ as in Theorem \ref{thm.nestedOperators}, and (i)-(ii) hold.
\end{theorem}

\begin{theorem}\label{theo.WFs}
Let $(x_0,\xi_0) \in \R^{2n}\setminus 0$, let $u_0 \in L^2(\R^n)$ and let $u$ be the solution of \eqref{mainprobintro} with initial datum $u_0 \in L^2(\R^n)$. Then, if $u_0$ has finite moment of order $s \in \R$ along the forward bicharacteristic starting from $(x_0,\xi_0)$, for all $t \in (0,T]$ we have that $(x_0,\xi_0) \notin WF^s(u(t,\cdot))$.
\end{theorem}
In particular, if the moment condition holds globally (see \eqref{eq.globalmomentK} and \eqref{eq.globalmoment} below), then, as a consequence of Theorem \ref{theo.WFs}, we obtain the following global result, which reveals a remarkably strong regularizing property for equations of the form \eqref{mainprobintro}.
\begin{theorem}\label{theo.glob}
Let $s\in [0,+\infty)$ and $u_0 \in L^2(\R^n)$ be such that 
\begin{equation}\label{eq.globalmomentK}
\langle x \rangle^s u_0 \in L^2(\R^n)
\end{equation}
Then, if $u$ solves \eqref{mainprobintro} with initial datum $u_0$, we have $u(t,\cdot) \in H_{\mathrm{loc}}^s(\R^n)$ for all $t \in (0,T]$. In particular, if \begin{equation}\label{eq.globalmoment}
\langle x \rangle^s u_0 \in L^2(\R^n), \quad \forall s \in [0,+\infty),
\end{equation}
we have $u(t,\cdot) \in C^\infty(\R^n)$ for all $t \in (0,T]$. 
\end{theorem}

\vspace{2mm}

The proof of these results relies on a deep geometric analysis
in phase space tailored to (the principal symbol of) our $ZK$
operators \eqref{eq.operatorsintro}. This analysis focuses on the
bicharacteristic flow associated with the principal symbol $a_3$
of $A$ in \eqref{eq.operatorsintro}. The first fundamental properties
are established in Lemmas \ref{lemma.nontrapping} and \ref{lemmaCKS},
which show, respectively, that under our hypotheses
\ref{hp1}-\ref{hp3} all points in the phase space are nontrapped
(see Definition \ref{nontrapped}), and that, under suitable
assumptions on the initial point $(x,\xi)$, the bicharacteristic
flow of $A$ behaves like the flow of the constant-coefficient operator $D_{x_1}\Delta$.
The latter result allows us to define the scattering regions
$\mathscr{S}_{\pm}$ (cf. Definition \ref{def.scatteringsets}),
which, roughly speaking, are the regions in the phase space where
the bicharacteristic flow behaves like its constant-coefficient
counterpart. We can then define the scattering variables
(see Lemma \ref{lem.scatteringV}) and use them to construct,
for each $(x_0,\xi_0)$, a pair of symbols $(b,c)$
(see Subsection \ref{subsec.bc}) supported in a neighborhood
of the forward bicharacteristic through $(x_0,\xi_0)$.
These symbols allow us to apply a positive commutator argument
in Section \ref{sec.smoothing} and obtain our smoothing results
using the a priori estimates provided in our scattering framework.

\vspace{2mm}

The paper is organized as follows.
In Section \ref{sec.pre} we recall the tools from pseudodifferential
calculus and microlocal analysis used throughout the paper.
In Section \ref{sec.ZK}, we introduce the class of operators under
consideration and describe its relation to the constant-coefficient
ZK operator.
Section \ref{sec.bicZK} is devoted to the analysis of the associated
bicharacteristic flow. Its geometric properties guide the construction
of the phase-space localizers used in the smoothing estimates.
In Section \ref{sec.smoothing}, we establish the microlocal smoothing
results and derive their consequences in the global setting.
In Section \ref{sec.futper}, we discuss some questions that
we plan to investigate in future work.
Finally, the appendix collects several auxiliary results used
in the proofs.

\medskip 

\textbf{Notation.} We work in $T^\ast \R^n \cong \R^{2n}$ and write $\R^{2n}\setminus 0=\R^n \times (\R^n \setminus \lbrace 0 \rbrace)$. 

We denote by $\mathbb{N}_0$ the set $\mathbb{N}_0=\mathbb{N} \cup \lbrace 0 \rbrace$. 

We use the notation $x \cdot y=\sum_{j=1}^nx_jx_j$ for the standard euclidean product on $\R^n$, $D_j$ for the first order differential operator $D_j=i^{-1}\partial_{y_j}$, $\delta_{ij}$ for the Kronecker delta with index $i$ and $j$,  and $I=(\delta_{ij})_{i,j=1}^n$ for the identity matrix. 

If $f,g:I\rightarrow \mathbb{R}$, with $I\subseteq \mathbb{R}$, we shall write $f(s) \approx g(s)$ whenever there exists a universal constant $C>0$ such that $C^{-1}g(s)\leq f(s)\leq Cg(s)$, for all $s \in I$. Furthermore, we write $f(s) \lesssim g(s)$ (resp. $f \gtrsim g$) if there exists a constant $C>0$ such that $f(s) \leq Cg(s)$ (resp. $g(s) \geq Cf(s)$) for all $s\in I$. More generally, for real-valued functions $a,b$ defined on more general sets, we write $a \lesssim b$ (resp. $a \gtrsim b$) if there exists a universal constant $C>0$ such that $a \leq b$ (resp. $a \geq b$).

For $F,G:\R^n\rightarrow \R^m$ we write $F(x) \simeq G(x)$ as $|x|\rightarrow +\infty$, if $F_j(x)-G_j(x)\rightarrow 0$ as $|x|\rightarrow+\infty$, for all $j=1,\dots,m$.

We denote by $\mathsf{M}_{n\times m}(\R)$ the set of $n\times m$-matrices with real entries, with $\mathsf{M}_{n}(\R)$ the set of $n\times n$ matrices with real entries and then by $0_n\in \mathsf{M}_{n}(\R)$ the matrix with all zero entries. 

We say that a set $V\subseteq \R^{2n}$ is conic if $(x,\xi)\in V$ implies $(x,\lambda \xi) \in V$ for all $\lambda>0$. 
A conic neighborhood $V$ of a set $U \subseteq \R^{2n}$ is a neighborhood of $U$ which is conic. 

We say that a set $W$ is conic at high frequencies (resp. conic neighborhood of $U$ at high frequencies) if there exists a conic set $V$  (resp. $V$ a conic neighborhood of $U$) and $R>0$ such that $W\cap \lbrace |\xi|>R \rbrace= V\cap \lbrace |\xi|>R \rbrace$.

We say that a function $f : \R^{2n} \rightarrow \R$ is homogeneous of degree $k\in \mathbb{N}_0$ if $f(x,\lambda\xi)=\lambda^kf(x,\xi)$ for all $\lambda>0$ and for every $(x,\xi) \in \R^{2n}$. 

The space $C^\infty_b$ is the space of $C^\infty$ bounded functions with bounded derivatives at any order.

\section{Preliminaries}\label{sec.pre}

In this section, we recall some preliminary notions that we will use throughout the work. 

\subsection{Basics of standard pseudodifferential calculus}
We start with the definition of the so-called \textit{standard} symbol class $S^m(\R^n):=S_{1,0}^m(\R^n)$.

    Let $m \in \R$ and $a\in C^\infty(\R^{2n})$. We say that $a \in S^{m}(\R^n)$, that is that \textit{$a$ is a symbol of order $m$}, if for all $\alpha,\beta \in \mathbb{N}_0^n$ there exists a constant $C_{\alpha,\beta}>0$ such that
\[
|\partial_x^\beta\partial_\xi^\alpha a(x,\xi)|\leq C_{\alpha,\beta}\langle \xi\rangle^{m-|\alpha|},\quad  \forall (x,\xi)\in\R^{2n}.
\]
The class $S^{m}(\R^n)$ of symbols of order $m$ is a Frechét space, with seminorms given by 
\begin{align}
    \gamma_{k,m}(a)=\sup_{(x,\xi)\in\R^n \times \R^n, |\alpha|+|\beta|\leq k} |(\partial_x^\beta\partial_\xi^\alpha a)(x,\xi)|\langle \xi\rangle^{-m+|\alpha|},\quad k\in\mathbb{N}.
\end{align}

Given a symbol $a\in S^m(\R^n), m\in \R$, we  may define a \textit{pseudodifferential operator} associated with $a$ by using some formulas called \textit{quantization formulas}. In this work we deal with the so-called Kohn-Nirenberg (or left) quantization, which is defined by 
\[
\mathrm{Op}(a)u(x)= (2\pi)^{-n}\int_{\R^n \times \R^n}e^{i(x-y) \cdot \xi} a(x,\xi)u(y)dyd\xi, \quad u \in \mathscr{S}(\R^n).
\]

Hence, we say that $A$ is a \textit{standard pseudodifferential operator of order $m$}, 
and write $A \in \Psi^m(\mathbb{R}^n)$, if there exists 
$a \in S^m(\mathbb{R}^n)$ such that $A = \mathrm{Op}(a)$. 

If $a \in S^m(\mathbb{R}^n)$ for some $m \in \mathbb{R}$, 
then it is easy to verify that $\mathrm{Op}(a)$ 
defines a linear operator 
acting continuously on $\mathscr{S}(\mathbb{R}^n)$, 
and extends by duality to a continuous linear operator 
on $\mathscr{S}'(\mathbb{R}^n)$.

\subsection{The scattering calculus} In this work we deal with the scattering calculus introduced in \cite{CKS}. As we will see below, this calculus will be useful to work in the forward region (see Subsection \ref{subsec.bc} below).
We recall here the definition of the symbol class under consideration and, for the sake of completeness, also the properties of the calculus. We refer the reader to \cite{CKS} for the proof of the results stated before. 
\begin{definition}\label{def.Smkrhodelta}
Let $m,k \in \mathbb{Z}$, $\rho,\delta \in [0,1]$. We say that $a \in C^\infty(\R^{2n})$ belongs to the \textit{symbol class $S(m,k;\rho,\delta)$} if for all $\alpha,\beta \in \mathbb{N}_0^n$ there exists $C_{\alpha,\beta}>0$ such that 
\[
|\partial_x^\beta \partial_\xi^\alpha a(x,\xi)| \leq C_{\alpha,\beta} \la \xi \ra^{m-|\alpha|}\la x \ra^{k+\rho|\alpha|-\delta|\beta|}, \quad \forall (x,\xi) \in \R^{2n}.
\]
\end{definition}
\begin{remark}
    Note that with $\rho=0$ and $\delta_1$ we recover the standard scattering symbol classes $S_{SC}^{m,k}(\R^n)$ on $\R^n$ (cf.  \cite{P,C,C0,Hi}), namely classes of functions $a \in C^\infty(\R^{2n})$ satisfying
    \[
    |\partial_x^\beta \partial_\xi^\alpha a(x,\xi)| \leq C_{\alpha,\beta} \la \xi \ra^{m-|\alpha|}\la x \ra^{k-|\beta|}, \quad \forall (x,\xi) \in \R^{2n}.
    \]
\end{remark}
Definition \ref{def.Smkrhodelta} immediately yields the subsequent property of the product of symbols: 
\begin{align}
 \text{if }  &\, a \in S(m_1,k_1;\rho_1,\delta_1), \; b \in S(m_2,k_2; \rho_2,\delta_2)\\
 \text{then } & \, ab\in S(m_1+m_2,k_1+k_2;\max\{\rho_1,\rho_2\},\min\{\delta_1,\delta_2\}). \label{propertyComposition}
\end{align}

Moreover, the following composition and adjoint formulas for this symbol class hold.

\begin{theorem}
Let $a \in S(m_1,k_1;\rho,\delta)$ and $b \in S(m_2,k_2; \rho,\delta)$, with $0\leq \rho < \delta \leq 1$ and let $A=\Op(a)$ and $B=\Op(b)$. 
Then, for each $N \in \mathbb{N}$
\[
AB=E+R_{N+1},
\]
for $E=\Op(e)$, with 
\[
e=\sum_{0 \leq |\alpha| \leq N}\frac{(-i)^{|\alpha|}}{\alpha!}(\partial_\xi^\alpha a)( \partial_x^\alpha b),
\]
and some $R_{N+1}:\mathscr{S}'(\R^n)\rightarrow \mathscr{S}'(\R^n)$. Moreover, if $N$ is sufficiently large, the operator  $R_{N+1}$ is bounded on $L^2(\R^n)$.
\end{theorem}

\begin{theorem}
Let $a \in S(m,k;\rho,\delta)$, with $0\leq \rho < \delta \leq 1$ and let $A=\Op(a)$. Then, for all $N \in \mathbb{N}$, the adjoint 
\[
A^\ast=\Op(a^\ast)+R_{N+1},
\]
where 
\[
a^\ast=\sum_{0 \leq |\alpha| \leq N}\frac{(-i)^{|\alpha|}}{\alpha!}\partial_\xi^\alpha  \partial_x^\alpha \widebar{a},
\]
and some $R_{N+1}:\mathscr{S}'(\R^n)\rightarrow \mathscr{S}'(\R^n)$. Moreover, if $N$ is sufficiently large, the operator  $R_{N+1}$ is bounded on $L^2(\R^n)$.
\end{theorem}

We also have the Sharp G\aa rding inequality for this calculus, that can be stated as follows. \begin{theorem}\label{thm.GardingInequality}
Let  $b \in S(m,k; \rho,\delta)$, with $0\leq \rho<\delta \leq 1$ be such that $b \geq 0$ with smooth square root $s=b^{1/2} \in S(m/2,k/2; \rho,\delta)$. Then, there exists $e \in S(m-2,k-2(\delta-\rho); \rho,\delta)$ with $ \mathrm{supp}(e)\subseteq \mathrm{supp}(b)$ such that, on denoting $B=\Op(b)$, $S=\Op(s)$ and $E=\Op(e)$ we have 
\[
\re (\psi, B \psi)_{L^2} \geq \|S^\ast \psi\|_{L^2}^2+\re(\psi,E\psi)-C\|\psi\|_{L^2}^2, \quad \mathscr{S}(\R^n),
\]
for some constant $C>0$. If $m \leq 2$ and $k \leq 2(\delta-\rho)$, the symbol $e$ may be taken to be zero. 
\end{theorem}
As a consequence of this Theorem we also have the next two results that will be useful in Section \ref{sec.smoothing}. 
\begin{theorem}\label{thm.L2boundedness}
    Let $b \in S(0,0;\rho,\delta)$, with $0 \leq \rho < \delta \leq 1$. Then $\Op(b)$ is bounded on $L^2(\R^n)$.
\end{theorem}
\begin{theorem}\label{thm.UpperBound}
Let  $b_0 \in S(m,k; \rho,\delta)$, with $0\leq \rho<\delta \leq 1$ be such that $b_0 \geq 0$ with smooth square root $s_0=b_0^{1/2} \in S(m/2,k/2; \rho,\delta)$. Then, if $b \in S(m,k;\rho,\delta)$ satisfies 
\[
\mathrm{supp}(b)\subseteq \mathscr{E}^0=\lbrace (x,\xi); b_0(x,\xi) \geq \langle x \rangle^{k} \langle \xi \rangle^{m} \rbrace, 
\]
on denoting $B=\Op(b)$ and $S_0=\Op(s_0)$ we have 
\[
|(\psi,B \psi)_{L^2}|\leq C((\psi, S_0^\ast S_0 \psi)_{L^2}+\|\psi\|_{L^2}^2), \quad \psi \in \mathscr{S}(\R^n),
\]
for some constant $C>0$. 
\end{theorem}

\begin{remark}\label{rmk.theo2.7}
By repeating the same proof (see \cite{CKS}) under the same notation of Theorem \ref{thm.UpperBound}, it is possible to get also 
\[
|(\psi,B \psi)_{L^2}|\leq C((\psi, S_0 S_0^\ast \psi)_{L^2}+\|\psi\|_{L^2}^2), \quad \psi \in \mathscr{S}(\R^n).
\]
\end{remark}

\subsection{Microlocal analysis}\label{sub.micro}
We next recall the basics of microlocal analysis we need in our work. For further details we refer to \cite{HoV3,Hi}.
\begin{definition}\label{def.elliptic}
    Let $(x_0,\xi_0) \in \R^{2n}\setminus 0 $ and $A=\Op(a) \in \Psi^m(\R^n)$. We say \textit{$A$ is elliptic at $(x_0,\xi_0)$} if there exist $c,C>0$ and $\varepsilon>0$ such that 
    \[
    |a(x,\xi)| \geq c |\xi|^m, \quad \text{if} \ |\xi| \geq C, \ |x-x_0|+\Bigl|\frac{\xi}{|\xi|}-\frac{\xi_0}{|\xi_0|}\Bigr|<\varepsilon.
    \]
    The set of elliptic points of $A$ is denoted by $\mathrm{Ell}(A)$ and is called \textit{the elliptic set} of $A$, while its complementary $\mathrm{Char}(A)=(\R^{2n}\setminus 0) \setminus \mathrm{Ell}(A)$ is called \textit{the characteristic set} of $A$.
\end{definition}
Through this definition it is possible to introduce the notion of (standard) wave front sets of a tempered distribution, which characterize the
$C^\infty$ and $H^s_{\mathrm{loc}}$ regularity of a distribution, where recall, for $s \in \R$,
\[
H^s_{\mathrm{loc}}(\R^n)=\lbrace u \in \mathscr{D}'(\R^n); \; \phi u \in H^s(\R^n), \, \forall \phi \in C_c^\infty(\R^n) \rbrace.
\]
\begin{definition}
    Let $u \in \mathscr{S}'(\R^n)$. Then, \textit{the wave front set of $u$} is defined by 
    \[
    WF(u)=\bigcap_{\substack{A\in\Psi^0(\R^n)\\ Au\in C^\infty(\R^n)}}
\mathrm{Char}(A),
\]
and, for $s \in \mathbb{R}$, \textit{the wave front set of $u$ of order $s \in \R$}, by 
\[
    WF^s(u)=\bigcap_{\substack{A\in\Psi^0(\R^n)\\ Au\in H_{\mathrm{loc}}^s(\R^n)}}
\mathrm{Char}(A).
\]
\end{definition}
 
\begin{remark}
Notice that the definition of the wave front set can be equivalently formulated as follows: a point $\mathbb{R}^{2n}\setminus{0} \ni (x_0,\xi_0) \notin WF(u)$ if and only if there exists an operator $A\in\Psi^0(\mathbb{R}^n)$, elliptic at $(x_0,\xi_0)$, such that $Au\in C^\infty(\mathbb{R}^n)$. Similarly, for any $s\in\mathbb{R}$, a point $ \mathbb{R}^{2n}\setminus{0} \ni (x_0,\xi_0) \notin WF^s(u)$ if and only if there exists an operator $A\in\Psi^0(\mathbb{R}^n)$, elliptic at $(x_0,\xi_0)$, such that $Au\in H_{\mathrm{loc}}^s(\mathbb{R}^n)$.
\end{remark}
\begin{theorem}\label{theo.singsupp}
    Let $u \in \mathscr{S}'(\R^n)$, then if $\pi_x: \R^{2n}\rightarrow \R^n_x$ is the projection on the $x$-variables, we have 
    \[
    \pi_x (WF(u))=\mathrm{singsupp}(u).
    \]
\end{theorem}
\begin{theorem}\label{theo.Hs}
Let $u \in \mathscr{S}'(\R^n)$ and $s \in \R$. Then $WF^s(u)=\emptyset$ if and only if $u \in H_{\mathrm{loc}}^s(\R^n)$.
\end{theorem}

\begin{theorem}\label{theo.WFWFs}
Let $u \in \mathscr{S}'(\R^n)$. Then 
\begin{equation}\label{eq.WFWFs}
WF(u)=\overline{\bigcup_{s \in \R} WF^s(u)}.
\end{equation}
\end{theorem}
\begin{remark}
Note also that for $s \in \R$ the wave front set of order $s$ of a tempered distribution $u$ can be equivalently defined by 
\[
    WF^s(u)=\bigcap_{\substack{A\in\Psi^0(\R^n)\\ Au\in H^s(\R^n)}}
\mathrm{Char}(A),
\]
(if $Au \in H^s_{\mathrm{loc}}(\R^n)$ and $A$ elliptic at $(x_0,\xi_0)$ then $(\chi A)u \in H^s(\R^n)$, with $\chi$ cutoff function near $x_0$) and moreover, via the use of parametrices (see \cite{Hi}) by 
\[    WF^s(u)=\bigcap_{\substack{A\in\Psi^s(\R^n)\\ Au\in L^2(\R^n)}}
\mathrm{Char}(A),
\]
and hence $ \mathbb{R}^{2n}\setminus{0} \ni(x_0,\xi_0) \notin WF^s(u)$ if and only if there exists an operator $A\in\Psi^s(\mathbb{R}^n)$ elliptic at $(x_0,\xi_0)$ such that $Au \in L^2(\R^n)$.  
\end{remark}

\begin{remark}\label{rmk.WFK}
Note that if $u \in \mathscr{S}'(\R^n)$, and $s,s' \in \R$, with $s'<s$, we have $WF^{s'}(u)\subseteq WF^s(u)$, and so we may rewrite \eqref{eq.WFWFs} as 
\[
WF(u)=\overline{\bigcup_{s \geq 0} WF^s(u)}.
\]
\end{remark}

\subsection{The bicharacteristic flow}
In order to give the notion of bicharacteristic flow we first recall the notion of Hamilton vector field. 

If $a \in C^\infty (\R^{2n})$ is a real-valued smooth function, $H_a$ denotes the \textit{Hamilton vector field} associated with $a$ (or Hamilton vector field of $a$), which, in standard symplectic coordinates $(x,\xi) \in \R^{2n}$, can be defined as 
\begin{equation}\label{eq.defHa}
H_a:=\sum_{j=1}^n \Bigl( (\partial_{\xi_j}a)\frac{\partial}{\partial x_j}-(\partial_{x_j}a)\frac{\partial}{\partial \xi_j} \Bigr)\in \Gamma(\R^{2n};T\R^{2n}),
\end{equation}
and then, if $b \in C^\infty (\R^{2n})$,
\begin{equation}\label{eq.defPoisson}
    \lbrace a, b \rbrace(x,\xi):=H_ab(x,\xi)=\sum_{j=1}^n \Bigl( (\partial_{\xi_j}a)\partial_{x_j}b-(\partial_{x_j}a)\partial_{\xi_j}b \Bigr)(x,\xi), \quad (x,\xi)\in \mathbb{R}^{2n},
\end{equation}
where $\lbrace \cdot, \cdot \rbrace$ are known as \textit{the Poisson brackets}.

The notions of bicharacteristic flow and of $X$-bicharacteristic flow (i.e. its projection onto the configuration space) are given as follows.

\begin{definition}
    Let $a \in C^\infty(\R^{2n};\R)$ and let $H_a \in \Gamma(\R^{2n};T\R^{2n})$ be the Hamilton vector field associated with $a$. \textit{The bicharacteristic curve} associated with $a$ starting from $(x,\xi)\in  \R^{2n}\setminus 0$  is the maximal solution of the ODE
    \begin{equation}\label{eq.Hameq}
    \begin{cases}
        \dot \gamma(s)=H_a(\gamma(s)) \\
        \gamma(0)=(x,\xi).
    \end{cases}
    \end{equation}
    
    In what follows, we shall denote such a curve by $\gamma_{(x,\xi)}(s)=(X(s;x,\xi),\Xi(s;x,\xi))$, with $s \in I=I(x,\xi)\subseteq \R$ maximal interval of definition, and call \textit{$X$-bicharacteristic curve} of $a$ starting from $(x,\xi)$ its projection on $\R^n_x$, that is the curve $s\mapsto X(s;x,\xi)$. Finally, denoting by $U=\lbrace (s,x,\xi) \in \R \times \R^{2n}; \ s \in I(x,\xi) \rbrace$, \textit{the bicharacteristic flow} and \textit{the $X$-bicharacteristic flow} of $a$, are the maps
    \begin{equation}
         \varphi:U \rightarrow \R^{2n}, \quad (s,x,\xi) \mapsto (X(s;x,\xi),\Xi(s;x,\xi))\in \R^{2n} \ \ 
    \end{equation}
    and 
    \begin{equation}
        X: U \mapsto \R^{2n}, \quad (s,x,\xi) \mapsto X(s;x,\xi) \in \R^{n},
    \end{equation}
    respectively.
    
\end{definition}

Now, let $A$ be a differential operator of order $m$, written as 
\[
A=\sum_{|\alpha|\leq m}a_\alpha(x)D_x^\alpha, 
\]
with $a_\alpha\in C_b^\infty(\R^n)$, for all $|\alpha|\leq m$. The \textit{principal symbol} of $A$ is defined as 
\[
a_m(x,\xi):=\sum_{|\alpha|=m}a_\alpha(x)\xi^\alpha, \quad (x,\xi) \in \R^{2n}.
\]
We recall here the notion of trapped and nontrapped points for $A$ (cf. \cite{CKS,FT}).

\begin{definition}\label{nontrapped}
Let $A$ be a differential operator of order $m$ and let $a_m$ be its principal symbol. We say that $(x,\xi)\in \R^{2n}\setminus 0$ is \textit{forward (backward) nontrapped} for $A$, if the flow $(X(s;x,\xi),\Xi(s;x,\xi))$ associated with its principal symbol $a_m$ is globally defined, and 
\[
\lim_{s\rightarrow +\infty}\abs{X(s;x,\xi)}=+\infty \quad (\lim_{s\rightarrow -\infty}\abs{X(s;x,\xi)}=+\infty).
\]
A point $(x,\xi) \in \R^{2n}\setminus 0$ is said to be \textit{nontrapped} if it is either forward or backward nontrapped. Finally, we say that a point is \textit{trapped} if it is not nontrapped, and we denote by 
\[
S^A_{\mathrm{trap}}=\Bigl\lbrace (x,\xi) \in \R^{2n}\setminus 0; \ (x,\xi) \ \mathrm{is} \ \mathrm{trapped }\Bigr\rbrace.
\]
\end{definition}

\section{ZK-type operators}\label{sec.ZK}
The aim of this work is to investigate the microlocal smoothing effect associated with a class of third-order differential operators with variable coefficients, tailored to the well-known linear ZK operator (see \eqref{eq.defZKcc} below). In this section, we first introduce the standard variable-coefficient ZK-type operator on $\mathbb{R}^2$, which serves as our toy model. We then state our assumptions in order to treat a generalization of this operator in $\mathbb{R}^n$, with $n \geq 2$.
\subsection{The toy model}
The standard $ZK$ operator is defined as the differential operator on $\R^2$ defined by 
\begin{equation}\label{eq.defZKcc}
A_0=D_{x_1}\Delta,
\end{equation}
where $\Delta$ is the positive Laplacian $\Delta=D_{x_1}^2+D_{x_2}^2$. Note that, by using the Fourier transform, we may write $A=\Op(a_0)$, with Hamiltonian $a_0 \in S^3(\R^2)$ of the form 
\[
a_0(x,\xi)=\xi_1a_2(\xi_1,\xi_2), \quad a_2(\xi_1,\xi_2)=\xi_1^2+\xi_2^2, \quad \xi \in \R^2.
\]
Hence, in this case, the Hamilton flow $(X(s;x,\xi),\Xi(s;x,\xi))$ associated with the Hamilton equations 
\[
\begin{cases}
        \dot \gamma(s)=H_{a_3}(\gamma(s)) \\
        \gamma(0)=(x,\xi).
    \end{cases}
\]
has the simple form 
\[
\begin{split}
&X(s;x,\xi)=(X_1(s;x,\xi),X_2(s;x,\xi))=(x_1+s(3\xi_1^2+\xi_2^2),x_2+s(2\xi_1\xi_2))\\
&\Xi(s;x,\xi)=(\Xi_1(s;x,\xi),\Xi_2(s;x,\xi))=(\xi_1,\xi_2),
\end{split}
\]
for $s \in \R$ and $(x,\xi) \in \R^{2}$. Then, for the standard $ZK$ operator the bicharacteristic flow is uniform straight line motion that can be written as 
\begin{equation}\label{eq.flowZk}
\varphi(s;x,\xi)=(x+s\xi^{(2)},\xi), \quad \text{with} \ \ \xi^{(2)}:=(3\xi_1^2+\xi_2^2,2\xi_1\xi_2),
\end{equation}
The purpose of this work is to study a (multidimensional) variable-coefficients version of the standard ZK-type operator, which in two spatial dimensions is defined as $A = \mathrm{Op}(a_3)$ with symbol
\[
a_3(x,\xi)=\xi_1\bigl(\xi_1^2+a_{22}(x)\xi_2^2\bigr), \quad (x,\xi) \in \mathbb{R}^{4},
\]
where $a_{22}$ satisfies suitable assumptions. In what follows, we study these operators in arbitrary dimensions, as explained below.

\subsection{The general case}
We consider $A=A(x,D)$ of the form
\begin{align}\label{eq.defop}
    A(x,D)=D_1A_2(x,D),\quad A_2(x,D)=\sum_{i,j}D_ia_{ij}(x)D_j,
\end{align}
with $a_{ij} \in C^\infty_b(\R^n;\R)$ for all $i,j=1,\dots,n$.
As before, the operator $A$ should be understood as a generalization of the linear (possibly multidimensional) $ZK$ operator 
\[
A_0=D_{x_1}\Delta,
\]
with $\Delta=\sum_{i=1}^nD_{x_i}^2$ the positive Laplacian, that can be written as $A_0=\Op(a_0)$, with $a_0 \in S^3(\R^n)$ of the form 
\[
a_0(x,\xi)=a_0(\xi)=\xi_1|\xi|^2, \quad \xi \in \R^n.
\]
In this work we assume that the coefficients matrix $(a_{ij}(x))_{i,j}$ which defines 
\begin{equation*}
    a_2(x,\xi):=\sum_{i,j=1}^n a_{ij}(x)\xi_i\xi_j, \quad (x,\xi)\in \R^{2n},
\end{equation*}
and consequently the principal symbol $a_3$ of $A$, defined as 
\begin{equation}\label{eq.defa3}
    a_3(x,\xi):=\xi_1 a_2(x,\xi), \quad (x,\xi) \in \R^{2n},
\end{equation}
satisfies the following properties:
\begin{hp}\label{hp1}(\textit{The block decomposition condition}).
The coefficients matrix $\mathrm{A}(x)=(a_{ij}(x))_{i,j=1,\dots,n}$ that defines $a_2$ is symmetric, positive definite and satisfies
\begin{equation}
a_{1j}=a_{j1}=\delta_{1j}, \quad \forall j=1,\ldots, n. 
\end{equation}
\end{hp}
\begin{hp}\label{hp2}(\textit{The uniform ellipticity condition}). There exists a universal constant $\Lambda>0$ such that
\begin{equation*}
\Lambda^{-1}|\xi|^2\leq a_2(x,\xi)\leq \Lambda |\xi|^2,\quad \forall (x,\xi)\in\R^{2n}.
\end{equation*}
\end{hp}
\begin{hp}\label{hp3}
(\textit{The asymptotically constant condition}). For all $\beta \in \mathbb{N}_0^n$ there exists $C_\beta>0$ such that  
\begin{align}
    |\partial^\beta_x(a_{ij}(x)-\delta_{ij})|\leq \frac{C_\beta}{\langle x\rangle^{|\beta|+2}},\quad \forall x \in \R^n.
\end{align}
\end{hp}
Throughout this work, when $(a_{ij}(x))_{i,j}$ satisfies these conditions, we shall sometimes say that the operator $A$ satisfies the hypotheses \ref{hp1}-\ref{hp3}.

Finally, note that as a consequence of such assumptions we also have that $a_3$ satisfies the following conditions.  
\begin{itemize}
    \item 
There exists a constant $C>0$ such that
    \begin{equation*}
        C^{-1}|\xi|^{2}\leq |\nabla_\xi a_3(x,\xi)|\leq C |\xi|^{2}, \quad \forall (x,\xi) \in \R^{2n}\setminus 0.
    \end{equation*}

   \item For all $\alpha,\beta \in \mathbb{N}_0^n$ there exists a constant $C_{\alpha\beta}>0$ such that 
    \begin{equation}
        |\partial_x^\beta \partial_\xi^\alpha (a_3(x,\xi)-a_0(\xi))| \leq C_{\alpha\beta} \langle x \rangle^{-|\beta|-2}\langle \xi \rangle^{3-|\alpha|}, \quad \forall (x,\xi) \in \R^{2n}.\label{cond.aderivatives}
    \end{equation}
    \end{itemize}

\section{The bicharacteristic flow of ZK-type operators}\label{sec.bicZK}
In this section we study the properties of the bicharacteristic flow of the principal symbol (see \eqref{eq.defa3}) of operators of the form \eqref{eq.defop}.

From now on
$(X(s;x,\xi),\Xi(s;x,\xi))$ will be the bicharacteristic flow of $a_3$. Moreover, we will write $(X(s),\Xi(s))$ when the starting point $(x,\xi)$ is fixed. 
\medskip

\subsection{Basic properties of the flow}
The first fundamental result we establish is the \emph{global} nontrapping property under our conditions \ref{hp1}-\ref{hp3}.

\begin{lemma}\label{lemma.nontrapping}
Let $A$ be a differential operator of order three satisfying \ref{hp1}-\ref{hp3}. Then every point $(x,\xi) \in \R^{2n}\setminus 0$ is both forward and backward nontrapped for $A$.
\end{lemma}
\begin{proof}
    Let $(x,\xi)\in \R^{2n}\setminus 0$. First of all since the principal symbol $a_3$ of $A$ has the form \eqref{eq.defa3},
    we have that by \ref{hp1} the $X_1$-component of the bicharacteristic curve $\gamma=\gamma_{(x,\xi)}$ associated with $a_3$ and starting from $(x,\xi)$, satisfies 
    \[
    \dot X_1(s)=a_2(X(s),\Xi(s))+2\Xi_1(s)^2, \quad s \in I_{(x,\xi)},
    \]
    where $I_{(x,\xi)}\subseteq \R$ is the interval of definition of the bicharacteristic curve, and hence, by \ref{hp2},
    \[
    \dot X_1(s) \approx |\Xi(s)|^2, \quad s \in I_{(x,\xi)}.
    \]
    We next note that $\Xi(s)\neq 0$, for all $s \in I_{(x,\xi)}$. Indeed, assume by contradiction there exists $s_\ast \in I_{(x,\xi)}$ such that $\Xi(s_\ast)=0$ and define $\gamma_\ast$ as the constant curve $\gamma_\ast(s)=(X(s_\ast),0)=:z_\ast$. The curve $\gamma_\ast$ solves the Hamilton equations \eqref{eq.Hameq} starting from $z^\ast$ and is defined for all $s \in \R$. We now set 
    \[
    E=\lbrace s \in I_{(x,\xi)}; \ \gamma(s)=z_\ast \rbrace.
    \]
    Since $s_\ast \in E$ we have that $E$ is not empty. Moreover, is closed due to $E=\gamma^{-1}(\lbrace z^\ast \rbrace)$ and is open, since for all $s_0 \in E$ there exists a sufficiently small neighborhood of $s_0$ which is also contained in $E$ by the (local) uniqueness of the solution of the Cauchy problem \eqref{eq.Hameq} (namely $\gamma_\ast$ is equal to $\gamma$ in such a neighborhood). Then, since $I$ is connected we must have $I=E$ and hence (since $0 \in I=E$) $\gamma(0)=(x,\xi)=z^\ast$ which yields $\xi=0$. This is a contradiction. Then  
    \begin{equation}\label{eq.X_1(s)}
    \dot X_1(s) \approx |\Xi(s)|^2 >0, \quad \forall s \in I.
    \end{equation}
    For what concerns the $\Xi$-component of the flow, by \ref{hp3} we have 
    \begin{equation}\label{eq.aprioriXidot}
    |\dot \Xi(s)|\lesssim \langle X(s) \rangle^{-3}|\Xi(s)|^3.
    \end{equation}
    Since $\Xi(s)\neq 0$, for all $s \in I$, we may consider the function 
    \[
    I \ni s \mapsto \log(|\Xi(s)|), 
    \]
    for which, by \eqref{eq.aprioriXidot},
    \[
    \Bigl|\frac{d}{ds}\log(|\Xi(s)|)\Bigr|= \Bigl|\frac{\Xi(s)\cdot \dot \Xi(s)}{|\Xi(s)|^2}\Bigr|\leq \frac{|\dot \Xi(s)|}{|\Xi(s)|} \lesssim \langle X(s) \rangle^{-3}|\Xi(s)|^2, \quad  \forall s \in I.
    \]
    By \eqref{eq.X_1(s)} we have that this implies that there exists $C>0$ universal constant such that 
    \[
    \Bigl|\frac{d}{ds}\log(|\Xi(s)|)\Bigr|\leq C \frac{\dot X_1(s)}{\langle X_1(s) \rangle^3}.
    \]
    We now take $s_1,s_2 \in I$, with $s_1<s_2$, integrate in the time interval $[s_1,s_2]$ and get 
    \[
    \begin{split}
    |\log(|\Xi(s_2)|)-\log(|\Xi(s_1)|)|&\leq \int_{s_1}^{s_2}\Bigl|\frac{d}{ds}\log(|\Xi(s)|)\Bigr|ds\\
    &\leq C \int_{s_1}^{s_2} \frac{\dot X_1(s)}{\langle X_1(s) \rangle^3}ds\\
    & \leq C \int_{X_1(s_1)}^{X_1(s_2)}\frac{1}{\langle y \rangle^3}dy \\
    &\leq C \int_{\R}\frac{1}{\langle y \rangle^3}dy\\
    &\leq  C',
    \end{split}
    \]
    for some $C'>0$. Therefore, by setting $s_1=0$ and $s_2=s$ we obtain 
    \[
    \Bigl|\log\Bigl(\frac{|\Xi(s)|}{|\xi|}\Bigr)\Bigr|\leq C', \quad \forall s \in I,
    \]
    which yields 
    \[
    M^{-1}|\xi|\leq |\Xi(s)|\leq M |\xi|, \quad \forall s \in I,
    \]
    for some universal $M>0$. This proves that $I_{(x,\xi)}=\R$, due to the Gronwall type argument, since $|\dot \Xi(s)|\lesssim |\xi|$ and $|\dot X(s)|\lesssim |\xi|$. Finally, for all $s>0$,
    \[
    X_1(s)-x_1=\int_0^s \dot X_1(s) ds \gtrsim |\xi|^2s,
    \]
    and so 
    \[
    X_1(s)\geq x_1+\mu |\xi|^2s,
    \]
    for some $\mu>0$. This proves that the point $(x,\xi)$ is forward non trapped. The fact that it is also backward nontrapped follows from 
    \[
    x_1-X_1(s)\gtrsim |\xi|^2|s|, \quad \forall s<0.
    \]
    
\end{proof}

We next prove the following lemma, which provides precise information about the behavior of the bicharacteristic flow associated with the principal symbol of the operators under consideration.

\begin{lemma}\label{st.lem1}
    There is a constant $R_0>0$ such that, for all $R\geq R_0$, if $(x,\xi)\in \R^{2n}\setminus 0$ satisfies \begin{equation}\label{eq.condlemmafond}
        |x|>R,\quad |\xi|=1, \quad 2 \ x \cdot \nabla_\Xi a_3(x,\xi)=2\sum_{j=1}^nx_j\partial_{\Xi_j}a_3(x,\xi)=:N\geq 0,
    \end{equation}
then the associated flow $(X(s),\Xi(s))=(X(s;x,\xi),\Xi(s;x,\xi))$ is defined for $s\geq 0$. Moreover, for all $s \geq 0$ there exist $C_1,C_2>0$, with $C_1<C_2$, such that
\begin{align}\label{lem-0.2}
  |x|^2+ C_1s^2+Ns\leq  |X(s)|^2\leq |x|^2+ C_2s^2+Ns,
\end{align}
\begin{align}\label{lem-1.2}
    2X(s)\cdot\nabla_\Xi a_3(X(s),\Xi(s))\geq N+\frac{C_1}{2}s\geq 0,
\end{align}
\begin{align}\label{lem-3.2}
    \frac 1 C |\xi|\leq |\Xi(s)|\leq C|\xi|,
\end{align}
for some constant $C>0$.
\end{lemma}

\proof
Let $(x,\xi) \in \R^{2n}\setminus 0$ be a point that satisfies the conditions in \eqref{eq.condlemmafond}. The idea of the proof is to estimate the second order derivative of the modulus squared of the $X$-bicharacteristic curve $s \mapsto X(s;x,\xi)$ and then integrate in time.

Since
\begin{align}
    \frac{d}{ds}|X(s)|^2=2 X(s)\cdot \nabla_\Xi a_3(X(s),\Xi(s))
\end{align}
we have
\begin{align}
    \frac{d^2}{ds^2}|X(s)|^2&=2 \frac{d}{ds}X(s)\cdot \nabla_\Xi a_3(X(s),\Xi(s))+ X(s)\cdot \frac{d}{ds}\nabla_\Xi a_3(X(s),\Xi(s))\\
    &=2|\nabla_\Xi a_3(X(s),\Xi(s))|^2+\sum_{j=1}^n X_j(s)\cdot\frac{d}{ds}\partial_{\Xi_j(s)} a_3(X(s),\Xi(s)),
\end{align}
where
\begin{align}
    &\frac{d}{ds}\partial_{\Xi_j} a_3(X(s),\Xi(s))\\
    &=\left(\ds \p_{\Xi_1}a_3(X(s),\Xi(s)),\ldots, \ds \p_{\Xi_n}a_3(X(s),\Xi(s))\right)\\
    &=\left(\nabla_X \p_{\Xi_1}a_3\cdot \nabla_\Xi a_3-\nabla_\Xi \p_{\Xi_1}a_3\cdot \nabla_X a_3, \ldots, \nabla_X \p_{\Xi_n}a_3\cdot \nabla_\Xi a_3-\nabla_\Xi \p_{\Xi_n}a_3\cdot \nabla_X a_3 \right),\label{eq d2x}
\end{align}
in which each component is a symbol of order $4$.

At this point, note that, by the decay conditions, if $|X(s)|\geq R$ then 
\begin{align}
    \left|X_j(s)\frac{d}{ds}\partial_{\Xi_j} a_3(X(s),\Xi(s))\right|\lesssim\langle X(s) \rangle^{-2}|\Xi(s)|^{2m-2}\lesssim \frac{c |\Xi(s)|^{2(m-1)}}{(1+R)^2},
\end{align}
 which leads to (with $c>0$ constant)
\begin{align}
    2|\nabla_\Xi a_3(X(s),\Xi(s))|^2-\frac{c |\Xi(s)|^4}{(1+R)^2}\leq \frac{d^2}{ds^2}|X(s)|^2\leq 2|\nabla_\Xi a_3(X(s),\Xi(s))|^2+\frac{c|\Xi(s)|^{4}}{(1+R)^2}.\label{eq.d2}
\end{align}
We now focus on $2|\nabla_\Xi a_3(X(s),\Xi(s))|^2$ and have (with $\Xi'=(\Xi_2,\dots,\Xi_n)$)
\begin{align}
   2|b_2(X,\Xi)|^2&\leq 2|a_2(X,\Xi)+2\Xi_1^2|^2\\
   &\leq 2|\nabla_\Xi a_3(X,\Xi)|^2\\
   &=2 \left|\Big(a_2(X,\Xi)+2\Xi_1\sum_{j=1}^n a_{1j}\Xi_j, \nabla_{\Xi'} \Xi_1a_2(X,\Xi)\Big)  \right |^2\\
     &=2 \left|\Big(a_2(X,\Xi)+2\Xi_1^2, \nabla_{\Xi'} \Xi_1a_2(x,\Xi)\Big)  \right |^2\\
     &\lesssim 2|a_2(X,\Xi)+2\Xi_1^2|^2\\
     &\lesssim  2|a_2(X,\Xi)|^2. \label{eq.est.nablaa}\\
\end{align}
Moreover, by \ref{hp2}, for all $(X,\Xi) \in \R^{2n}\setminus 0$
\begin{align}
   \Lambda^{-2}|\Xi|^4\leq a_2(X,\Xi)^2
 \leq \Lambda^2|\Xi|^4.
\end{align} 
Then, from \eqref{eq d2x} we get, for some $C_1,C_2>0, C_2>C_1$,
\begin{align}
2C_2|\Xi(s)|^4&\geq 2|\nabla_\Xi a_3(X(s),\Xi(s))|^2+\frac{c|\Xi(s)|^4 }{(1+R)^2} \\
&\geq
   \frac{d^2}{ds^2}|X(s)|^2\\
    &\geq  2|\nabla_\Xi a_3(X(s),\Xi(s))|^2-\frac{c |\Xi(s)|^4}{(1+R)^2}\\
    &\geq \Big(C_1-\frac{c }{(1+R)^2}\Big) |\Xi(s)|^4\\
    &\geq C_1/2|\Xi(s)|^4.\label{eq.dx2}
\end{align}
Hence, in particular 
\begin{equation}
    \frac{d}{ds}|X(s)|^2 \geq N+ C_1/2\int_0^s|\Xi(\sigma)|^4 d\sigma \geq 0,
\end{equation}
which shows that $|X(s)|$ is increasing, so $|X(s)|>R$ for all $s\geq 0$.
Since, by assumption $|\xi|=|\Xi(0)|=1$, we have $|\Xi(s)|\geq 1/2$ for $s$ sufficiently small, thus
\begin{equation}\label{eq.defsast}
s^*:=\sup\{s: |X(\tau)|>R \,\, \text{and}\,\, \frac{1}{2} |\Xi(0)|^4\leq |\Xi(\tau)|^4\leq 2|\Xi(0)|^4 ,\quad 0<\tau<s\}>0.
\end{equation}
Now, for all $s\in [0,s^*)$, if $R$ is sufficiently large
    \begin{align}
   2C_2=2C_2|\Xi(0)|^4\geq  \frac{d^2}{ds^2}|X(s)|^2\geq C_1/2|\Xi(0)|^4=C_1/2. \label{eq:1811-1}
\end{align}
We then integrate \eqref{eq:1811-1} in time and get that for all $s \in [0,s^\ast)$
\begin{align}
   N+2C_2s&\geq \frac{d}{ds}|X(s)|^2\\
   &\geq \frac{d}{ds}|X(s)|^2|_{s=0}+\frac{C_1}{2}s=2X(0)\cdot\nabla_{\Xi}a_3(X(s),\xi(s))|_{s=0}+\frac{C_1}{2}s|\Xi(0)|^4=N+\frac{C_1}{2}s> 0,\label{eq.dx}
\end{align}
which gives \eqref{lem-1.2} in $[0,s^\ast)$, that is
\begin{align}
   2X(s)\cdot\nabla_\Xi a_3(X(s),\Xi(s))\geq N+\frac{C_1}{2}s.\label{eq.1811-2}
\end{align}
An integration in time of \eqref{eq.1811-2} now yields
\begin{align}
 |X(0)|^2+Ns+C_2s^2   \geq |X(s)|^2\geq |X(0)|^2+Ns+\frac{C_1s^2}{4},\label{eq.1811-3}
\end{align}
which proves the lower bound \eqref{lem-0.2} in $[0,s^*)$. Note that the validity of \eqref{lem-3.2} in $[0,s^\ast)$ follows from the definition of $s^\ast$.

We finally prove that $s^\ast=+\infty$. Let us assume by contradiction $s^\ast<+\infty$.
Since $|X(s)|$ is increasing $\in [0,s^\ast)$ we have
\begin{align}
    \frac{d}{ds}|\Xi(s)|^2\leq 2|\sum_k\Xi_k(s) \frac{d}{ds}\Xi_k(s)|&=2|\sum_k \Xi_k(s) \partial_X a(X(s),\Xi(s)) |\lesssim \frac{1}{\langle X(s)\rangle^{3}}|\Xi(s)|^4\\
    & \lesssim \frac{1}{(X(0)^2+s^2)^{3/2}})|\Xi(s)|^4\\
     & \lesssim \frac{1}{(R+s|\Xi(0)|)^3}|\Xi(s)|^4,
\end{align}
and then
\begin{align}
  \Big|  |\Xi(s)|^2-|\Xi(0)|^2\Big|=\Big|\int_{0}^s \frac{d}{ds'}|\Xi(s')|^2 ds'\Big|&\lesssim |\Xi(0)|^3 \int_0^s\frac{|\Xi(0)|}{(R+|\Xi(0)|s')^3}ds'\\
  &\lesssim |\Xi(0)|^3\Big( \frac{1}{R^2}-\frac{1}{(R+|\Xi(0)|s)^2}\Big)\lesssim \frac{|\Xi(0)|^3}{R^2}.
\end{align}
Since $|\Xi(0)|^2=|\Xi(0)|^3=1$, we get
\begin{align}
 (1-\frac{C}{R^2})=(1-\frac{C}{R^2})|\Xi(0)|^2 \leq|\Xi(s)|^2\leq(1+\frac{C}{R^2})|\Xi(0)|^2=(1+\frac{C}{R^2}).
\end{align}
Therefore, for $R$ sufficiently large,
$$ \frac 3 4 |\Xi(0)|^4< |\Xi(s)|^4<\frac 5 4 |\Xi(0)|^4.$$
Hence, by continuity (in time) of the flow $(X(s),\Xi(s))$ we have that there exists $\delta>0$ sufficiently small such that for all $s \in [0,s^*+\delta)$
\[
|X(s)|>R \quad \text{and} \quad \frac 1 2 |\Xi(0)|^4< |\Xi(s)|^4<2 |\Xi(0)|^4.
\]
This contradicts the definition of $s^\ast$ as the supremum in \eqref{eq.defsast}. Therefore, \eqref{lem-0.2}, \eqref{lem-1.2}, and \eqref{lem-3.2} hold for all $s \geq 0$, and the proof is concluded. 
\endproof
\begin{remark}\label{rem.s<0}
    Note that Lemma \ref{st.lem1} holds for $s\leq 0$ if the starting point $(x,\xi)$ satisfies
    \begin{equation}\label{eq.condlemmafond-s<0}
        |x|>R,\quad |\xi|=1, \quad  2\, x \cdot \nabla_\Xi a_3(x,\xi)=2\sum_{j=1}^nx_j\partial_{\Xi_j}a_3(x,\xi)=:-N\leq 0,
    \end{equation}
    for some $N>0$. 
Namely, in this case the flow $(X(s;x,\xi),\Xi(s;x,\xi))$ is defined for $s\leq 0$,  and \eqref{lem-0.2} and \eqref{lem-1.2} assume the form
\begin{align}
     |x|^2+ C_1s^2-Ns\leq  |X(s)|^2\leq |x|^2+ C_2s^2-Ns
\end{align}
    \begin{align}\label{lem-1.2-s<0}
    2X(s)\cdot\nabla_\Xi a_3(X(s),\Xi(s))\leq -N+2C_1s\leq 0,
\end{align}
where $C_2,C_1$ are universal positive constant such that $C_2>C_1$.

    The proof follows mutatis mutandi from that of Lemma \ref{st.lem1} by simply noticing that we will have $-N$ in place of $N$ in the whole proof, and that $-Ns\geq 0$ since $s\leq 0$.
    
\end{remark}

Some consequences of the previous lemma are in order.
\begin{itemize}
 
    \item If $(x,\xi)$ is such that $|x|> R_0$, $|\xi|=1$,  and $2\sum_{j}x_j\partial_{\xi_j}a_3(x,\xi)=N,\,\, N\geq0,$ then there exists a constant $M>0$ such that
    \begin{equation}\label{eq.xsscatt}
        \frac 1 M (|x|+s)\leq |X(s)|\leq M(|x|+s), \quad \forall 0<s<s^*,
    \end{equation}
    that is, in other words, 
    $$\langle X(s)\rangle\approx \langle |x|+s\rangle,$$ so $X(s)$ grows linearly in time for $s>0$. Note that this is a useful feature of flows of constant coefficient operators.
    \item It is possible to conclude the same property obtained for outgoing orbits ($s>0$) also for incoming orbits starting at points $(x,\xi)\in B_{R_0}^c\times\mathbb{S}^{n-1}$ such that $2\sum_{j}x_j\partial_{\xi_j}a_3(x,\xi)=-N,\,\, N\geq0$. In this case we consider
    \begin{align}
        s_*=\inf\{s: |X(\tau)|\geq R_0\,\, \text{and}\,\, \frac 1 2 \leq |\Xi(\tau)|^4\leq 2,\,\, \, s<\tau<0 \},
    \end{align}
    which, as in the outcoming case, can be proven to be $-\infty$.
    \item Repeating Lemma \ref{st.lem1} for incoming orbits (see Remark \ref{rem.s<0}, we obtain that there exists $M>0$ such that
     \begin{equation}\label{eq.xsscatt-s<0}
        \frac{1}{M}(|x|-s)\leq |X(s)|\leq M(|x|-s), \quad \forall s_*<s<0, \quad |\xi|=1.
    \end{equation} 
\end{itemize}

\subsection{The scattering regions}
Following \cite{CKS} we define the \textit{scattering (or outgoing) regions} of $\R^{2n}$, for forward and backward time, as follows. 
\begin{definition}\label{def.scatteringsets} Let $R_0>0$ be as in Lemma \ref{st.lem1}.
    \textit{The scattering region of $\R^{2n}$ for forward time and backward time} of $a_3$ are defined, respectively, by
\begin{align}
    \mathscr{S}_+&:=\{(x,\xi)\in \R^{2n}; \ \ |x|\geq R_0, \ \sum_{j=1}^n x_j\partial_{\Xi_j}a_3(x,\xi)\geq 0 \},\\
    \mathscr{S}_-&:=\{(x,\xi)\in \R^{2n}; \ \ |x|\geq R_0, \ \sum_{j=1}^n x_j\partial_{\Xi_j}a_3(x,\xi)\leq 0 \}.
\end{align}
We also denote by $\mathscr{S}_{+,1}:= B_{R_0}^c\times \mathbb{S}^{n-1}\cap\mathscr{S}_+$ and by $\mathscr{S}_{-,1}:=B_{R_0}^c\times \mathbb{S}^{n-1}\cap\mathscr{S}_+$.
\end{definition}
The motivation behind these definitions is the following.
A point $(x,\xi)\in\mathscr{S}_+$ (resp.\ $\mathscr{S}_-$)
satisfies $|x|\geq R_0$ and
\[
\left.\frac{d}{ds}|X(s;x,\xi)|^2\right|_{s=0}
=2x\cdot\nabla_\xi a_3(x,\xi)
\geq 0
\qquad (\text{resp.}\ \leq 0).
\]
The latter condition means that the spatial component of the
Hamiltonian flow has a nonnegative (resp.\ nonpositive)
radial velocity at $s=0$, namely that such points are instantaneously
directed away from  the origin (resp. toward the origin). Thus, $\mathscr{S}_+$ and $\mathscr{S}_-$ represent the
instantaneously outgoing and incoming regions of phase space,
respectively, for points in the configuration space lying outside
a large ball centered at the origin.

\begin{remark}
    Note that we are free to choose $R_0$ in Definition \ref{def.scatteringsets} larger than the optimal one given in Lemma \ref{st.lem1}. This will not alter the content of the forthcoming results, as scattering points are used just as a bridge to obtain results for nontrapped points. Thus, restrictions on the choice of $R_0$ are not needed. This avoids the unnecessary and unpleasant use of subsets of the "optimal" $\mathscr{S}_-$. In view of these observations, we will take $R_0$ large enough in the coming statements.
\end{remark}

We will now show that, as in the case treated in \cite{CKS}, the flow starting at scattering points is asymptotically linear far from the ball $B_{R_0}$ ($R_0$ as in Lemma \ref{st.lem1} or larger). To prove this result, we will need the following preliminary lemma on the homogeneity properties of the flow of homogeneous operators of order three (cf. Lemma 2.3 in \cite{CKS}).

\begin{lemma}\label{lem.hom.prop}
    The principal symbol  $a_3(X,\Xi)$ of $A$ is homogeneous of degree $3$, i.e.
    \begin{equation}\label{eq.homa3}
        a_3(X,r\Xi)=r^3a_3(X,\Xi), \quad (X,\Xi) \in \R^{2n}\setminus 0, \ r>0.
    \end{equation}
    Moreover, the associated flow $(X(s;x,\xi),\Xi(s
    ;x,\xi))$ satisfies 
    \begin{equation}\label{eq.homflow}
    \begin{split}
    X(s;x,r\xi)&=X(r^2s;x,\xi)\\
        \Xi(s;x,r\xi)&=r\Xi(r^2s;x,\xi),
        \end{split}
    \end{equation}
    and, for all multiindices $\alpha,\beta \in \mathbb{N}_0^n$,
    \begin{equation}\label{eq.homflowder}
        \begin{split}(\partial^\beta_x\partial^\alpha_\xi X)(s;x,r\xi)&=r^{-|\alpha|}\partial^\beta_x\partial^\alpha_\xi X(r^2s;x,\xi)\\
        (\partial^\beta_x\partial^\alpha_\xi \Xi)(s;x,r\xi)&=r^{1-|\alpha|}\partial^\beta_x\partial^\alpha_\xi \Xi(r^2s;x,\xi),
        \end{split}
    \end{equation}
    for all $(x,\xi) \in \R^{2n}\setminus 0$, $r>0$ and $s \in  \R$ such that the integral curves are defined.
\end{lemma}
\begin{proof}
    The property \eqref{eq.homa3} follows directly by the definition of $a_3$ (recall \eqref{eq.defa3}). 
    
    \noindent As for the homogeneous property of the flow \eqref{eq.homflow}, it suffices to note that, if $r>0$, then
    $(\tilde{X}(s;x,\xi),\tilde{\Xi}(s;x,\xi)):=(X(r^2s;x,\xi),r\Xi(r^2s;x,\xi))$ solves the ODE
    \begin{equation}
        \begin{cases}
            \dot \gamma(s)=H_{a_3}(\gamma(s)), \\
            \gamma(0)=(x,r\xi).
        \end{cases}
    \end{equation}
    Finally, for $r>0$, denoting $\tilde{\xi}=r\xi$, we have
    \begin{equation}
    \begin{split}
        &\frac{1}{r}\frac{\partial}{\partial \xi} X(r^2s;x,\xi)=\frac{\partial}{\partial \tilde{\xi}}X(s;x,\tilde\xi) \\
        &\frac{\partial}{\partial \xi}\Xi(r^2s;x,\xi)=\frac{\partial}{\partial \tilde{\xi}}\Xi(s;x,\tilde{\xi}),
        \end{split}
    \end{equation}
    Iterating this argument for higher order derivatives and taking into account derivatives with respect to the $x$-variables, we obtain \eqref{eq.homflowder}.
\end{proof}

\begin{remark}\label{rmk.behavscattpoints}
    Note that by \eqref{eq.xsscatt} and \eqref{eq.homflow}, if $(x,\xi) \in \mathscr{S}_{+}$, we have that $X(s)=X(s;x,\xi)$ satisfies 
    \begin{align}\label{asypt.xflow}
         \langle X(s)\rangle  \approx \langle x \rangle +s|\xi|^2,\quad s\geq 0,
    \end{align}
    and, moreover, by \eqref{lem-3.2} and \eqref{eq.homflow}, $\Xi(s)=\Xi(s;x,\xi)$ fulfills
    \begin{align} \label{asypt.xiflow}
        |\Xi(s)|\lesssim |\xi|.
    \end{align}
    Similarly, if $(x,\xi)\in \mathscr{S}_-$, by \eqref{eq.xsscatt-s<0} and \eqref{eq.homflow}, we have
     \begin{align}\label{asypt.xflow-} 
         \langle X(s)\rangle  \approx \langle x \rangle -s|\xi|^2,\quad s\leq 0.
    \end{align}
\end{remark}

We are now ready to describe the asymptotic behavior of the flow as it moves away from the ball of radius $R_0$.

We are going to state and prove the results for $(x,\xi)\in\mathscr{S}_+$ when $s\geq 0$. The analogous statements hold for $(x,\xi)\in\mathscr{S}_-$ when $s\leq 0$.
\begin{lemma}\label{lem.scatteringV}
    Let $(x,\xi)\in \mathscr{S}_+$. Then there exist (forward) \textit{scattering variables} $(\widebar{X}(x,\xi),\widebar{\Xi}(x,\xi))$ such that, for all $s\geq 0$,
    \begin{equation}\label{st.lem2-eq1}
        |\Xi(s;x,\xi)-\widebar{\Xi}(x,\xi)|\lesssim \frac{|\xi|}{\langle X(s)\rangle^2},
    \end{equation}
    \begin{equation}\label{st.lem2-eq2}
        \Bigl|X(s;x,\xi)-\big(\widebar{X}(x,\xi)+s\widebar{\Xi}^{(2)}(x,\xi)\big)\Bigr|\lesssim \frac{1}{\langle X(s)\rangle},
    \end{equation}
    where 
    \[
    \widebar{\Xi}^{(2)}(x,\xi):=\big(\sum_{j=1}^nc_j\widebar{\Xi}^2_j(x,\xi), 2\widebar{\Xi}_1(x,\xi)\widebar{\Xi}_2(x,\xi),\ldots,2\widebar{\Xi}_1(x,\xi)\widebar{\Xi}_n(x,\xi)\big),
    \]
    for $c_1=3$ and $c_j=1$ for all $j=2,\ldots,n$.
    
    \noindent In particular, \begin{equation}\label{eq.scattbehav}
    \widebar{\Xi}(x,\xi)=\xi+O(|\xi|\langle x\rangle^{-2}), \quad  \widebar{X}(x,\xi)=x+O(\langle x\rangle^{-1}).
    \end{equation}
\end{lemma}

\proof
Below we give the proof when $(x,\xi)\in \mathscr{S}_{+,1}$, since,
 once the result is proved in this case, then  the general case $(x,\xi) \in \mathscr{S}_+$ follows from the homogeneity properties established in Lemma \ref{lem.hom.prop}. Hence, let $(x,\xi) \in \mathscr{S}_{+,1}$ and denote $(X(s),\Xi(s))=(X(s;x,\xi),\Xi(s;x,\xi))$. By \ref{hp3}, \eqref{eq.xsscatt} and Remark \ref{rmk.behavscattpoints},  for all $s\geq 0$
\[
|\partial_{x_j}a(X(s),\Xi(s))|\lesssim (|x|+s|\xi|^2)^{-3}.
\]
Define now
\begin{equation}\label{eq.defxibar}
    \widebar{\Xi}(x,\xi):=\xi+\int_0^{+\infty} \frac{d}{dt}\Xi(t)dt,
\end{equation}
then 
\begin{equation}\label{eq.scattxi}
\begin{split}
    |\Xi(s;x,\xi)-\widebar{\Xi}(x,\xi)|&=\left|-\int_s^{+\infty}\frac{d}{dt}\Xi(t)\,dt\right|\\
    &\leq\int_s^{+\infty}|\nabla_X a_3(X(t),\Xi(t))|\,dt\\
    &\lesssim\int_s^{+\infty}\frac{|\xi|^3}{(|x|+t|\xi|^2)^3}\,dt\\
    &\lesssim |\xi|(|x|+s|\xi|^2)^{-2}
    \lesssim\langle X(s)\rangle^{-2},
\end{split}
\end{equation}
which proves \eqref{st.lem2-eq1}, when $(x,\xi) \in \mathscr{S}_{+,1}$.

Next, we prove \eqref{st.lem2-eq2}. Recall that
$$\nabla_\Xi a_3(X,\Xi)=\Big(b_2(X,\Xi)+2\Xi_1^2, \nabla_{\Xi'} \Xi_1b_2(x,\Xi)\Big), $$
where $$b_2(X,\Xi)=\sum_{ij}b_{ij}\Xi_i\Xi_j$$
is asymptotically constant (see \ref{hp3}). Note also that
\begin{align}
    \ds \left(X_j(s;x,\xi)-s\widebar{\Xi}^2_j(x,\xi)\right)= \partial_{\Xi_j}a_3(X(s),\Xi(s))-\widebar{\Xi}^2_j(x,\xi).
\end{align}
Then, due to $|\Xi(s)|\leq C$ and \eqref{st.lem2-eq1}, we have (recall $c_j=1$ for $j\neq 1$ and $c_1=3$)
\begin{align}
  &\left| \ds \left(X_1(s;x,\xi)-s\sum_{j=1}^nc_j\widebar{\Xi}^2_j(x,\xi)\right)\right| \\
  &=|a_2(X(s),\Xi(s))+2\Xi_1(s)^2-\sum_{j=1}^nc_j\widebar{\Xi}^2_j(x,\xi)| \label{eq.Flow_1}\\
  &=|\sum_{i,j=2}^n(a_{ij}(X(s))-\delta_{ij}))\Xi_i(s)\Xi_j(s)+\sum_{j=2}^n(\Xi^{(2)}_j(s)-\widebar{\Xi}_j^2(x,\xi))+3(\Xi^{(2)}_1(s)-\widebar{\Xi}_1^2(x,\xi))|\\
  &\lesssim \frac{1}{\la X(s)\ra^2},
\end{align}
where the last inequality follows from \ref{hp3} and 
\begin{align}
    |\Xi^{2}_j(s)-\widebar{\Xi}_j^2(x,\xi)|&\lesssim|\Xi_j(s)+\widebar{\Xi}_j(x,\xi)|\la X(s)\ra^{-2}\\
    &\lesssim \la X(s)\ra^{-2}(2|\Xi(s)|+\la X(s)\ra^{-2})\\
    &\lesssim\la X(s)\ra^{-2}.
\end{align}

On the other hand, for all $j\neq 1$,
\begin{align}
    &\left| \ds \left(X_j(s;x,\xi)-2s \widebar{\Xi}_1(x,\xi)\widebar{\Xi}_j(x,\xi)\right)\right| \\
    &=|\Xi_1(s)\partial_{\Xi_j}a_2(X(s),\Xi(s))- 2\widebar{\Xi}_1(x,\xi)\widebar{\Xi}_j(x,\xi)|\\
    &=\left| 2\sum_{i=2}^n (a_{ij}(X(s))-\delta_{ij})\Xi_1(s)\Xi_i(s)+2(\Xi_1(s)\Xi_j(s)-\widebar{\Xi}_1(x,\xi)\widebar{\Xi}_j(x,\xi))\right|\label{eq.dsYj}\\
    &\lesssim \frac{1}{\la X(s)\ra^2},
\end{align}
where, in this case, in the last inequality we have used that
\begin{align}
\Xi_1\Xi_j-\widebar{\Xi}_1\widebar{\Xi}_j=(\Xi_1-\widebar{\Xi}_1)\Xi_j\cancel{+\widebar{\Xi}_1 \Xi_j}-\widebar{\Xi}_1(\widebar{\Xi}_j-\Xi_j)\cancel{-\widebar{\Xi}_1\Xi_j}\lesssim\frac{1}{\la X(s)\ra^2}.
\end{align}
In turn, the previous estimates yield
\begin{align}
    \left| \ds (X(s;x,\xi)-s \big(\sum_{j=1}^nc_j\widebar{\Xi}^2_j(x,\xi), 2\widebar{\Xi}_1(x,\xi)\widebar{\Xi}_2(x,\xi),\ldots,2\widebar{\Xi}_1(x,\xi)\widebar{\Xi}_n(x,\xi)\big) )\right |\lesssim \frac{1}{\la X(s)\ra^2}.
\end{align}
 Finally, we define

$$ Y(s):= X(s;x,\xi)-s\big(\, \sum_{j=1}^nc_j\widebar{\Xi}^2_j(x,\xi), 2\widebar{\Xi}_1(x,\xi)\widebar{\Xi}_2(x,\xi),\ldots,2\widebar{\Xi}_1(x,\xi)\widebar{\Xi}_n(x,\xi)\,\big) $$
and 
\begin{equation}\label{eq.defxbar}
    \widebar{X}(x,\xi):=x+\int_0^{+\infty} \ds Y(s) ds,
\end{equation}
so that 
\begin{equation}\label{eq.scattx}
    |Y(s)-\widebar{X}(x,\xi)|=\left| \int_s^{+\infty}\dt Y(t)\,\,dt\right|\lesssim\int_s^{+\infty} ( 1+|x|+t)^{-2}dt\lesssim \la|x|+s\ra^{-1},
\end{equation}
which concludes the proof of the statement in the case $(x,\xi) \in \mathscr{S}_{+,1}$. The proof for an arbirtary point in $\mathscr{S}_+$, as mentioned at the beginning of the proof, follows by the homogeneity of the flow (see Lemma \ref{lem.hom.prop}). Finally \eqref{eq.scattbehav} follows by definition \eqref{eq.defxibar} and \eqref{eq.defxbar} by applying the estimates \eqref{eq.scattxi} and \eqref{eq.scattx} at time $s=0$.
\endproof
\begin{remark}\label{rmk.Xi2}
We emphasize that the key information provided by Lemma \ref{lem.scatteringV} is that the asymptotic behavior, that is as $s \to +\infty$,  of the bicharacteristic flow starting at a point in the scattering region $\mathscr{S}_{+}$, is the same as that of the flow of the constant-coefficient (multidimensional) standard $ZK$ operator defined in \eqref{eq.flowZk}. 
\end{remark}

\begin{remark}\label{rmk.deflemmabackwardsv}
The analogous backward scattering variables are defined on $\mathscr{S}_-$ by
\[
\widebar\Xi_{-}(x,\xi)=\xi-\int_{-\infty}^0\dot\Xi(s;x,\xi)\,ds,
\qquad
\widebar X_{-}(x,\xi)=x-\int_{-\infty}^0\dot Y_{-}(s;x,\xi)\,ds,
\]
where $Y_{-}(s)=X(s)-s\widebar\Xi_{-}^{(2)}$.
They satisfy the same estimates for $s\leq0$. Below, the scattering variables without a subscript refer to the forward construction.
\end{remark}

A fundamental consequence of the existence of scattering regions is that, on denoting by $\varphi(s;x,\xi)$ the bicharacteristic starting at $(x,\xi)$, then
\begin{align}
    \varphi(s;\varphi(\tau;x,\xi))=\varphi(s+\tau;x,\xi).
\end{align}
Then, using the notation of Lemma \ref{lem.scatteringV}, if
$$\varphi(s;x,\xi)\simeq(\widebar{X}(x,\xi)+s\widebar{\Xi}^{(2)}(x,\xi),\widebar{\Xi}(x,\xi)),\quad \text{as}\,\, s\rightarrow +\infty,$$
we find that another point in the orbit, say $\varphi(\tau;x,\xi)$, has scattering variables
\begin{align}\label{pro.scatt.var}
(\widebar{X}(\varphi(\tau;x,\xi),\widebar{\Xi}(\varphi(\tau;x,\xi))   = (\widebar{X}(x,\xi)+\tau\widebar{\Xi}^{(2)}(x,\xi),\widebar{\Xi}(x,\xi)).
\end{align}
\begin{remark}\label{rmk.homscattvar}
    We also have that, for all $j=1,\dots, n$, $\widebar{\Xi}_j$ is homogeneous of degree $1$ and $\widebar{X}_j$ is homogeneous of degree $0$. Indeed, by Lemma \ref{lem.hom.prop}, for $\lambda>0$ we have (recall \eqref{eq.defxibar})
    \[
    \begin{split}
    \widebar{\Xi}_j(x,\lambda \xi)&= \lambda \xi_j+\int_0^{+\infty}\dot \Xi_j(s;x,\lambda\xi)ds \\
    &=\lambda \xi_j-\int_0^{+\infty}\partial_{X_j}a(X(s;x,\lambda \xi),\Xi(s;x,\lambda \xi))ds \\
    &=\lambda \xi_j-\int_0^{+\infty}\partial_{X_j}a(X(\lambda^2s;x, \xi),\lambda\Xi(\lambda^2s;x, \xi))ds \\
    &=\lambda \xi_j-\lambda^3\int_0^{+\infty}\partial_{X_j}a(X(\lambda^2s;x,\xi),\Xi(\lambda^2s;x, \xi))ds \\
    &= \lambda \xi_j-\lambda\int_0^{+\infty}\partial_{X_j}a(X(\tau;x, \xi),\Xi(\tau;x, \xi))d\tau \\
    &=\lambda \widebar{\Xi}_j(x,\xi),
    \end{split}
    \]
    and similarly (recall \eqref{eq.defxbar})
    \[
    \begin{split}
    \widebar{X}_j(x,\lambda\xi)&=x_j+\int_0^{+\infty} \dot Y_j(s;x,\lambda\xi)ds \\
    &=x_j+\int_0^{+\infty} \dot Y_j(\tau;x,\xi)d\tau \\
    &=\widebar{X}_j(x,\xi).
    \end{split}
    \]
\end{remark}

\begin{lemma}[Estimates of derivatives of the flow]\label{lem.estderflow}
    Let $(x,\xi)\in\mathscr{S}_+$ and $0\leq s<+\infty$. Then, there exists $C_1>0$ such that, for all $k=1,\dots,n$,
    \begin{align}
        |\p_{x_k} X(s;x,\xi)|&\leq C_1\left( 1+\frac{|s||\xi|^2}{\la x\ra^2}\right)  \label{eq.estderflow1}\\  
        |\p_{x_k}\Xi(s;x,\xi)|&\leq C_1\frac{|\xi|}{\la x\ra^2}  \label{eq.estderflow2},    \end{align}
    and 
     \begin{align}
        |\p_{\xi_k} X(s;x,\xi)|&\leq C_1|s||\xi| \label{eq.estderflow3} \\  
        |\p_{\xi_k}\Xi(s;x,\xi)|&\leq C_1. \label{eq.estderflow4} \\  
    \end{align}
    Moreover, if $\alpha,\beta \in \mathbb{N}_0^n$ with $|\alpha|+|\beta|>1$, we have that, if $|\beta|=0$ and $|\alpha|=2$, 
    \begin{align}
        |\p_x^\beta\p_\xi^\alpha X(s;x,\xi)|\leq C_{\alpha,\beta} \, |s|, \\
        |\p_x^\beta\p_\xi^\alpha \Xi(s;x,\xi)|\leq C_{\alpha,\beta}\frac{|\xi|^{1-|\alpha|}}{\la x\ra^{k(\alpha,\beta)}},
    \end{align}
    and otherwise 
    \begin{align}
        |\p_x^\beta\p_\xi^\alpha X(s;x,\xi)|&\leq C_{\alpha,\beta}\frac{|s||\xi|^{2-|\alpha|}}{\la x\ra^{k(\alpha,\beta)}},\\  
        |\p_x^\beta\p_\xi^\alpha \Xi(s;x,\xi)|&\leq C_{\alpha,\beta}\frac{|\xi|^{1-|\alpha|}}{\la x\ra^{k(\alpha,\beta)}},
    \end{align}
    where $C_{\alpha,\beta}>0$ and $k(\alpha,\beta)>|\beta|$.
\end{lemma}

\proof
Recall that
$$\nabla_\Xi a_3(X,\Xi)=\Big(a_2(X,\Xi)+2\Xi_1^2, \nabla_{\Xi'} \Xi_1a_2(x,\Xi)\Big) $$
where $$a_2(X,\Xi)=\sum_{i,j=1}^na_{ij}(X)\Xi_i\Xi_j=\Xi_1^2+\sum_{i,j=2}^na_{ij}(X)\Xi_i\Xi_j$$
is asymptotically constant. 

Then, denoting $(X,\Xi)=(X(s;x,\xi),\Xi(s;x,\xi))$, we have
\begin{align}
    \ds \p_{x_k}X_1&=3\p_{x_k}(\Xi_1^2)+\sum_{i,j=2}^n\nabla_Xa_{ij}(X)\cdot (\p_{x_k}X) \Xi_i\Xi_j+\sum_{i,j=2}^na_{ij}(X)\p_{x_k}(\Xi_j\Xi_i)\\
    &=2\p_{x_k}(\Xi_1^2)+\sum_{i,j=2}^n\nabla_Xa_{ij}(X)\cdot (\p_{x_k}X) \Xi_i\Xi_j+\sum_{i,j=1}^na_{ij}(X)\p_{x_k}(\Xi_j\Xi_i)\\
    &=2\p_{x_k}(\Xi_1^2)+\sum_{i,j=2}^n\nabla_Xa_{ij}(X)\cdot (\p_{x_k}X) \Xi_i\Xi_j+\sum_{i,j=1}^n(a_{ij}-\delta_{ij})\p_{x_k}(\Xi_j\Xi_i)+\sum_{j=1}^n\p_{x_k}\Xi_j^2,
\end{align}
and, for $\ell\neq 1$,
\begin{align}
    \ds \p_{x_k}X_\ell&=2\sum_{i=2}^n\nabla_Xa_{i\ell}\cdot (\p_{x_k}X) (\Xi_i\Xi_1)+2\sum_{i=2}^na_{i\ell}\p_{x_k}(\Xi_i\Xi_1)\\
    &=2\sum_{i=2}^n\nabla_Xa_{i\ell}\cdot (\p_{x_k}X) (\Xi_i\Xi_1)+2\sum_{i=2}^n(a_{i\ell}-\delta_{i\ell})\p_{x_k}(\Xi_i\Xi_1)+2\p_{x_k}(\Xi_1\Xi_\ell).
\end{align}
Moreover, for all $i,j=1,\dots,n$,
\begin{align}
    &\ds \p_{x_k}(\Xi_i\Xi_j)=\\
    &=-\p_{x_k}\Big(\Xi_j\Xi_1\sum_{\ell,m=1}^n (\p_{X_i}a_{\ell m})(X)\Xi_m\Xi_\ell+\Xi_i\Xi_1\sum_{\ell,m=1}^n (\p_{X_j}a_{\ell m})(X)\Xi_m\Xi_\ell\Big)\\
    &=-\sum_{\ell,m=1}^n \Big( (\p_{x_k}\Xi_1\Xi_j) (\p_{X_i}a_{\ell m})(X)\Xi_m\Xi_\ell+ \Xi_1\Xi_j \nabla_X (\p_{X_i}a_{\ell m})(X)\cdot (\p_{x_k}X)\Xi_m\Xi_\ell\\
    &+ \Xi_1\Xi_j (\p_{X_i}a_{\ell m})(X)(\p_{x_k}\Xi_m\Xi_\ell)  \Big)\\
    &-\sum_{\ell,m=1}^n \Big( (\p_{x_k}\Xi_1\Xi_i) (\p_{X_j}a_{\ell m})(X)\Xi_m\Xi_\ell+ \Xi_1\Xi_i \nabla_X (\p_{X_j}a_{\ell m})(X)\cdot (\p_{x_k}X)\Xi_m\Xi_\ell\\
    &+ \Xi_1\Xi_i (\p_{X_j}a_{\ell m})(X)(\p_{x_k}\Xi_m\Xi_\ell) 
    \Big)\\
    &=-\sum_{\ell,m=1}^n \Big( (\p_{x_k}\Xi_1\Xi_j) (\p_{X_i}(a_{\ell m}-\delta_{\ell m}))(X)\Xi_m\Xi_\ell+ \Xi_1\Xi_j \nabla_X (\p_{X_i}(a_{\ell m}-\delta_{\ell m}))(X)\cdot (\p_{x_k}X)\Xi_m\Xi_\ell\\
    &+ \Xi_1\Xi_j (\p_{X_i}(a_{\ell m}-\delta_{\ell m}))(X)(\p_{x_k}\Xi_m\Xi_\ell) \Big)\\
    &-\sum_{\ell,m=1}^n \Big( (\p_{x_k}\Xi_1\Xi_i) (\p_{X_j}(a_{\ell m}-\delta_{\ell m}))(X)\Xi_m\Xi_\ell+ \Xi_1\Xi_i \nabla_X (\p_{X_j}(a_{\ell m}-\delta_{\ell m}))(X)\cdot (\p_{x_k}X)\Xi_m\Xi_\ell\\
    &+ \Xi_1\Xi_i (\p_{X_j}(a_{\ell m}-\delta_{\ell m}))(X)(\p_{x_k}\Xi_m\Xi_\ell) \Big).
\end{align}
Let now
\begin{align}
    {\bf \Xi}^{(2)}:=(\Xi_1^2,\dots,\Xi_1\Xi_n, \Xi_2^2, \Xi_2\Xi_3,\ldots , \Xi_2\Xi_n,\ldots, \Xi_n^2)^T\in \R^{(n+1)n/2}.
\end{align}
Then, on denoting by $\mathfrak{I}_{1\times (n(n-1)/2)}=(1, \dots, 0, \dots, 1, 0,\dots, 0, \dots, 1) \in \R^{n(n-1)/2}$ the vector whose entries are all 1 corresponding to the position $\Xi_j^2$ for $j=2,\dots, n$ and all $0$ corresponding to the position $\Xi_i\Xi_j$ for $i,j\geq 2$, $i \neq j$. 


\begin{equation}\label{eq.dsdx}
\begin{split}
    &\ds \left(\begin{array}{c}
         \p_{x_k}X \\
         \p_{x_k}{\bf \Xi}^{(2)}
    \end{array} \right)=\\
    &=\left(\begin{array}{cccc}
    0_{1\times n} & 3 & 0_{1\times (n-1)} & \mathfrak{I}_{1\times n(n-1)/2}\\
        0_{n-1\times n} & 0_{n-1\times 1}& 2I_{n-1} & 0_{n-1\times n(n-1)/2} \\
        0_n & 0_{n\times 1} & 0_{n\times n-1}& 0_{n\times n(n-1)/2 }\\
        0_{n-1\times n} & 0_{n-1\times1} & 0_{n-1\times n-1} & 0_{n-1\times (n(n-1)/2)}\\
        0_{n-2\times n} & 0_{n-2\times 1}& 0_{n-2\times n-1} & 0_{n-2\times (n(n-1)/2)}\\
        \vdots & \vdots &\vdots & \vdots\\
        0_{1\times n} & 0 & 0_{1\times n-1}& 0_{1\times (n(n-1)/2)}\\
    \end{array} \right)
    \left( \begin{array}{c}
        \p_{x_k} X \\
         \p_{x_k}{\bf \Xi}^{(2)}
    \end{array}\right) +P(X,\Xi)  \left( \begin{array}{c}
        \p_{x_k} X \\
         \p_{x_k}{\bf \Xi}^{(2)}
    \end{array}\right)\\
    &=:M \left( \begin{array}{c}
        \p_{x_k} X \\
         \p_{x_k}{\bf \Xi}^{(2)}
    \end{array}\right)+P(X,\Xi) \left( \begin{array}{c}
        \p_{x_k} X \\
         \p_{x_k}{\bf \Xi}^{(2)}
    \end{array}\right)
    \end{split}
\end{equation}
where, setting $n^\ast:=n+(n+1)n/2$, $M \in \mathsf{M}_{n^\ast}(\R)$ and, for $(X,\Xi) \in \R^{2n}$, $P(X,\Xi)$ is a matrix with the following block structure

\begin{align}
P(X,\Xi)
=\left(
\begin{array}{cccccc}
P^{1,1}
& P^{1,2}
& P^{1,3}
& P^{1,4}
& \cdots
& P^{1, n^\ast}
\\
P^{2,1}
& P^{2,2}
& P^{2,3}
& P^{2,4}
& \cdots
& P^{2,n^\ast}
\\
\vdots
& \vdots
& \vdots
& \vdots
& \vdots
& \vdots
\\
P^1_{n(n-1)/2\times 1}
& P^2_{ n(n-1)/2 \times 1}
& P^3_{n(n-1)/2\times 1}
& P^4_{ n(n-1)/2 \times 1}
& \cdots
& P^{n^\ast}_{n(n-1)/2\times 1}
\end{array}
\right)
\in \mathsf{M}_{n^\ast}(\R),
\end{align}
where the entries $P^{i,j}=P^{i,j}(X,\Xi)$ satisfy

\begin{align}
    P^{i,j}&=O(\la X\ra^{-3}|\Xi|^2),\,\,\forall i,j=1,\ldots,n;\\
    P^{i,j}&=O(\la X\ra^{-2}),\,\,\forall i=1,\ldots,n,\,\,\forall j=n+1,\ldots,n^*;\\
    P^{i,j}&=O(\la X\ra^{-4}|\Xi|^4),\,\,\forall i=n+1,\ldots,n^*,\,\,\forall j=1,\ldots,n;\\
    P^{i,j}&=O(\la X\ra^{-3}|\Xi|^2),\,\,\forall i=n+1,\ldots,n^*,\,\,\forall j=n+1,\ldots,n^*. 
\end{align}

Hence, by applying the variation of constants formula to \eqref{eq.dsdx} we obtain
\begin{equation}
\begin{split}
     \left(\begin{array}{c}
         \p_{x_k}X(s) \\
         \p_{x_k}{\bf \Xi}^{(2)}(s)
    \end{array} \right)=e^{sM}
    \left( \begin{array}{c}
        \p_{x_k} X(0) \\
         \p_{x_k}{\bf \Xi}^{(2)}(0)
    \end{array}\right)+\int_0^s e^{(s-\sigma)M}P(X(\sigma),\Xi(\sigma))  \left( \begin{array}{c}
        \p_{x_k} X(\sigma) \\
         \p_{x_k}{\bf \Xi}^{(2)}(\sigma)
    \end{array}\right)d\sigma,
    \end{split}
\end{equation}
and then, since $M^2=0$,
\begin{equation}
\begin{split}
     \left(\begin{array}{c}
         \p_{x_k}X(s) \\
         \p_{x_k}{\bf \Xi}^2(s)
    \end{array} \right)&=(I_{n^\ast}+sM)
    \left( \begin{array}{c}
        \p_{x_k} X(0) \\
         \p_{x_k}{\bf \Xi}^{(2)}(0)
    \end{array}\right)\\
    &\quad +\int_0^s (I_{n^\ast}+(s-\sigma)M)P(X(\sigma),\Xi(\sigma)) \left( \begin{array}{c}
        \p_{x_k} X(\sigma) \\
         \p_{x_k}{\bf \Xi}^{(2)}(\sigma)
    \end{array}\right)d\sigma.
    \end{split}
\end{equation}

Next, we observe that
\begin{equation}
\begin{split}
    M \cdot P(X,\Xi)&=
    \left(
\begin{array}{cccccc}
3P^{n+1,1}
& 3P^{n+1,2}
& 3P^{n+1,3}
& 3P^{n+1,4}
& \cdots
& 3P^{n+1,n^\ast}\\
0
& 0
& 0
& 0
& \cdots
& 0
\\
\vdots
& \vdots
& \vdots
& \vdots
& \vdots
& \vdots
\\
0
& 0
& 0
& 0
& \cdots
& 0
\end{array}
\right)\\
&\quad +\left(
\begin{array}{cccccc}
0
& 0
& 0
& 0
& \cdots
& 0
\\
2P^{n+2,1}
& 2P^{n+2,2}
& 2P^{n+2,3}
& 2P^{n+2,4}
& \cdots
& 2P^{n+2,n^\ast}
\\
2P^{n+3,1}
& 2P^{n+3,2}
& 2P^{n+3,3}
& 2P^{n+3,4}
& \cdots
& 2P^{n+3,n^\ast}
\\
\vdots
& \vdots
& \vdots
& \vdots
& \vdots
& \vdots\\
2P^{2n,1}
& 2P^{2n,2}
& 2P^{2n,3}
& 2P^{2n,4}
& \cdots
& 2P^{2n,n^\ast}
\\
0
& 0
& 0
& 0
& \cdots
& 0
\\
\vdots
& \vdots
& \vdots
& \vdots
& \vdots
& \vdots
\\
0
& 0
& 0
& 0
& \cdots
& 0
\end{array}
\right) \\
&\quad +
\left(
\begin{array}{cccccc}
P^{1}_{\mathfrak{I}}
& P^{2}_{\mathfrak{I}}
& P^{3}_{\mathfrak{I}}
& P^{4}_{\mathfrak{I}}
& \cdots
& P^{n^\ast}_{\mathfrak{I}}\\
0
& 0
& 0
& 0
& \cdots
& 0
\\
\vdots
& \vdots
& \vdots
& \vdots
& \vdots
& \vdots
\\
0
& 0
& 0
& 0
& \cdots
& 0
\end{array}
\right),
\end{split}
\end{equation}
where 
\[
\begin{split}
P_{\mathfrak{I}}^j&=O(\langle X \rangle^{-4}|\Xi|^4) \quad \forall j=1,\dots, n,\\
P_{\mathfrak{I}}^j&=O(\langle X \rangle^{-3}|\Xi|^2) \quad \forall j=n+1,\dots, n^\ast.
\end{split}
\]

We now set 
\begin{align}
    y(s):=\left( 1+\frac{|s||\xi|^2}{\la x\ra^{2}}\right)^{-1}|\p_xX(s)|+\frac{\la x\ra^2}{|\xi|^2}|\p_x{\bf \Xi}^{(2)}(s)|,
\end{align}
and have, for some $c_0>0$,
\begin{align}
    &y(s)\leq \left( 1+\frac{|s||\xi|^2}{\la x\ra^{2}}\right)^{-1}(|\p_xX(0)|+3
    |s||\p_x{\bf\Xi}^{(2)}(0)| )+\\
    &+ c_0\int_0^s  \left( 1+\frac{|s||\xi|^2}{\la x\ra^{2}}\right)^{-1}\left( 1+\frac{|\sigma||\xi|^2}{\la x\ra^{2}}\right)
\left( \frac{|\Xi(\sigma)|^2}{\la X(\sigma)\ra^3}+3|s-\sigma||\Xi(\sigma)|^4\la X(\sigma)\ra^{-4}\right)\\
&\quad \quad \times \left(1+ \frac{|\sigma||\xi|^2}{\la x\ra^{2}}\right)^{-1}|\p_xX(\sigma)|d\sigma\\
    &+c_0\int_0^s  \left(1+ \frac{|s||\xi|^2}{\la x\ra^{2}}\right)^{-1}\left( \frac{|\xi|^2}{\la x\ra^{2}}\right)
    \left( |s-\sigma||\Xi(\sigma)|^2\la X(\sigma)\ra^{-3} +\la X(\sigma)\ra^{-2}\right)
    \left( \frac{|\xi|^2}{\la x\ra^{2}}\right)^{-1}|\p_x{\bf\Xi}^{(2)}(\sigma)| d\sigma\\
   &+\frac{\la x\ra^2}{|\xi|^2} |\p_x{\bf\Xi}^{(2)}(0)|\\
   &+c_0\int_0^s  \frac{\la x\ra^2}{|\xi|^2} \left( 1+\frac{|\sigma||\xi|^2}{\la x\ra^{2}}\right)
    \frac{|\Xi(\sigma)|^4}{\la X(\sigma) \ra^4}  \left( 1+\frac{|\sigma||\xi|^2}{\la x\ra^{2}}\right)^{-1} |\p_{x}X(\sigma)| d\sigma\\
    &+ c_0\int_0^s  \frac{\la x\ra^2}{|\xi|^2} \la X(\sigma)\ra^{-3}|\Xi(\sigma)|^2 |\p_x{\bf\Xi}^{(2)}(\sigma)|d\sigma\\
    \leq & y(0)+\int_0^s  k(\sigma)y(\sigma) d\sigma.
\end{align} 

Next, we identify $k(\sigma)$. By the previous lemma, for $(x,\xi)\in\mathscr{S}_+$, we have
$$\la X(s)\ra\approx \la x \ra +|s||\xi|^2,$$
which is increasing in $|s|$. Moreover,  $\left( 1+\frac{|s||\xi|^2}{\la x\ra^{2}}\right)$ is increasing in $|s|$ and bounded from below by 1, while $|\Xi(s)|\approx |\xi|$. Therefore, for some $c_1>0$, possibly different from line to line, 
\begin{align}
    \left( 1+\frac{|s||\xi|^2}{\la x\ra^{2}}\right)^{-1}\left( 1+\frac{|\sigma||\xi|^2}{\la x\ra^{2}}\right)
\left( \frac{|\Xi(\sigma)|^2}{\la X(\sigma)\ra^3}\right) \leq&
\left( \frac{\la x\ra^2}{\la x\ra^{2}+|\sigma||\xi|^2}\right)\left( 1+\frac{|\sigma||\xi|^2}{\la x\ra^{2}}\right)
\left( \frac{|\xi|^2}{(\la x\ra+|\sigma||\xi|^2)^3}\right)\\
&\leq c_1\frac{|\xi|^2}{(\la x\ra+|\sigma||\xi|^2)^3}
\end{align}
and

\begin{align}
     &\left( 1+\frac{|s||\xi|^2}{\la x\ra^{2}}\right)^{-1}\left( 1+\frac{|\sigma||\xi|^2}{\la x\ra^{2}}\right)|s-\sigma||\Xi(\sigma)|^4\frac{1}{\la X(\sigma) \ra^4} \\
     \leq& 
     c_1\left( \frac{\la x\ra^2}{\la x\ra^{2}+|s||\xi|^2}\right)\left( 1+\frac{|\sigma||\xi|^2}{\la x\ra^{2}}\right)|s||\xi|^4\frac{1}{(\la x \ra+|\sigma||\xi|^2)^4} \\
     \leq & c_1\frac{|\xi|^2 }{(\la x \ra+|\sigma||\xi|^2)^2} + c_1\frac{|\sigma||\xi|^2}{\la x\ra^{2}} \frac{\la x\ra^2}{\la x\ra^{2}+|s||\xi|^2} |s||\xi|^4\frac{1}{(\la x \ra+|\sigma||\xi|^2)^4}\\
     \leq & c_1\frac{|\xi|^2 }{(\la x \ra+|\sigma||\xi|^2)^2} + c_1\frac{|\xi|^2 }{(\la x \ra+|\sigma||\xi|^2)^3} \\
     \leq & c_1\frac{|\xi|^2 }{(\la x \ra+|\sigma||\xi|^2)^2}.
\end{align}

On the other hand
\begin{align}
    \left(1+ \frac{|s||\xi|^2}{\la x\ra^{2}}\right)^{-1}\left( \frac{|\xi|^2}{\la x\ra^{2}}\right)
    \left(|s-\sigma||\Xi(\sigma)|^2\la X(\sigma)\ra^{-2} +\la X(\sigma)\ra^{-2}\right)\leq c_1\frac{|\xi|^2}{(\la x\ra+|s||\xi|^2)^2},
\end{align}

\begin{align}
    \frac{\la x\ra^2}{|\xi|^2}  \left( 1+\frac{|\sigma||\xi|^2}{\la x\ra^{2}}\right)
    \frac{|\Xi(\sigma)|^4}{\la X(\sigma) \ra^4} \leq c_1\frac{|\xi|^2}{(\la x\ra+|\sigma||\xi|^2)^2},
\end{align}
and
\begin{align}
    |\Xi(\sigma)|^2\la X(\sigma)\ra^{-2}\leq c_1\frac{|\xi|^2}{(\la x\ra+|\sigma||\xi|^2)^2}.
\end{align}
Therefore, for $C>0$ sufficiently large, we set
$$k(\sigma):=C\frac{|\xi|^2}{(\la x\ra+|\sigma||\xi|^2)^2}.$$
We have that 
$$\int_0^{+\infty} k(\sigma)d\sigma\lesssim \la x \ra^{-1}.$$
Moreover,
$$y(s)\leq y(0)+\int_0^s k(\sigma)y(\sigma)d\sigma.$$
which yields, by Propositon \ref{prop.Gronwall}, applied after the change of variable $s\mapsto -s$,
\begin{align}
    |\p_{x_k}X(s)|&\lesssim \left(1+ \frac{|s||\xi|^2}{\la x\ra^{2}}\right)\\
    |\p_{x_k}{\bf \Xi}^{(2)}(s)|&\lesssim \frac{|\xi|^2}{\la x\ra^2},
\end{align}
for all $k=1,\ldots,n$. Finally, by Proposition \ref{eq.usefineq}, we have
$$|\p_{x_k}\Xi(s)|\leq C |\p_{x_k}{\bf \Xi}^{(2)}(s)|/|\Xi(s)|\lesssim  \frac{|\xi|}{\la x\ra^2},$$
which concludes the proof of the first part of the statement.

To prove the bound for the $\xi$-derivatives, one proceeds similarly.
In fact, by repeating the same argument as above we obtain the system 

\begin{equation}\label{eq2801}
\begin{split}
     \left(\begin{array}{c}
         \p_{\xi_k}X(s) \\
         \p_{\xi_k}{\bf \Xi}^{(2)}(s)
    \end{array} \right)&=(I_{n^\ast}+sM)
    \left( \begin{array}{c}
        \p_{\xi_k} X(0) \\
         \p_{\xi_k}{\bf \Xi}^{(2)}(0)
    \end{array}\right)\\
    &\quad +\int_0^s (I_{n^\ast}+(s-\sigma)M)P(X(\sigma),\Xi(\sigma)) \left( \begin{array}{c}
        \p_{\xi_k} X(\sigma) \\
         \p_{\xi_k}{\bf \Xi}^{(2)}(\sigma)
    \end{array}\right)d\sigma,
    \end{split}
\end{equation}
where the matrix $P$ is the same as the one above (one has to repeat the above calculations replacing $\p_{x_k}X, \p_{x_k}\Xi^{(2)}$ with $\p_{\xi_k}X, \p_{\xi_k}\Xi^{(2)}$).
In this case, to obtain \eqref{eq.estderflow3} and \eqref{eq.estderflow4} for a fixed $k \in \lbrace 1,\dots,n \rbrace $, the idea is to use Proposition \ref{prop.Gronwall} on
\[
y_1(s):=\frac{1}{\la s \ra |\xi|}|\p_{\xi_k} X(s)|+\frac{1}{|\xi|}|\p_{\xi_k} {\bf\Xi}^{(2)}(s)|, \quad s \in [0,+\infty),
\]
in the case $|\xi|=1$ 
to get \eqref{eq.estderflow3} and \eqref{eq.estderflow4} in this case first.
The general case $|\xi|\neq 1$ will be obtained by using the homogeneity properties of the flow.

We first observe that
$$|\xi|\leq|\p_{\xi_k}{\bf\Xi}^{(2)}(0)|\leq 2n|\xi|,$$
so, we have
\begin{align}
    y_1(s)&\leq\frac{1}{|\xi|\langle s \rangle}(|\p_{\xi_k} X(0)|+3|s||\p_{\xi_k} {\bf \Xi}^{(2)}(0)|) \\
    &\quad + c_0\int_0^s  \frac{\la \sigma\ra}{|\xi|\langle s \rangle}
\left( \frac{|\Xi(\sigma)|^2}{\la X(\sigma)\ra^3}+|s-\sigma||\Xi(\sigma)|^4\la X(\sigma)\ra^{-4}\right)
    \frac{1}{\la \sigma\ra}|\p_{\xi_k} X(\sigma)| d\sigma\\
    &\quad +c_0\int_0^s  \frac{1}{|\xi|\langle s \rangle}
    \left( \frac{1}{\la X(\sigma)\ra^2}+|s-\sigma||\Xi(\sigma)|^2\la X(\sigma)\ra^{-3}\right)
   |\p_{\xi_k}{\bf\Xi}^{(2)}(\sigma)| d\sigma\\
   &\quad +\frac{1}{|\xi|}|\p_{\xi_k}{\bf \Xi}^{(2)}(0)| \\
   &\quad +c_0\int_0^s  \frac{\la \sigma\ra}{|\xi|}
    \frac{|\Xi(\sigma)|^4}{\la X(\sigma) \ra^4}  \frac{1}{\la \sigma\ra}|\p_{\xi_k}X(\sigma)| d\sigma+ c_0\int_0^s  \frac{1}{|\xi|}\la X(\sigma)\ra^{-3}|\Xi(\sigma)|^2 |\p_{\xi_k}{\bf\Xi}^{(2)}(\sigma)|d\sigma\\
    &\leq c +\int_0^s  k_1(\sigma)y_1(\sigma) d\sigma,
\end{align}
for some $c,c_0>0$.
Next, we restrict ourselves to $|\xi|=1$ and determine $k_1$ in this case.

We note that for $c_1>0$ sufficiently large, 
\begin{align}
    &\frac{\la \sigma\ra}{\langle s \rangle}
 \frac{|\Xi(\sigma)|^2}{\la X(\sigma)\ra^3}\leq c_1\frac{|\xi|^2}{\la |x|+|\sigma||\xi|^2\ra^3}=c_1\frac{1}{\la |x|+|\sigma|\ra^3}\\
 & \frac{\la \sigma\ra}{\langle s \rangle} |s-\sigma|\frac{|\Xi(\sigma)|^4}{\la X(\sigma)\ra^{4}}\leq  c_1 \frac{|\xi|^2}{( \la x\ra+|\sigma||\xi|^2)^3}= c_1 \frac{1}{( \la x\ra+|\sigma|)^3} \\
 &\frac{1}{\la s\ra \la X(\sigma)\ra^2}\leq c_1\frac{1}{ ( \la x\ra+|\sigma|)^2},\
\end{align}
have uniformly bounded integral over $[0,+\infty)$.
Hence, by taking 
\[
k_1(\sigma):=C \frac{1}{( \la x\ra+|\sigma|)^2}
\]
with $C>0$ sufficiently large, as before, by using Proposition \ref{prop.Gronwall}, we conclude
\begin{align}
    &|\p_{\xi_k}X(s)|\lesssim \la s\ra , \label{eq.s}\\
    &|\p_{\xi_k}\Xi^{(2)}(s)|\lesssim  1\\
    &|\p_{\xi_k}\Xi(s)|\lesssim  1\label{eq.xi}
\end{align}
for every $k=1,\ldots,n$.
These imply \eqref{eq.estderflow3}, with $\la s\ra$ in place of $|s|$ in the right hand side, and \eqref{eq.estderflow4}, whenever $|\xi|=1$.
To obtain \eqref{eq.estderflow3} (with $|s|$ in place of $\la s\ra$) from \eqref{eq.s}  we distinguish the cases $s\geq 1$ and $s\in[0,1)$.

If $s\geq 1$, we have $\la s\ra \lesssim |s|$ so \eqref{eq.s} gives \eqref{eq.estderflow3} when $|\xi|=1$.

If $s\in [0,1)$, then \eqref{eq.s} and $\eqref{eq.xi}$ give 

\[|\p_{\xi_k}X(s)|, |\p_{\xi_k}\Xi(s)| \leq C,\quad \forall k=1,\ldots,n\]
therefore, for all $\ell=1,\ldots,n,$
\[
\begin{aligned}
\frac{d}{ds}\partial_{\xi_k}X_\ell(s)
={}&\sum_{j=1}^n
(\partial_{x_j}\partial_{\xi_\ell}a_3)(X(s),\Xi(s))
\,\partial_{\xi_k}X_j(s)\\
&+\sum_{j=1}^n
(\partial_{\xi_j}\partial_{\xi_\ell}a_3)(X(s),\Xi(s))
\,\partial_{\xi_k}\Xi_j(s).
\end{aligned}
\]
Since the coefficients in the sums above are uniformly bounded when $|\xi|=1$ due to the properties of $A_2$, we get
\[
\left|\frac{d}{ds}\partial_{\xi_k}X(s)\right|
\le C\bigl(
|\partial_{\xi_k}X(s)|+
|\partial_{\xi_k}\Xi(s)|
\bigr)\le C_1,
\qquad 0\le s\le1.
,\quad \text{for } |\xi|=1,\]
which yields
\[
|\partial_{\xi_k}X(s)|
\le\int_0^s|\frac{d}{d\sigma}\,\partial_{\xi_k}X(\sigma;x,\xi)|\,d\sigma
\le C_1(s-0)\leq C_1|s|,
\qquad 0\le s\le1,
\]
where, note, we have used $\left.\partial_{\xi_k}X_j(s;x,\xi)\right|_{s=0}=
\partial_{\xi_k}\bigl(X_j(0;x,\xi)\bigr)=\partial_{\xi_k} x_j=0$. This, finally, gives \eqref{eq.estderflow3} when $s\in[0,1)$, so, adding the case $s\geq 1$, we conclude 
\begin{equation}\label{eq.ss}
    |\partial_{\xi_k}X(s)|\lesssim |s|,\quad \text{when } |\xi|=1,
\end{equation}
for all $s\geq 0$.

To show \eqref{eq.estderflow3} and \eqref{eq.estderflow4} when $|\xi|\neq 1$, we use \eqref{eq.ss} and \eqref{eq.xi} as follows.

Write $r=|\xi|$ and $\omega=\xi/ r$. Differentiating the homogeneity identities with $r$ fixed, we obtain
\[
\partial_{\xi_k}X(s;x,r\omega)
=
r^{-1}
\left.\partial_{\eta_k}X(r^2s;x,\eta)\right|_{\eta=\omega},
\]
\[
\partial_{\xi_k}\Xi(s;x,r\omega)
=
\left.\partial_{\eta_k}\Xi(r^2s;x,\eta)\right|_{\eta=\omega}.
\]
Then, since $|\omega|=1$, estimates \eqref{eq.ss} and \eqref{eq.xi} at time \(r^2s\) give
\[
|\partial_{\xi_k}X(s;x,\xi)|\le C|s||\xi|,
\qquad
|\partial_{\xi_k}\Xi(s;x,\xi)|\le C,
\]
uniformly for all \(s\geq 0\) and \((x,\xi)\in\mathscr S_+\). This completes the proof of the first-order estimates.
\medskip

The estimates for higher-order derivatives follow from by combinining the multivariate Fa\'a di Bruno formula stated in Theorem \ref{theo.fdb}, and, once again, a Gronwall-type inequality. The essential difference with the case in \cite{CKS} is that we work with a different homogeneity, induced by the principal symbol $a_3$ of $A$ of order three, whereas in \cite{CKS} a second-order symbol is considered. Accordingly, as in the previous estimates, the argument is adapted by choosing a suitable function to which the Gronwall inequality is applied, reflecting this different homogeneity.

Specifically, for each fixed $\alpha,\beta \in \mathbb{N}_0^n$ with $|\alpha|+|\beta|>1$, by Theorem \ref{theo.fdb}, as in \eqref{eq.dsdx}, we have 
\begin{equation}\label{eq.dersup}
\begin{split}
    \ds \left(\begin{array}{c}
         \p_x^\beta \p_\xi^\alpha X \\
         \p_x^\beta \p_\xi^\alpha{\bf \Xi}^{(2)}
    \end{array} \right)&=M \left( \begin{array}{c}
        \p_x^\beta \p_\xi^\alpha X \\
         \p_x^\beta \p_\xi^\alpha{\bf \Xi}^{(2)}
    \end{array}\right)+P(X,\Xi) \left( \begin{array}{c}
        \p_x^\beta \p_\xi^\alpha X \\
        \p_x^\beta \p_\xi^\alpha{\bf \Xi}^{(2)}
    \end{array}\right)\\
    &\quad + \sum c_{\alpha,\beta,\ell,\lambda,\mu} (\p_X^\gamma \p_\Xi^\delta P)(X,\Xi)
    \left(\begin{array}{c}
        \prod_j [\p_x^{\ell_j^1}\p_\xi^{\ell_j^2} X]^{\lambda_j^1} [\p_\xi^{\ell_j^3}\p_\xi^{\ell_j^4}{\bf \Xi}^{(2)}]^{\lambda_j^2} \\
        \vdots \\
        \prod_j [\p_x^{\ell_j^1}\p_\xi^{\ell_j^2} X]^{\mu_j^1} [\p_\xi^{\ell_j^3}\p_\xi^{\ell_j^4}{\bf \Xi}^{(2)}]^{\mu_j^2} \\
        \vdots
    \end{array}\right),
    \end{split}
\end{equation}
where $\partial_X^\gamma \partial_\Xi^\delta P$ denotes the matrix (of the same size as $P$) whose $(i,j)$-entry is given by $\partial_X^\gamma \partial_\Xi^\delta P_{ij}$, and the multiindexes on the sum satisfy suitable relations (see Theorem \ref{theo.fdb}).
Therefore, by applying the variation of constants formula to \eqref{eq.dersup}, as before we get 
\[
\begin{split}
    &\left(\begin{array}{c}
         \p_x^\beta \p_\xi^\alpha X \\
         \p_x^\beta \p_\xi^\alpha{\bf \Xi}^{(2)}
    \end{array} \right)\\
    &= (I_{n^\ast}+sM)
    \left( \begin{array}{c}
        \p_x^\beta \p_\xi^\alpha X(0) \\
         \p_x^\beta \p_\xi^\alpha{\bf \Xi}^{(2)}(0)
    \end{array}\right)\\
    &+\int_0^s (I_{n^\ast}+(s-\sigma)M)P(X(\sigma),\Xi(\sigma)) \left( \begin{array}{c}
        \p_x^\beta \p_\xi^\alpha X(\sigma) \\
        \p_x^\beta \p_\xi^\alpha{\bf \Xi}^{(2)}(\sigma)
    \end{array}\right)d\sigma\\
    &\quad +\int_0^s (I_{n^\ast}+(s-\sigma)M)\sum c_{\alpha,\beta,\ell,\lambda,\mu} (\p_X^\gamma \p_\Xi^\delta P)(X(\sigma),\Xi(\sigma))
    \left(\begin{array}{c}
        \prod_j [\p_x^{\ell_j^1}\p_\xi^{\ell_j^2} X(\sigma)]^{\lambda_j^1} [\p_\xi^{\ell_j^3}\p_\xi^{\ell_j^4}{\bf \Xi}^{(2)}(\sigma)]^{\lambda_j^2} \\
        \vdots \\
        \prod_j [\p_x^{\ell_j^1}\p_\xi^{\ell_j^2} X(\sigma)]^{\mu_j^1} [\p_\xi^{\ell_j^3}\p_\xi^{\ell_j^4}{\bf \Xi}^{(2)}(\sigma)]^{\mu_j^2} \\
        \vdots
    \end{array}\right)d\sigma.
    \end{split}
\]
We next observe that 
\[
\left( \begin{array}{c}
        \p_x^\beta \p_\xi^\alpha X(0) \\
         \p_x^\beta \p_\xi^\alpha{\bf \Xi}^{(2)}(0)
    \end{array}\right)= \left( \begin{array}{c}
        0 \\
         \p_x^\beta \p_\xi^\alpha{\bf \Xi}^{(2)}(0)
    \end{array}\right)=\left( \begin{array}{c}
        0 \\
         1_{\lbrace|\beta|=0, |\alpha=2|\rbrace}\p_\xi^\alpha{\bf \Xi}^{(2)}(0)
    \end{array}\right)
\]
At this point,  the estimates for the higher order derivatives are obtained by applying the same Gronwall inequality as in \cite{CKS} to the function
\[
y(s):=\frac{\langle x \rangle^{k(\alpha,\beta)}}{\langle s \rangle |\xi|^{2-|\alpha|}}|\partial_x^\beta \partial_\xi^\alpha X(s;x,\xi)|+\frac{\langle x \rangle^{k(\alpha,\beta)}}{|\xi|^{1-|\alpha|}}|\partial_x^\beta \partial_\xi^\alpha \Xi(s;x,\xi)|, \quad s \in [0,+\infty),
\]
where $k(\alpha,\beta)>|\beta|$.

As above, the regime of small $|s|$ near zero must be treated separately, by working directly with the system induced by the Hamilton system associated with the Hamiltonian $a_3$.
\endproof

Next, we derive estimates for the scattering variables.

\begin{theorem}\label{theo.dersv}
If $(x,\xi) \in \mathscr{S}_+$ the scattering variables $(\widebar{X}(x,\xi),\widebar{\Xi}(x,\xi))$ defined in Lemma \ref{lem.scatteringV} satisfy
    \begin{align}
    \left|\p_{\xi_i}\widebar{\Xi}_j(x,\xi)-\delta_{ij}\right|&\lesssim \frac{1}{\langle x \rangle^2}\label{eq.lem.xibar2}\\
    |\partial_{x_i}\widebar{\Xi}_j(x,\xi)\big| &\lesssim \frac{|\xi|}{\langle x\rangle^3}\label{eq.lem.xibar1}\\
     |\partial_{\xi_i}\widebar{X}_j(x,\xi)\big| &\lesssim \frac{1}{|\xi|\langle x\rangle}\label{eq.lem.xbar2}\\
    |\p_{x_i}\widebar{X}_j(x,\xi)-\delta_{ij}|&\lesssim \frac{1}{\la x\ra^2}. \label{eq.lem.xbar1} 
\end{align}
 Moreover, if $\alpha,\beta \in \mathbb{N}_0^n$ with $|\alpha|+|\beta|>1$, there exists $C_{\alpha,\beta}>0$ such that 
    \begin{align}
        |\p_x^\beta\p_\xi^\alpha \widebar X(x,\xi)|&\leq C_{\alpha,\beta}\frac{|\xi|^{-|\alpha|}}{\la x\ra^{|\beta|+1}}\,\label{eq.xbarhighder}\\  
        |\p_x^\beta\p_\xi^\alpha \widebar \Xi(x,\xi)|&\leq C_{\alpha,\beta}\frac{|\xi|^{1-|\alpha|}}{\la x\ra^{|\beta|+2}}\label{eq.xibarhighder}.
    \end{align}
\end{theorem}

\proof
We start by proving \eqref{eq.lem.xibar1} and \eqref{eq.lem.xibar2}.

Recall that
\begin{align}
     \widebar{\Xi}_j(x,\xi):=\xi_j+\int_0^{+\infty} \frac{d}{ds}\Xi_j(s)ds.\label{eq.xibar}
\end{align}
Then, on denoting by $\Xi(s):=\Xi(s;x,\xi)$ and $X(s)=X(s;x,\xi)$, for $s\geq 0$, we have
\begin{align}
\partial_{x_i} \ds \Xi_j(s)
&=-
\partial_{x_i}\Bigg[
\Xi_1(s)
\sum_{k,\ell=1}^n
\partial_{X_j} a_{k,\ell}(X(s))\,
\Xi_k(s)\,
\Xi_\ell(s)
\Bigg] \\
&=-
\partial_{x_i}\Xi_1(s) b(s)
- \Xi_1(s)\cdot
\Bigg\{
\sum_{k,\ell,m=1}^n
\partial^2_{X_m X_j} a_{k,\ell}(X(s))\,
(\partial_{x_i} X_m(s))\,
\Xi_k(s)\,\Xi_\ell(s)\\
&+ \sum_{k,\ell=1}^n
(\partial_{X_j}a_{k,\ell})(X(s))\,
\Big[
(\partial_{x_i}\Xi_k(s))\,\Xi_\ell(s)
+ \Xi_k(s)\,(\partial_{x_i}\Xi_\ell(s))
\Big]
\Bigg\},
\end{align}
with 
\[
b(s):=
\sum_{k,\ell=1}^n
\partial_{X_j} a_{k,\ell}(X(s))\,
\Xi_k(s)\,
\Xi_\ell(s).
\]
Hence,
\begin{equation}\label{eq.pxiXij}
\begin{split}
\big|\partial_{x_i}\frac{d}{ds}\Xi_j(s)\big|
&\le
\big|\partial_{x_i}\Xi_1(s)\big|\,|b(s)|
+ |\Xi(s)|\sum_{k,\ell,m}
\big|\partial^2_{X_m X_j} a_{k,\ell}(X(s))\big|\,
\big|\partial_{x_i}X_m(s)\big|\,
|\Xi_k(s)|^2\\
&
+ \sum_{k,\ell=1}^n
|(\partial_{X_j} a_{k,\ell})(X(s))|
\Big(
|\partial_{x_i}\Xi_k(s)|\,|\Xi_\ell(s)|
+ |\Xi_k(s)|\,|\partial_{x_i}\Xi_\ell|
\Big)\\
&\lesssim
\frac{|\xi|}{\langle x\rangle^2}\,
\frac{|\xi|^2}{\langle X(s)\rangle^3}
+ \frac{|\xi|}{\langle X(s)\rangle^4}
\Big(1 + |s||\xi|^2\Big)
\frac{|\xi|^2}{\langle x\rangle^2}
+ \frac{|\xi|}{\langle x\rangle^2}\,
\frac{|\xi|^2}{\langle X(s)\rangle^4}.
\end{split}
\end{equation}
Now, using that
\[
\langle X(s)\rangle \approx \langle x \rangle + |s||\xi|^2\quad 
\]
we have
\[
\int_0^{+\infty} \frac{|\xi|^2}{ \langle X(s)\rangle^\alpha} ds\lesssim
\frac{1}{ \langle x\rangle^{\alpha-1}},\quad \alpha\geq 2,\]
and
\[
\int_0^{+\infty} \frac{|s||\xi|^4}{ \langle X(s)\rangle^4} ds\lesssim \int_0^{+\infty} \frac{|\xi|^2}{ \langle X(s)\rangle^2} ds.
\]
Therefore, by integrating in \eqref{eq.pxiXij},
we conclude that 
\[
\int_0^{+\infty} |\partial_{x_i}\frac{d}{ds}\Xi_j(s)\big|ds \lesssim \frac{|\xi|}{\langle x\rangle^3},
\]
and thus
\[|\partial_{x_i}\widebar{\Xi}_j(x,\xi)\big| \lesssim \frac{|\xi|}{\langle x\rangle^3},\]
which gives \eqref{eq.lem.xibar1}.

By using similar strategies we estimate
\begin{align}
\p_{\xi_i}\widebar{\Xi}_j(x,\xi)-\delta_{ij}&=\int_0^\infty
(\partial_{\xi_i}\Xi_1(s)) b(s)
+ \Xi_1(s)\cdot
\Bigg\{
\sum_{k,\ell,m=1}^n
\partial^2_{X_m X_j} a_{k,\ell}(X(s))\,
(\partial_{\xi_i} X_m(s))\,
\Xi_k(s)\,\Xi_\ell(s)\\
&+ \sum_{k,\ell}
(\partial_{X_j}a_{k,\ell})(X(s))\,
\Big[
(\partial_{\xi_i}\Xi_k(s))\,\Xi_\ell(s)
+ \Xi_k(s)\,(\partial_{\xi_i}\Xi_\ell(s))
\Big]
\Bigg\},
\end{align}
and get, using that $|\Xi(s)|\approx |\xi|$ for $(x,\xi)\in\mathscr{S}_+,$
\begin{align}
    \left|\p_{\xi_i}\widebar{\Xi}_j(x,\xi)-\delta_{ij}\right|\lesssim \int_0^\infty \frac{|\xi|^2}{\langle X(s)\rangle^3} ds+\int_0^\infty \frac{s|\xi|^4}{\langle X(s)\rangle^4} ds\lesssim \frac{1}{\langle x \rangle^2},
\end{align}
that is \eqref{eq.lem.xibar2}.

Next, we focus on showing \eqref{eq.lem.xbar1} and \eqref{eq.lem.xbar2}. 
Recall that $\widebar{X}$ is defined as
\begin{equation}\label{eq.overlinex}
    \widebar{X}(x,\xi):=x+\int_0^{+\infty} \frac{d}{ds}Y(s) ds,
\end{equation}
where
\begin{equation}\label{eq.dtYs}
    \frac{d}{ds}Y(s):= \frac{d}{ds}X(s;x,\xi)-\big(\, \sum_{k=1}^nc_k\widebar{\Xi}^2_k(x,\xi), 2\widebar{\Xi}_1(x,\xi)\widebar{\Xi}_2(x,\xi),\ldots,2\widebar{\Xi}_1(x,\xi)\widebar{\Xi}_n(x,\xi)\,\big).
\end{equation}
We now estimate $\partial_{\xi_i} \widebar X_j$ for $i,j=1,\dots, n$.  To do so, we consider the cases $j=1$ and $j=2,\dots,n,$ separately, due to the different structure in \eqref{eq.dtYs}. 

For $j=1$ and $i=1,\dots, n,$ we have
(recall $c_1=3$ and $c_k=1$, for all $k=1,\dots, n$)
\begin{equation}\label{eq.dY1}
\begin{split}
\frac{d}{ds}Y_1(s)&=\frac{d}{ds}X_1(s;x,\xi)- \sum_{k=1}^nc_k\widebar{\Xi}^2_k(x,\xi)\\
&=\sum_{k,\ell=2}^n(a_{k\ell}(X(s))-\delta_{k\ell}))\Xi_k(s)\Xi_\ell(s)+\sum_{k=1}^nc_k(\Xi^{2}_k(s)-\widebar{\Xi}_k^2(x,\xi))\\
&=:I(s)+II(s)
\end{split}
\end{equation}
that is,
\begin{align}
    |\p_{\xi_i}\widebar{X}_1(x,\xi)|\leq \int_0^{+\infty} |\p_{\xi_i}I(s)|ds+\left|\int_0^{+\infty} \p_{\xi_i} II(s)ds\right|.\label{eq.partial1X1bar}
\end{align}
We then estimate the two integrals in \eqref{eq.partial1X1bar}.

Observe that
\begin{equation}
\begin{split}
    |\p_{\xi_i}I(s)| &\leq \Bigl| \Bigl(\partial_{\xi_i}\sum_{k,\ell=2}^n(a_{k\ell}(X(s))-\delta_{k\ell})\Bigr)\Xi_k(s)\Xi_\ell(s)\Bigr|+\Bigl|\sum_{k,\ell=2}(a_{k\ell}(X(s))-\delta_{k\ell})\partial_{\xi_i} \Xi_k(s)\Xi_\ell(s)\Bigr| \\
    &\leq \Bigl| \sum_{k,\ell=2}^n\Bigl(\nabla_X a_{k\ell}(X(s)) \cdot \partial_{\xi_i} X(s)\Bigr)\Xi_k(s)\Xi_\ell(s)\Bigr|\\
    &\quad  + \Bigl|\sum_{k,\ell=2}(a_{k\ell}(X(s))-\delta_{k\ell})\Bigl( \Xi_\ell(s) \partial_{\xi_i}\Xi_k(s)+\Xi_k(s)\partial_{\xi_i}\Xi_\ell(s)\Bigr)\Bigr|.
\end{split}
\end{equation}
Hence, by using \ref{hp3},  \eqref{lem-3.2} and \eqref{eq.estderflow3} on the first term and \ref{hp3}, \eqref{lem-3.2} and \eqref{eq.estderflow4} on the second one, we get
\[
|\p_{\xi_i}I(s)|\lesssim \frac{|s| |\xi|^3}{\langle X(s) \rangle^3}+\frac{|\xi|}{\langle X(s) \rangle^2},
\]
and thus, by Remark \ref{rmk.behavscattpoints}, we obtain
\begin{align}
    \int_0^{+\infty} |\p_{\xi_i}I(s)|ds\lesssim \int_0^{+\infty}\Bigl( \frac{|s| |\xi|^3}{\langle X(s) \rangle^3}+\frac{|\xi|}{\langle X(s) \rangle^2}\Bigr) ds\lesssim \frac{1}{|\xi|\la x\ra} \label{eq.pxiXb1},
\end{align}
which gives an estimate for the first term on the RHS of \eqref{eq.partial1X1bar}.

To estimate the second term on the RHS of \eqref{eq.partial1X1bar} we use \eqref{eq.xibar} to rewrite, for all $k=1,\dots,n$,
\[
\begin{split}
&\Xi_k^2(s)-\widebar{\Xi}_k^2=\\
&=(\Xi_k(s)-\widebar{\Xi}_k)(\Xi_k(s)+\widebar{\Xi}_k)\\
&=\Bigl(-\int_s^{+\infty}\dt \Xi_k(\sigma)d\sigma\Bigr)\Bigl(\Xi_k(s)+\xi_k+\int_0^{+\infty} \dt \Xi_k(s)ds \Bigr) \\
& = -\Xi_k(s)\int_s^{+\infty}\dt \Xi_k(\sigma)d\sigma-\xi_k\int_s^{+\infty}\dt \Xi_k(s)ds-\bigl(\int_0^{+\infty}\dt \Xi_k(s)ds\Bigr)\Bigl(\int_s^{+\infty}\dt \Xi_k(\sigma)d\sigma\Bigr),
\end{split}
\]
so that,
\begin{align}
    \int_0^\infty \p_{\xi_i}II(s) ds=&\sum_{k=1}^nc_k \Bigl[ 
    -\int_0^\infty(\p_{\xi_i}\Xi_k(s))\int_s^{+\infty}\dt \Xi_k(\sigma)d\sigma\, ds
    -\int_0^\infty\Xi_k(s)\int_s^{+\infty}\p_{\xi_i}\dt \Xi_k(\sigma)d\sigma\,ds\\
    &-\int_0^\infty(\p_{\xi_i}\xi_k)\int_s^{+\infty}\dt \Xi_k(\sigma)d\sigma\,ds -\int_0^\infty\xi_k\int_s^{+\infty}\dt \p_{\xi_i}\Xi_k(\sigma)d\sigma\,ds\\
    &-\bigl(\int_0^{+\infty}\dt \Xi_k(s)ds\Bigr)\Bigl(\int_0^\infty\int_s^{+\infty}\p_{\xi_i}\dt \Xi_k(\sigma)d\sigma\, ds\Bigr)\\
    &-\bigl(\int_0^{+\infty}\p_{\xi_i}\dt \Xi_k(s)ds\Bigr)\Bigl(\int_0^\infty\int_s^{+\infty}\dt \Xi_k(\sigma)d\sigma\, ds\Bigr)
    \Bigr]\\
    &= I_1+I_2+I_3+I_4+I_5+I_6.
\end{align}
We recall that, by Remark \ref{rmk.behavscattpoints}, we have
$\la X(s;x,\xi)\ra\approx \la x \ra+s|\xi|^2$. Moreover, 
\[\int_0^\infty \dt \Xi_k(s)ds\lesssim \frac{|\xi|}{\la x \ra},\quad \forall k=1,\ldots, n.\]
These property, together with \ref{hp3}, \eqref{lem-3.2}, \eqref{eq.estderflow3} and \eqref{eq.estderflow4}, give
\begin{align}
   \int_s^\infty|\dt \Xi_k(\sigma)|d\sigma\lesssim\frac{|\xi|}{\la X(s)\ra^2}\,\, \forall k=1,\ldots,n &\Longrightarrow |I_1|,|I_3|\lesssim \frac{1}{|\xi|\la x\ra}\label{est.int.xi}  
  \end{align}
  \begin{align}
   \int_s^\infty |\p_{\xi_i}\dt \Xi_k(\sigma)|d\sigma&\lesssim
   \int_s^\infty \left(\frac{\sigma |\xi|^4}{\la X(\sigma)\ra^4}+\frac{ |\xi|^2}{\la X(\sigma)\ra^3}\right) d\sigma&\\
   &\lesssim\frac{1}{\la X(s)\ra^2}\,\, \forall k=1,\ldots,n\Longrightarrow  |I_2|,|I_4|,|I_5|,|I_6|\lesssim \frac{1}{|\xi|\la x\ra},\label{est.int.derxi}
\end{align}
so we have
\begin{align}
   \left| \int_0^\infty \p_{\xi_i}II(s)ds\right|\leq |I_1|+|I_2|+\cdots+|I_6|\lesssim \frac{1}{|\xi|\la x\ra} \label{eq.pxiXb2},
\end{align}
and finally, plugging \eqref{eq.pxiXb1}  and \eqref{eq.pxiXb2} into \eqref{eq.partial1X1bar} we conclude that
\[|\p_{\xi_i}\widebar{X}_1(x,\xi)|\lesssim \frac{1}{|\xi|\la x\ra},\]
that is \eqref{eq.lem.xbar2} when $j=1$.

To treat the case $j\neq 1$ we go back to \eqref{eq.dtYs} and note that, by \eqref{eq.dsYj},

\begin{align}
    \ds Y_j(s;x,\xi)=&\ds \left(X_j(s;x,\xi)-2s \widebar{\Xi}_1(x,\xi)\widebar{\Xi}_j(x,\xi)\right) \\
    &= 2\sum_{i=2}^n (a_{ij}(X(s))-\delta_{ij})\Xi_1(s)\Xi_i(s)+2(\Xi_1(s)\Xi_j(s)-\widebar{\Xi}_1(x,\xi)\widebar{\Xi}_j(x,\xi))\\
    &= 2\sum_{i=2}^n (a_{ij}(X(s))-\delta_{ij})\Xi_1(s)\Xi_i(s)+2\Big( (\Xi_1-\widebar{\Xi}_1)\Xi_j-\widebar{\Xi}_1(\widebar{\Xi}_j-\Xi_j)\Big)\\
    &=I'(s)+II'(s),\label{eq.dYj}
\end{align}
so, for all $j\neq 1$ and for all $i=1,\ldots,n,$
\begin{align}
    |\p_{\xi_i}\widebar{X}_j(x,\xi)|\leq \int_0^{+\infty} |\p_{\xi_i}I'(s)|ds+\left|\int_0^{+\infty} \p_{\xi_i} II'(s)ds\right|.\label{eq.partial1Xjbar}
\end{align}

To estimate the first integral in \eqref{eq.partial1Xjbar}, we repeat the steps in the estimate of the first term in \eqref{eq.partial1X1bar} (where now $\ell=1$), and get
\begin{align}
    \int_0^{+\infty} |\p_{\xi_i}I'(s)|ds\lesssim\frac{1}{|\xi|\la x\ra}.\label{est.II1}
\end{align}

To estimate the second term in \eqref{eq.partial1Xjbar} we first use we use \eqref{eq.xibar} to write
\begin{align}
    (\Xi_1-\widebar{\Xi}_1)\Xi_j-\widebar{\Xi}_1(\widebar{\Xi}_j-\Xi_j)=&-\Xi_j(s)\int_s^\infty \dt \Xi_1(\sigma) d\sigma\\
    &-\Big(\xi_1+\int_0^\infty \dt \Xi_1(\sigma) d\sigma\Big)\int_s^\infty \dt \Xi_j(\sigma) d\sigma,
\end{align}
therefore, as above we get
\begin{align}
    &\int_0^\infty \p_{\xi_i} II'(s)ds\\
    &=2 \Bigl[-\int_0^\infty (\p_{\xi_i}\Xi_j(s))\int_s^\infty \dt \Xi_1(\sigma) d\sigma ds-\int_0^\infty \Xi_j(s)\int_s^\infty \p_{\xi_i}\dt \Xi_1(\sigma) d\sigma ds\\
    &-\int_0^\infty (\p_{\xi_i}\xi_1)\int_0^\infty\int_s^\infty \dt \Xi_j(\sigma) d\sigma  ds-\int_0^\infty \int_0^\infty \p_{\xi_i}\dt \Xi_1(\sigma) d\sigma \int_0^\infty\int_s^\infty \dt \Xi_j(\sigma) d\sigma  ds\\
     &-\Big(\xi_1+\int_0^\infty \dt \Xi_1(\sigma) d\sigma\Big) \int_0^\infty\int_s^\infty \p_{\xi_i}\dt \Xi_j(\sigma) d\sigma  ds\Bigr]\\
     &=J_1+J_2+J_3+J_4+J_5.
\end{align}
Finally, using the properties in \eqref{est.int.xi} and \eqref{est.int.derxi} we reach the following
\[|J_i|\lesssim\frac{1}{|\xi|\la x \ra}, \quad \forall i=1,\ldots,5,\]
that give
\begin{align}
   \left|\int_0^\infty \p_{\xi_i} II'(s)ds\right|\lesssim \frac{1}{|\xi|\la x \ra} \label{est.II2}.
\end{align}
We then have, using \eqref{est.II1} and \eqref{est.II2} into \eqref{eq.partial1Xjbar}, that for all $j\neq 1, i=1,\ldots,n$,
\begin{align}
    |\p_{\xi_i}\widebar{X}_j(x,\xi)|\lesssim \frac{1}{|\xi|\la x \ra},
\end{align}
which concludes the proof of \eqref{eq.lem.xbar2} for all $i,j=1,\ldots,n$.

We are now left with the proof of \eqref{eq.lem.xbar1}.
We start again with the case $j=1, i=1,\ldots,n$ and observe that
\begin{align}
    |\p_{x_i}\widebar{X}_1(x,\xi)-\delta_{i1}|\leq \int_0^{+\infty} |\p_{x_i}I(s)|ds+\left|\int_0^{+\infty} \p_{x_i} II(s)ds\right|,\label{eq.partialx1X1bar}
\end{align}
where
\begin{equation}
\begin{split}
    |\p_{x_i}I(s)| &\leq \Bigl| \sum_{k,\ell=2}^n\Bigl(\nabla_X a_{k\ell}(X(s)) \cdot \partial_{x_i} X(s)\Bigr)\Xi_k(s)\Xi_\ell(s)\Bigr|\\
    &\quad +\Bigl|\sum_{k,\ell=2}(a_{k\ell}(X(s))-\delta_{k\ell})\Bigl( \Xi_\ell(s) \partial_{x_i}\Xi_k(s)+\Xi_k(s)\partial_{x_i}\Xi_\ell(s)\Bigr)\Bigr|.
\end{split}
\end{equation}
By \eqref{eq.estderflow1} and \eqref{eq.estderflow2} we have
\begin{align}
    \int_0^\infty |\p_{x}I(s)|ds&\lesssim \int_0^{+\infty} \frac{|\xi|^2}{\la X(s)\ra^3}\left( 1+\frac{|s||\xi|^2}{\la x\ra^2}\right) ds+\int_0^{+\infty} \frac{|\xi|^2}{\la x\ra^2\la X(s)\ra^2}ds\lesssim\frac{1}{\la x \ra^2}.\label{est.pxI}
\end{align}

Arguing as above we also have
\begin{align}
    \int_0^\infty \p_{x_i}II(s) ds=&\sum_{k=1}^nc_k \Bigl[ 
    -\int_0^\infty(\p_{x_i}\Xi_k(s))\int_s^{+\infty}\dt \Xi_k(\sigma)d\sigma\, ds
    -\int_0^\infty\Xi_k(s)\int_s^{+\infty}\p_{x_i}\dt \Xi_k(\sigma)d\sigma\,ds\\
    &-\int_0^\infty(\p_{x_i}\xi_k)\int_s^{+\infty}\dt \Xi_k(\sigma)d\sigma\,ds -\int_0^\infty\xi_k\int_s^{+\infty}\dt \p_{x_i}\Xi_k(\sigma)d\sigma\,ds\\
    &-\bigl(\int_0^{+\infty}\dt \Xi_k(s)ds\Bigr)\Bigl(\int_0^\infty\int_s^{+\infty}\p_{x_i}\dt \Xi_k(\sigma)d\sigma\, ds\Bigr)\\
    &-\bigl(\int_0^{+\infty}\p_{x_i}\dt \Xi_k(s)ds\Bigr)\Bigl(\int_0^\infty\int_s^{+\infty}\dt \Xi_k(\sigma)d\sigma\, ds\Bigr)
    \Bigr]\\
    &= L_1+L_2+L_3+L_4+L_5+L_6.
\end{align}
Recall that
$\la X(s;x,\xi)\ra\approx\la x \ra+s|\xi|^2$ and
\[\int_0^\infty \dt \Xi_k(s)ds\lesssim \frac{|\xi|}{\la x \ra},\quad \forall k=1,\ldots, n.\]
Therefore, by \eqref{eq.estderflow2}

\begin{align}
   \int_s^\infty|\dt \Xi_k(\sigma)|d\sigma\lesssim\frac{|\xi|}{\la X(s)\ra^2}\,\, \forall k=1,\ldots,n &\Longrightarrow |L_1|,|L_3|\lesssim \frac{1}{\la x\ra^2}\label{est.int.x2}  
  \end{align}
  \begin{align}
   \int_s^\infty |\p_{x_i}\dt \Xi_k(\sigma)|d\sigma&\lesssim
   \int_s^\infty \left(\frac{|\xi|^3}{\la X(\sigma)\ra^4} \left( 1+\frac{\sigma|\xi|^2}{\la x\ra^2}\right)
   +\frac{ |\xi|^3}{\la x \ra^2\la X(\sigma)\ra^3}\right) d\sigma&\\
   &\lesssim\frac{|\xi|}{\la x\ra \la X(s)\ra^2}\,\, \forall k=1,\ldots,n\Longrightarrow  |L_2|,|L_4|,|L_5|,|L_6|\lesssim \frac{1}{\la x\ra^2},\label{est.int.derx2}
\end{align}
so we have
\begin{align}
   \left| \int_0^\infty \p_{x_i}II(s)ds\right|\leq |L_1|+|L_2|+\cdots+|L_6|\lesssim \frac{1}{\la x\ra^2} \label{eq.pxXb2},
\end{align}
and finally, plugging \eqref{eq.pxXb2}  and \eqref{est.pxI} into \eqref{eq.partialx1X1bar} we conclude that
\[|\p_{x_i}\widebar{X}_1(x,\xi)-\delta_{i1}|\lesssim \frac{1}{\la x\ra^2},\]
that is \eqref{eq.lem.xbar2} when $j=1$.

In the case  $j\neq 1$ and $i=1,\ldots,n,$ we have
\begin{align}
    |\p_{x_i}\widebar{X}_j(x,\xi)-\delta_{ij}|\leq \int_0^{+\infty} |\p_{x_i}I'(s)|ds+\left|\int_0^{+\infty} \p_{x_i} II'(s)ds\right|.\label{eq.partial1xXjbar}
\end{align}
The first term can be estimated as in the case $j\neq 1$, so we have
\begin{align}
    \int_0^{+\infty} |\p_{x_i}I'(s)|ds\lesssim\frac{1}{\la x\ra^2}.\label{est.II1.2}
\end{align}

To estimate the second term in \eqref{eq.partial1xXjbar} we estimate
\begin{align}
    &\int_0^\infty \p_{x_i} II'(s)ds\\
    &=2 \Bigl[-\int_0^\infty (\p_{x_i}\Xi_j(s))\int_s^\infty \dt \Xi_1(\sigma) d\sigma ds-\int_0^\infty \Xi_j(s)\int_s^\infty \p_{x_i}\dt \Xi_1(\sigma) d\sigma ds\\
    &-\int_0^\infty \int_0^\infty \p_{x_i}\dt \Xi_1(\sigma) d\sigma \int_0^\infty\int_s^\infty \dt \Xi_j(\sigma) d\sigma  ds- \Big(\xi_1+\int_0^\infty\dt \Xi_1(\sigma) d\sigma\Big) \int_0^\infty\int_s^\infty \p_{x_i}\dt \Xi_j(\sigma) d\sigma  ds\Bigr]\\
     &=M_1+M_2+M_3+M_4.
\end{align}
Using the properties in  \eqref{est.int.x2} and \eqref{est.int.derx2}, and \eqref{eq.estderflow2}, we reach
\begin{align}
   \left|\int_0^\infty \p_{x_i} II'(s)ds\right|\lesssim \frac{1}{\la x \ra^2} \label{est.II2.2},
\end{align}
and finally \eqref{est.II1.2} and \eqref{est.II2.2} give
\[|\p_{x_i}\widebar{X}_j(x,\xi)-\delta_{ij}|\lesssim \frac{1}{\la x\ra^2},\]
for all $j\neq1, i=1,\ldots,n$. This shows \eqref{eq.lem.xbar1}.

The estimates \eqref{eq.xbarhighder} and \eqref{eq.xibarhighder} follows as in \cite{CKS}, by using the estimates of the higher order derivatives established in Lemma \ref{lem.estderflow}. 

\endproof

\subsection{Construction of symbols}\label{subsec.bc}
Let now $(x_0,\xi_0) \in \R^{2n}\setminus 0$. 
The goal of this subsection is the construction of a neighborhood $\mathscr{E}_{+}(x_0,\xi_0)$ of the outgoing bicharacteristic $\{ \varphi(s;x_0,\xi_0), 0\leq s<+\infty\}$ that we call \textit{the forward region} (associated with $(x_0,\xi_0)$) and a couple of \emph{forward} symbols $(b(x,\xi),c(x,\xi))$ supported in $\mathscr{E}_{+}(x_0,\xi_0)$ (cf. \cite{CKS}) satisfying suitable conditions. 

Specifically, the couple of symbols we are looking for will belong to the symbol classes $S(m,k;\rho,\delta)$ (see Definition \ref{def.Smkrhodelta}) and must be such that
    \begin{align}\label{propertyPB}
\{a_3,b\}(x,\xi)=c(x,\xi),
    \end{align}
which, up to an integration constant, is equivalent to
\begin{align}\label{propertyPB2}
    b(x,\xi)=\int_{-\infty}^0 c(\varphi(s;x,\xi))ds,
\end{align}
and
     \begin{align}\label{propertyb-c}
b(x,\xi)\geq 0,\quad c(x,\xi)\geq 0,\quad \forall (x,\xi)\in \R^n\times\R^n.
    \end{align}

Observe that, due to \eqref{propertyPB}, we have
\begin{align}
    \mathrm{supp}(c)\subseteq \mathrm{supp}(b).
\end{align}
Moreover, due to \eqref{propertyPB2}, the forward flow maps $\mathrm{supp}(b)$ into itself, that is,
\begin{align}
    \varphi(s;\mathrm{supp}(b))\subseteq \mathrm{supp}(b),\quad  0\leq s<+\infty.
\end{align}
Additionally, we want
    \begin{align}\label{property-abmk}
        \{ a_3(x,\xi),\la x\ra^{-k}\la \xi\ra^{-m}b(x,\xi)\}\geq 0, \quad \forall (x,\xi)\in\mathscr{E}_+^0(x_0,\xi_0)\cap \mathscr{S}_+,
    \end{align}
where  $\mathscr{E}_+^0(x_0,\xi_0):=\{(x,\xi): \ \la x\ra^{-k}\la \xi\ra^{-m}b(x,\xi)\geq 1\}\subseteq \mathscr{E}_+(x_0,\xi_0)$ is a conic neighborhood of $(x_0,\xi_0)$.
\begin{remark}
In what follows, since the scattering variables are only defined locally in certain regions of phase space, we initially construct symbols in a conic neighborhood $V$ of a fixed point $(x,\xi)\in\mathbb{R}^{2n}\setminus 0$, with support contained in a smaller conic neighborhood $V'\subseteq V$. These symbols can then be extended to the whole phase space by multiplying them by a smooth cut-off function $\phi\in C^\infty(\mathbb{R}^{2n})$ satisfying $\phi\equiv 1$ on $V'$ and $\operatorname{supp}\phi\subseteq\widetilde V$, where $\widetilde V$ is a larger neighborhood containing $V$. To simplify the exposition, in what follows we shall define the symbols only in the regions of the phase space under consideration, and consider them as functions on the whole $\mathbb{R}^{2n}$ through the above cut-off procedure whenever necessary.
\end{remark}


Below, for $\xi \in \R^{n}$, we shall use the notation (cf. Lemma \ref{lem.scatteringV} and Remark \ref{rmk.Xi2})
\[
{\xi}^{(2)}:=(3\xi_1^2+|\xi'|^2,2\xi_1\xi_2,\ldots,2\xi_1\xi_n)\in\R^n,
\]
where, recall, $\xi'=(\xi_2,\dots,\xi_n)$. 

\begin{theorem}\label{thm.constrE-}
Let $\delta=1-\rho$, $k \geq 0$ and let $(x_0,\xi_0) \in \R^{2n}\setminus 0$.
Then, if $k>0$, there exist two large constants $M_0,R_0>0$, there exist two conic neighborhoods at the high frequencies of the forward bicharacteristic $\lbrace \varphi(s;x_0,\xi_0); \ s>0 \rbrace$, $\mathscr{E}_{+}^0(x_0,\xi_0)$ and $\mathscr{E}_{+}(x_0,\xi_0)$, with $\mathscr{E}_{+}^0(x_0,\xi_0) \subseteq \mathscr{E}_{+}(x_0,\xi_0)$, and there exist two symbols $b \in S(m,k;\rho,\delta)$ and $c \in S(m+2,k-1;\rho,\delta)$, with $\mathrm{supp}(c)\subseteq \mathrm{supp}(b) \subseteq \mathscr{E}_{+}(x_0,\xi_0)$, satisfying \eqref{propertyPB}, \eqref{propertyb-c}, \eqref{property-abmk} on the whole $\mathscr{E}_+^0(x_0,\xi_0)$,
and
\begin{equation}\label{eq.propk>0} 
b(x,\xi) \gtrsim  \langle \xi \rangle^m, \quad c(x,\xi)\gtrsim \la x\ra^{k-1}\la \xi\ra^{m+2} +\frac{M}{2} \la x\ra^{k-2}\la \xi\ra^{m+2} , \quad \forall (x,\xi) \in \mathscr{E}_{+}^0(x_0,\xi_0),
\end{equation}
where $C$ is a uniform constant such that $ C^{-2}|\xi|^2 \leq |\xi^{(2)}|\leq C^2 |\xi|^2$, and $M,R$ are fixed positive constants such that $M>M_0$ and $R>R_0$.
When $k=0$ the same statement holds with  \eqref{eq.propk>0} replaced by  
\begin{equation}
    b(x,\xi) \gtrsim  \langle \xi \rangle^m, \quad c(x,\xi) \gtrsim \langle x \rangle^{-k'}\langle \xi \rangle^{m+2}+\frac{M}{2}\langle x \rangle^{-k'-2} \langle \xi \rangle^{m+2}, \quad \forall (x,\xi) \in \mathscr{E}_{+}^0(x_1,\xi_1), 
\end{equation}
for all $k'>1$.

Finally, there exists $\widebar{\Xi}_1\in\R^n$ such that, on denoting 
\[
\Gamma_r= \Bigl\lbrace (x,\xi)\in \R^{2n}\setminus 0; \    \Bigl|\frac{\xi\cdot \widebar{\Xi}_1}{|\xi| | \widebar{\Xi}_1|}-1 \Big|<r, \ \ \Big|\frac{x\cdot \xi^{(2)}}{|x| |\xi^{(2)}|}-1  \Big|<\frac{r}{\la x \ra^\rho} \Bigr\rbrace, \quad r>0,
\]
we have that there exists $M>0$ such that 
\begin{equation}\label{lemE.eq1}
\Gamma_{\varepsilon'} \cap \lbrace |\xi|\geq M, |x|>R \rbrace \subseteq \mathscr{E}_+(x_0,\xi_0) \cap \lbrace |\xi|\geq M, |x|>R \rbrace \subseteq \Gamma_{\varepsilon} \cap \lbrace |\xi|\geq M, |x|>R \rbrace,
\end{equation}
and 
\begin{equation}\label{lemE.eq2}
\Gamma_{\varepsilon_0'} \cap \lbrace |\xi|\geq M, |x|>R \rbrace \subseteq \mathscr{E}^0_+(x_0,\xi_0) \cap \lbrace |\xi|\geq M, |x|>R \rbrace \subseteq \Gamma_{\varepsilon_0} \cap \lbrace |\xi|\geq M, |x|>R \rbrace,
\end{equation}
for all $R>R_0$ with $R_0>0$ sufficiently large, and some $0<\varepsilon_0'<\varepsilon_0<\varepsilon'<\varepsilon$.
\end{theorem}

\begin{remark}
The properties \eqref{lemE.eq1} and \eqref{lemE.eq2} are crucial for our purposes, as they encode the asymptotic behavior of the forward sets $\mathscr{E}_+^0(x_0,\xi_0)\subseteq \mathscr{E}_+(x_0,\xi_0)$ as $|x|\to+\infty$. From now on, for $(x_\ast,\xi_\ast) \in \R^{2n}\setminus 0$ we rewrite \eqref{lemE.eq1} and \eqref{lemE.eq2} more concisely as
\begin{align}\label{eq.asympset}
\mathscr{E}_+(x_\ast,\xi_\ast)\simeq 
\left\{ (x,\xi);  \; \Big| \frac{\xi\cdot \widebar{\Xi}_1}{|\xi| \, | \widebar{\Xi}_1|}-1 \Big|<\varepsilon, \; \Big|\frac{x\cdot \xi^{(2)}}{|x| \, |\xi^{(2)}|}-1  \Big|<\frac{\varepsilon}{\langle x \rangle^\rho} \right\}, \,\,\text{as } |x|\to+\infty.
\end{align}
and 
\begin{align}\label{eq.asympset0}
\mathscr{E}^0_+(x_\ast,\xi_\ast)\simeq 
\left\{ (x,\xi);  \; \Big| \frac{\xi\cdot \widebar{\Xi}_1}{|\xi| \, | \widebar{\Xi}_1|}-1 \Big|<\varepsilon_0, \; \Big|\frac{x\cdot \xi^{(2)}}{|x| \, |\xi^{(2)}|}-1  \Big|<\frac{\varepsilon_0}{\langle x \rangle^\rho} \right\}, \,\,\text{as } |x|\to+\infty,
\end{align}
where $0<\varepsilon_0<\varepsilon$.
\end{remark}
    To prove Theorem \ref{thm.constrE-} we need two preliminary lemmas. One, the first, provides the construction of two symbols $b_0\in S(0,0,\rho,1-\rho)$ and $c_0\in S(1,-1,\rho,1-\rho)$ satisfying \eqref{propertyPB} and \eqref{propertyb-c} in  $\mathscr{E}_+(x_0,\xi_0)$ when $(x_0,\xi_0)$ is a forward scattering point, that is  $(x_0,\xi_0)\in \mathscr{S}_+$. As for the second lemma, it uses $b_0,c_0$ to define two symbols $b_+ \in S(m,k;\rho,\delta)$ and $c_{+} \in S(m+2,k-1;\rho,\delta)$  satisfying \eqref{propertyPB}, \eqref{propertyb-c} and \eqref{property-abmk} in  $\mathscr{E}_+(x_0,\xi_0)$, once again when$(x_0,\xi_0)\in \mathscr{S}_+$. These lemmas provide the existence of $b$ and $c$ with the required properties \eqref{propertyPB}, \eqref{propertyb-c} and \eqref{property-abmk} when $(x_0,\xi_0)\in\mathscr{S}_+$. Then, to complete the proof of Theorem \ref{thm.constrE-} we just have to show how to deal with the construction of the symbols in $\mathscr{E}_+(x_0,\xi_0)$ when
    $(x_0,\xi_0)\not\in \mathscr{S}_+$.

\begin{lemma}\label{lem.E.1}

Let $(x_1,\xi_1)\in\mathscr{S}_+$. Then, there is a conic neighborhood $\mathscr{E}_+(x_1,\xi_1)$ of $\{ \varphi(s;x_1,\xi_1),s>0\}$ and there are two symbols $b_0\in S(0,0;\rho,1-\rho)$ and $c_0\in S(2,-1;\rho,1-\rho)$ that satisfy \eqref{propertyPB} and \eqref{propertyb-c}. In addition, there exists $\widebar{\Xi}_1\in\R^n$ such that, $\mathscr{E}_+(x_1,\xi_1)$ satisfies \eqref{eq.asympset}. 
\end{lemma}

\proof
Let $(x_1,\xi_1)\in\mathscr{S}_+$ and let $\widebar{\Xi}_1=\widebar{\Xi}(x_1,\xi_1), \widebar{X}_1=\widebar{X}(x_1,\xi_1)$ be the \emph{forward}  scattering variables defined as in Lemma \ref{lem.scatteringV}. Let also $V(x_1,\xi_1)$ be a conic neighborhood of $\lbrace \varphi(s;x_1,\xi_1), s\geq 0 \rbrace$ where the scattering variables $\widebar{X},\widebar{\Xi}$ are well defined and where, for all $(x,\xi)\in V(x_1,\xi_1)$,
\[
\varphi(s;x,\xi)\approx (\widebar{X}(x,\xi)+s\widebar{\Xi}^{(2)}(x,\xi),\widebar{\Xi}(x,\xi)),\quad \text{as}\,\, s\rightarrow +\infty.
\]
Now, consider $\phi\in C_0^\infty([0,1))$ such that $\phi'(s)\leq 0$, $\phi(s)=1$ when $s$ is close to 0, and $\phi(s)\equiv 0$ when $s$ is close to 1. We then define, for $\rho\in[0,1]$,
\begin{equation}\label{eq.defb1}
\begin{split}
    \tilde{b}_0(x,\xi):=&\phi\left(\frac{1}{\varepsilon}\Big(1-\frac{\widebar{\Xi}_1\cdot\widebar{\Xi}(x,\xi)}{|\widebar{\Xi}_1| |\widebar{\Xi}(x,\xi)|}  \Big) \right)\times \\
    &\times \phi \left( \frac{\la \widebar{X}(x,\xi)-\widebar{X}_1\ra^\rho}{\varepsilon} \Big(1-\frac{(\widebar{X}(x,\xi)-\widebar{X}_1)\cdot \widebar{\Xi}^{(2)}(x,\xi)}{ |\widebar{X}(x,\xi)-\widebar{X}_1||\widebar{\Xi}^{(2)}(x,\xi)|}  \Big) \right),
\end{split}
\end{equation}
where $\varepsilon>0$ is a sufficiently small parameter to be chosen in order to make sure that the support of $\tilde{b}_0$ is included in $V(x_1,\xi_1)$.

The singularities of $\tilde{b}_0$ at $\widebar{\Xi}=0, \widebar{X}=0$ will be removed by using a suitable cutoff function.
Finally, note also that, by Remark \ref{rmk.homscattvar}, $\tilde{b}_0$ is bounded and homogeneous of degree 0 on its domain.

We now define 
\begin{equation}\label{eq.defE-x1xi1}
\mathscr{E}_{+}(x_1,\xi_1):=\mathrm{supp}(\tilde{b}_0) \subseteq \R^{2n}\setminus 0,
\end{equation}
and show that $\mathscr{E}_{+}(x_1,\xi_1)$ is a conic neighborhood of $\lbrace \varphi(s;x_1,\xi_1); \ s>0 \rbrace$ that satisfies \eqref{lemE.eq1}.

First, by \eqref{pro.scatt.var}, for $s>0$ we have 
\[
\begin{split}
\tilde{b}_0(\varphi(s;x_1,\xi_1))&=\phi\left(\frac{1}{\varepsilon}\Big(1-\frac{\widebar{\Xi}_1\cdot\widebar{\Xi}_1}{|\widebar{\Xi}_1| |\widebar{\Xi}_1|}  \Big) \right) \phi \left( \frac{\la s\widebar{\Xi}^{(2)}_1\ra^\rho}{\varepsilon} \Big(1-\frac{(s\widebar{\Xi}^{(2)}_1)\cdot \widebar{\Xi}_1^{(2)}}{ |s\widebar{\Xi}^{(2)}_1||\widebar{\Xi}^{(2)}_1|}  \Big) \right)\\
&=\phi(0)\phi \Bigl( \frac{\la s\widebar{\Xi}^{(2)}_1\ra^\rho}{\varepsilon} \Big(1-\frac{s}{|s|}\Bigr)\Bigr)\\
&=\phi(0)^2=1.
\end{split}
\]
Therefore, since $\mathscr{E}_{+}(x_1,\xi_1)$ is conic due to the homogeneity of $\tilde{b}_0$, it follows immediately that $\mathscr{E}_{+}(x_1,\xi_1)$ is a conic neighborhood of $ \lbrace \varphi (s;x_1,\xi_1); \ s>0 \rbrace$. 
To prove that $\mathscr{E}_{+}(x_1,\xi_1)$ satisfies \eqref{lemE.eq1}, it suffices to use the properties
\begin{equation}
    \widebar{\Xi}(x,\xi)=\xi+O(|\xi|\langle x\rangle^{-2}), \quad  \widebar{X}(x,\xi)=x+O(\langle x\rangle^{-1}),
\end{equation}
which are due to Lemma \ref{lem.scatteringV}, inside the definition of $\mathscr{E}_+(x_1,\xi_1)$. 
Hence, \eqref{lemE.eq1} is a consequence of the definition of $\mathscr{E}_{+}(x_1,\xi_1)=\mathrm{supp}(\tilde{b}_0)$ and of the properties of the support of $\phi$.

We next define 
\[
\tilde{c}_0(x,\xi):=\lbrace a_3, \tilde{b}_0 \rbrace(x,\xi), 
\]
and show that $\tilde{c}_0 \geq 0$.  Since (recall \eqref{pro.scatt.var})
\[
(\widebar{X}(\varphi(s;x,\xi)), \widebar{\Xi}(\varphi(s;x,\xi)) )= (\widebar{X}(x,\xi)+s\widebar{\Xi}^{(2)}(x,\xi),\widebar{\Xi}(x,\xi)),
\]
we have
$\frac{d}{ds}\widebar{\Xi}(\varphi(s;x,\xi)) )=0$ and
\begin{equation}\label{deriv.flow}
\{a_3,f(\widebar{X},\widebar{\Xi})\}=\Bigl(\dfrac{d}{ds}f(\widebar{X}(\varphi(s)),\widebar{\Xi}(\varphi(s)))\Bigr)\Bigl|_{s=0}=\widebar{\Xi}^{(2)}\cdot \p_{\widebar{X}}f(\widebar{X},\widebar{\Xi}).
\end{equation}
Using this property we get
\begin{align}
    \{a_3,\tilde{b}_0\}(x,\xi)&=\phi\left(\frac{1}{\varepsilon}\Big(1-\frac{\widebar{\Xi}_1\cdot\widebar{\Xi}(x,\xi)}{|\widebar{\Xi}_1| |\widebar{\Xi}(x,\xi)|}  \Big) \right) \\
    &\times \p_s\phi \left( \frac{\la \widebar{X}(x,\xi)-\widebar{X}_1\ra^\rho}{\varepsilon} \Big(1-\frac{(\widebar{X}(x,\xi)-\widebar{X}_1)\cdot \widebar{\Xi}^{(2)}(x,\xi)}{ |\widebar{X}(x,\xi)-\widebar{X}_1||\widebar{\Xi}^{(2)}(x,\xi)|}  \Big) \right)\\
    &\times \p_{\widebar{X}}\left( \frac{\la \widebar{X}(x,\xi)-\widebar{X}_1\ra^\rho}{\varepsilon} \Big(1-\frac{(\widebar{X}(x,\xi)-\widebar{X}_1)\cdot \widebar{\Xi}^{(2)}(x,\xi)}{ |\widebar{X}(x,\xi)-\widebar{X}_1||\widebar{\Xi}^{(2)}(x,\xi)|}  \Big) \right)\cdot \widebar{\Xi}^{(2)}\\
    &=\phi\left(\frac{1}{\varepsilon}\Big(1-\frac{\widebar{\Xi}_1\cdot\widebar{\Xi}(x,\xi)}{|\widebar{\Xi}_1| |\widebar{\Xi}(x,\xi)|}  \Big) \right)\\
    &\times \p_s\phi \left( \frac{\la \widebar{X}(x,\xi)-\widebar{X}_1\ra^\rho}{\varepsilon} \Big(1-\frac{(\widebar{X}(x,\xi)-\widebar{X}_1)\cdot \widebar{\Xi}^{(2)}(x,\xi)}{ |\widebar{X}(x,\xi)-\widebar{X}_1||\widebar{\Xi}^{(2)}(x,\xi)|}  \Big) \right)\\
    &\quad \times \Big\{
\frac{\la \widebar{X}(x,\xi)-\widebar{X}_1\ra^\rho}{\varepsilon} \Big(  -\frac{\widebar{\Xi}^{(2)}\cdot\widebar{\Xi}^{(2)}}{|\widebar{X}(x,\xi)-\widebar{X}_1||\widebar{\Xi}^{(2)}(x,\xi)|}+\frac{((\widebar{X}-\widebar{X}_1)\cdot \widebar{\Xi}^{(2)})^2}{|\widebar{X}(x,\xi)-\widebar{X}_1|^3|\widebar{\Xi}^{(2)}|} +\\
&\quad \quad + \frac{\rho}{\varepsilon} \la \widebar{X}-\widebar{X}_1\ra ^{\rho-2}\widebar{\Xi}^{(2)}\cdot(\widebar{X}-\widebar{X}_1)\Big(1-\frac{(\widebar{X}(x,\xi)-\widebar{X}_1)\cdot \widebar{\Xi}^{(2)}(x,\xi)}{ |\widebar{X}(x,\xi)-\widebar{X}_1||\widebar{\Xi}^{(2)}(x,\xi)|}
\Big)
    \Big\}.
\end{align}
Hence, to show that $\tilde{c}_0 \geq 0$, since $\phi'\leq 0$, $0\leq\phi\leq 1$, it suffices to prove that the quantity inside the curly braces is nonpositive. To this end, we rewrite such a quantity as
\begin{align}
    &\Big\{
\frac{\la \widebar{X}(x,\xi)-\widebar{X}_1\ra^\rho}{\varepsilon} \Big(  -\frac{\widebar{\Xi}^{(2)}\cdot\widebar{\Xi}^{(2)}}{|\widebar{X}(x,\xi)-\widebar{X}_1||\widebar{\Xi}^{(2)}(x,\xi)|}+\frac{((\widebar{X}-\widebar{X}_1)\cdot \widebar{\Xi}^{(2)})^2}{|\widebar{X}(x,\xi)-\widebar{X}_1|^3|\widebar{\Xi}^{(2)}|} +\\
&\quad \quad + \frac{\rho}{\varepsilon} \la \widebar{X}-\widebar{X}_1\ra ^{\rho-2}\widebar{\Xi}^{(2)}\cdot(\widebar{X}-\widebar{X}_1)\Big(1-\frac{(\widebar{X}(x,\xi)-\widebar{X}_1)\cdot \widebar{\Xi}^{(2)}(x,\xi)}{ |\widebar{X}(x,\xi)-\widebar{X}_1||\widebar{\Xi}^{(2)}(x,\xi)|}\Big)
    \Big\}\\
    &=-\frac{\la \widebar{X}(x,\xi)-\widebar{X}_1\ra^{\rho-2}}{\varepsilon} \Big(1-\frac{(\widebar{X}(x,\xi)-\widebar{X}_1)\cdot \widebar{\Xi}^{(2)}(x,\xi)}{ |\widebar{X}(x,\xi)-\widebar{X}_1||\widebar{\Xi}^{(2)}(x,\xi)|}
\Big)\\
&\times \Big( \frac{\la \widebar{X}(x,\xi)-\widebar{X}_1\ra^2}{|\widebar{X}(x,\xi)-\widebar{X}_1|^2} \big(|\widebar{\Xi}^{(2)}||\widebar{X}-\widebar{X}_1|+\widebar{\Xi}^{(2)}\cdot(\widebar{X}-\widebar{X}_1)\big) -\rho \widebar{\Xi}^{(2)}\cdot(\widebar{X}-\widebar{X}_1)\Big)
\end{align}
where, since $|\widebar{\Xi}^{(2)}||\widebar{X}-\widebar{X}_1|\geq\widebar{\Xi}^{(2)}\cdot(\widebar{X}-\widebar{X}_1)$ and $\la \widebar{X}(x,\xi)-\widebar{X}_1\ra^{(2)}\geq |\widebar{X}(x,\xi)-\widebar{X}_1|^2$, 
\begin{align}
    &\frac{\la \widebar{X}(x,\xi)-\widebar{X}_1\ra^2}{|\widebar{X}(x,\xi)-\widebar{X}_1|^2} \big(|\widebar{\Xi}^{(2)}||\widebar{X}-\widebar{X}_1|+\widebar{\Xi}^{(2)}\cdot(\widebar{X}-\widebar{X}_1)\big) -\rho \widebar{\Xi}^{(2)}\cdot(\widebar{X}-\widebar{X}_1)\\
    &\geq |\widebar{\Xi}^{(2)}||\widebar{X}-\widebar{X}_1|+(1-\rho)\widebar{\Xi}^{(2)}\cdot(\widebar{X}-\widebar{X}_1)\\
    &\geq(1+(1-\rho)(1-\varepsilon))|\widebar{\Xi}^{(2)}||\widebar{X}-\widebar{X}_1|\geq 0
\end{align}
which proves that $\tilde{c}_0=\{a_3,\tilde{b}_0\}\geq 0$.

We now show that, for $\widebar{X}$ away from $\widebar{X}_1$ and $\widebar{\Xi}$ away from 0, $\tilde{c}_0(x,\xi)$ satisfies
$$\tilde{c}_0(x,\xi)\lesssim \la \widebar{X}-\widebar{X}_1\ra^{-1}\la \xi\ra^2\lesssim \langle x \rangle^{-1}\la \xi \ra ^2.$$
Since $\tilde{b}_0$ and $\tilde{c}_0$ are not smooth at $\widebar{\Xi}^{(2)}= 0$ and at  $\widebar{X}=\widebar{X}_1$, once we prove the bound for $\tilde{c}_0$ we will define $c_0\in S(2,-1;\rho,1-\rho)$ by  means of $\tilde{c}_0$ and a suitable cutoff function.

To see the bound on $\tilde{c}_0$ away from the singular points, we assume that $|\widebar{X}-\widebar{X}_1|>\delta/2$ for some $\delta>0$ sufficiently small and $\widebar{\Xi}\neq 0$, so we have
\begin{align}
    \tilde{c}_0(x,\xi)\lesssim \frac{|\widebar{\Xi}^{(2)}|}{|\widebar{X}-\widebar{X}_1|}+\frac{\rho|\widebar{\Xi}^{(2)}|}{\la \widebar{X}-\widebar{X}_1\ra}\lesssim \frac{\la \widebar{\Xi}\ra^2}{\la \widebar{X}-\widebar{X}_1\ra}\lesssim\la \xi\ra^2 \la x \ra^{-1},
\end{align}
where we have used the fact that
\begin{equation}\label{eq.usCx1}
\la \widebar{X}-\widebar{X}_1\ra\geq C(x_1)\la \widebar{X}(x,\xi)\ra\geq C'(x_1)\la x\ra,
\end{equation}
where the last inequality is due to Proposition \ref{prop.App.2}.

Next, we define $b_0$ and $c_0$ as 
\begin{align}
    &b_0(x,\xi):=\psi(|\widebar{\Xi}^{(2)}(x,\xi)|) \psi\Big( \frac{1}{\varepsilon'}|\widebar{X}(x,\xi)-\widebar{X}_1|\Big)\tilde{b}_0(x,\xi),\label{def.b0}\\
    & c_0(x,\xi):=\lbrace a_3,b_0 \rbrace (x,\xi),
\end{align}
where $\psi\in C^\infty(\R)$ is a smooth cutoff function supported in $[1/2,+\infty)$ such that $\psi'\geq 0$ and $\psi(s)=1$ for $s>1$, and where $\varepsilon'$ is a sufficiently small constant. The cutoff function cuts the singularities of $\tilde{b}_0$ and makes $b_0$ smooth and nonnegative. Note, in particular, that $b_0(x_1,\xi_1)=0$ due to $\psi\Big( \frac{1}{\varepsilon'}|\widebar{X}-\widebar{X}_1|\Big)$. 

Moreover, the forward flow maps $(x_1,\xi_1)$ into a point $\varphi(s;x_1,\xi_1)$ such that 
\[
b_0(\varphi(s;x_1,\xi_1))>0, \quad  \forall s>\varepsilon'/|\Xi_1^{(2)}|.
\]
To see this we observe that by \eqref{pro.scatt.var}
\begin{align}
    (\widebar{X}(\varphi(s;x_1,\xi_1)),\widebar{\Xi}(\varphi(s;x_1,\xi_1)))=\big(\widebar{X}(x_1,\xi_1)+s\widebar{\Xi}^{(2)}(x_1,\xi_1), \widebar{\Xi}(x_1,\xi_1)\big)=(\widebar{X}_1+s\widebar{\Xi}_1^{(2)},\widebar{\Xi}_1),
\end{align}
therefore,
\begin{align}
  \frac{1}{\varepsilon'} |\widebar{X}(\varphi(s;x_1,\xi_1))-\widebar{X}_1|=\frac{1}{\varepsilon'}|\widebar{X}_1+s\widebar{\Xi}_1^2-\widebar{X}_1|\geq\frac{|s||\widebar{\Xi}_1^2|}{\varepsilon'}>1,\quad \forall s>\varepsilon'/|\widebar{\Xi}_1^2|,
\end{align}
meaning that $\psi\Big( \frac{1}{\varepsilon'}|\widebar{X}-\widebar{X}_1|\Big)>0$ and $b_0(\varphi(s;x_1,\xi_1))>0$. Moreover, thanks to \eqref{deriv.flow}, it is easy to check (cf. \cite{CKS}) that the nonnegativity of $\{a_3,b_0\}=c_0$ is also preserved.

What is left is to see is that $b_0$ and $c_0$ belong to the suitable symbol classes. To investigate the symbol properties of $b_0$, it suffices to control each term of $\tilde{b}_0$. After application of multivariate Faa di Bruno's formula (see Theorem \ref{theo.fdb}) the argument is the same as that in \cite{CKS}, to which we refer for further details. In our case, by Theorem \ref{theo.dersv}, we have that for each $\alpha,\beta \in \mathbb{N}_0^n$ there exists $C_{\alpha,\beta}>0$ such that 
\[
\Bigl|\p_x^\beta \p_\xi^\alpha \Bigl(\frac{\widebar{\Xi}_1\cdot\widebar{\Xi}(x,\xi)}{|\widebar{\Xi}_1| |\widebar{\Xi}(x,\xi)|}\Bigr)\Bigr|\leq C_{\alpha,\beta} \la x \ra^{-|\beta|}|\xi|^{-|\alpha|}, 
\]
\[
\Bigl|\p_x^\beta \p_\xi^\alpha \Bigl(\frac{(\widebar{X}(x,\xi)-\widebar{X}_1)\cdot \widebar{\Xi}^{(2)}(x,\xi)}{ |\widebar{X}(x,\xi)-\widebar{X}_1||\widebar{\Xi}^{(2)}(x,\xi)|}\Bigr)\Bigr|\leq C_{\alpha,\beta} \la x \ra^{-|\beta|}|\xi|^{-|\alpha|},
\]
and 
\[
|\p_x^\beta \p_\xi^\alpha\langle \widebar{X}(x,\xi)-\widebar{X}_1 \rangle^\rho | \leq C_{\alpha,\beta}\langle x \rangle^{\rho-|\beta|}|\xi|^{-|\alpha|},
\]
where $|\widebar X-X_1|>\varepsilon'/2$. We then examine the second term of $\tilde{b}_0$, for the first term the same argument applies. To this end, we denote 
\[
F(x,\xi)
=
\frac{\langle\bar X(x,\xi)-\bar X_1\rangle^\rho}{\varepsilon}
\left(
1-
\frac{(\bar X(x,\xi)-\bar X_1)\cdot\bar\Xi^{(2)}(x,\xi)}
{|\bar X(x,\xi)-\bar X_1|\,|\bar\Xi^{(2)}(x,\xi)|}
\right).
\]
By the previous estimates,
\[
\left|
\partial_x^\beta\partial_\xi^\alpha F(x,\xi)
\right|
\leq
C_{\alpha\beta}
\langle x\rangle^{\rho-|\beta|}
|\xi|^{-|\alpha|}.
\]
For $|\alpha|+|\beta|\geq1$, the iterated chain rule expresses
$\partial_x^\beta\partial_\xi^\alpha(\phi(F))$
as a finite sum of terms of the form $\phi^{(j)}(F)
\prod_{\ell=1}^{j}
\partial_x^{\beta_\ell}\partial_\xi^{\alpha_\ell}F$, where $\sum_{\ell=1}^{j}\alpha_\ell=\alpha$, $
\sum_{\ell=1}^{j}\beta_\ell=\beta$, $|\alpha_\ell|+|\beta_\ell|\geq1.$
In particular, $1\leq j\leq|\alpha|+|\beta|$.
Since $\phi \in C_c^\infty([0,+\infty))$, each such term
satisfies
\[
\begin{aligned}
\left|
\phi^{(j)}(F)
\prod_{\ell=1}^{j}
\partial_x^{\beta_\ell}\partial_\xi^{\alpha_\ell}F
\right|
&\leq
C_{\alpha\beta}
\prod_{\ell=1}^{j}
\left(
\langle x\rangle^{\rho-|\beta_\ell|}
|\xi|^{-|\alpha_\ell|}
\right)\\
&=
C_{\alpha\beta}
\langle x\rangle^{j\rho-|\beta|}
|\xi|^{-|\alpha|}\\
&\leq
C_{\alpha\beta}
\langle x\rangle^{(\rho-1)|\beta|+\rho|\alpha|}
|\xi|^{-|\alpha|},
\end{aligned}
\]
where the last inequality follows from $\rho\geq0$ and
$j\leq|\alpha|+|\beta|$. Summing over all terms yields
\[
\left|
\partial_x^\beta\partial_\xi^\alpha
\bigl(\phi(F(x,\xi))\bigr)
\right|
\leq
C_{\alpha\beta}
\langle x\rangle^{(\rho-1)|\beta|+\rho|\alpha|}
|\xi|^{-|\alpha|}.
\]
The same estimate holds for $\alpha=\beta=0$ by the boundedness
of $\phi$.

Hence, as in \cite{CKS}, by Theorem \ref{theo.fdb} we obtain that $b_0\in S(0,0;\rho,1-\rho)$ and $c_0\in S(2,-1;\rho,1-\rho)$, where the latter follows from the fact that in our case $a_3$ is a standard symbol of order $3$.
\endproof

\begin{remark}\label{rmk.suppb0}
    Recall that, due to Lemma \ref{lem.scatteringV}, $|\widebar{\Xi}(x,\xi)| \approx |\xi|$, therefore, by definition of $\psi$ (see \eqref{def.b0}) we have $|\xi|\geq c$ for some $c>0$, in $\mathrm{supp}(b_0)$. Moreover, since the (forward) scattering variables $\widebar{X}, \widebar{\Xi}$ are defined on a small conic neighborhood of $\mathscr{S}_{+}$, we also that $|x| \geq R$, for a large $R>0$, on $\mathrm{supp}(b_0)$. Summarizing, we have that
    \[|
    \xi|\geq c, \quad |x|>R,\quad \text{on } \mathrm{supp}(b_0).
    \]
    Note also that, if $(x,\xi)$ is such that $b_0(x,\xi)\geq 0$, then $c_0(x,\xi)\geq 0$.

\end{remark}

Building on the construction in Lemma \ref{lem.E.1}, we can find a pair of symbols $(b,c)$ of order $m$ and $m+2$ respectively, with $m\in\R$, satisfying suitable growth conditions. The existence of such symbols is given in the following lemma.

\begin{lemma}\label{lem.E.2}
Let $\delta=1-\rho$, $k \geq 0$ and $(x_1,\xi_1) \in \mathscr{S}_{+}$. Then, if $k>0$, there exist two neighborhoods $\mathscr{E}_{+}^0(x_1,\xi_1)$ and $\mathscr{E}_{+}(x_1,\xi_1)$, conic at high frequencies, of the forward bicharacteristic $\lbrace \varphi(s;x_1,\xi_1); \ s >s_\ast \rbrace$, for some $s_\ast>0$, with $\mathscr{E}_{+}^0(x_1,\xi_1) \subseteq \mathscr{E}_{+}(x_1,\xi_1)$, satisfying \eqref{eq.asympset}, and there exist two symbols $b_{+} \in S(m,k;\rho,\delta)$ and $c_{+} \in S(m+2,k-1;\rho,\delta)$, with $\mathrm{supp}(c_+)\subseteq \mathrm{supp}(b_+) \subseteq \mathscr{E}_{+}(x_1,\xi_1)$, satisfying \eqref{propertyPB}, \eqref{propertyb-c}, \eqref{property-abmk}
and 
\begin{equation}\label{eq.propk>0bis}
b_+(x,\xi) \gtrsim  \langle x \rangle^k \langle \xi \rangle^m, \quad c_+(x,\xi) \gtrsim \langle x \rangle^{k-1}\langle \xi \rangle^{m+2}+\frac{M}{2}\langle x \rangle^{k-2} \langle \xi \rangle^{m+2}, \quad \forall (x,\xi) \in \mathscr{E}_{+}^0(x_1,\xi_1),
\end{equation}
where $M>0$ is a universal constant that can be taken arbitrarily large. 

When $k=0$ the same statement holds with  \eqref{eq.propk>0bis} replaced by 
\begin{equation}
    b_+(x,\xi) \gtrsim  \langle \xi \rangle^m, \quad c_+(x,\xi) \gtrsim \langle x \rangle^{-k'}\langle \xi \rangle^{m+2}+\frac{M}{2}\langle x \rangle^{-k'-2} \langle \xi \rangle^{m+2}, \quad \forall (x,\xi) \in \mathscr{E}_{+}^0(x_1,\xi_1), 
\end{equation}
for all $k'>1$.
\end{lemma}
\begin{proof}
Let $k>0$. In this case on denoting $b_0,c_0$ the symbols constructed in Lemma \ref{lem.E.1}. First note that, for all $(x,\xi) \in \mathrm{supp}(b_0)$, we have 
\begin{equation}\label{eq.usneg0}
    (\widebar{X}(x,\xi)- \widebar{X}_1(x,\xi)) \cdot \widebar{\Xi}^{(2)}(x,\xi)\geq (1-\varepsilon)  |\widebar{X}(x,\xi)-\widebar{X}_1(x,\xi)| |\widebar{\Xi}^{(2)}(x,\xi)|, \quad 0<\varepsilon<1/2
\end{equation}
Furthermore, by Proposition \ref{prop.App.2},
    \[
    \quad \la \widebar{X}(x,\xi)\ra \approx \la x\ra , \quad |\widebar{\Xi}^{(2)}(x,\xi)|\approx |\xi|^2. 
    \]
    Define
    \[
    b_+(x,\xi):=e^{-M/\langle \widebar{X}-\widebar{X_1}\rangle}\langle \widebar{X}(x,\xi) \rangle^k \langle \widebar{\Xi}(x,\xi) \rangle^m b_0(x,\xi), \quad c_+(x,\xi):=\lbrace a_3,b_+ \rbrace(x,\xi), 
    \]
    where $M>0$ is a suitable constants to be chosen later.
    
    The positivity of $b_+$ follows by the positivity of $b_0$. Next, we check the positivity of $c_+$.
    We have  
    \[
    \begin{split}
    \lbrace a_3,b_+ \rbrace &=e^{-M/\langle \widebar{X}-\widebar{X_1}\rangle} \langle \widebar{\Xi} \rangle^m b_0\lbrace a_3,\langle \widebar{X} \rangle^k \rbrace +e^{-M/\langle \widebar{X}-\widebar{X_1}\rangle}\langle \widebar{X} \rangle^kb_0\lbrace a_3,\langle \widebar{\Xi}\rangle^m \rbrace\\
    &\quad +e^{-M/\langle \widebar{X}-\widebar{X_1}\rangle} \langle \widebar{X} \rangle^k \langle \widebar{\Xi} \rangle^m \lbrace a_3,b_0 \rbrace+b_0\langle \widebar{X} \rangle^k \langle \widebar{\Xi}\rangle^m\lbrace a_3, e^{-M/\langle \widebar{X}-\widebar{X_1}\rangle} \rbrace
    \end{split}
    \]
    Now, by \eqref{pro.scatt.var}, we get
    \[
    \begin{split}
    \lbrace a_3, \langle \widebar{X} \rangle^k \rbrace (x,\xi)&=\frac{d}{dt} \Bigl(\langle \widebar{X}(\varphi(t;x,\xi))\rangle^k\Bigr)\Bigl|_{t=0}\\
    &=\frac{d}{dt} \Bigl(\langle \widebar{X}(x,\xi)+t\widebar{\Xi}^{(2)}(x,\xi)\rangle^k\Bigr)\Bigl|_{t=0}\\
    &=k\la \widebar{X}(x,\xi) \ra^{k-2} \widebar X(x,\xi) \cdot \widebar{\Xi}^{(2)}(x,\xi)
    \end{split}
    \]
    and 
    \[
    \begin{split}
    \lbrace a_3, \langle \widebar{\Xi} \rangle^m \rbrace (x,\xi)&=\frac{d}{dt} \Bigl(\langle \widebar{\Xi}(\varphi(t;x,\xi))\rangle^m\Bigr)\Bigl|_{t=0}\\
    &=\frac{d}{dt} \Bigl(\langle \widebar{\Xi}(x,\xi)\rangle^m\Bigr)\Bigl|_{t=0}\\&=0,
    \end{split}
    \]
    and, moreover, 
    \[
    \lbrace a_3, -M/\langle \widebar{X}-\widebar{X_1}\rangle \rbrace=M\frac{(\widebar{X}-\widebar{X}_1)\cdot \widebar{\Xi}^{(2)}}{\langle \widebar{X}-\widebar{X}_1 \rangle^3}.
    \]
    Therefore, we have 
    \begin{equation}\label{eq.clowbound}
    \begin{split}
    c_+&=\lbrace a_3,b_+ \rbrace \\
    &=e^{-M/\langle \widebar{X}-\widebar{X_1}\rangle}\langle \widebar{X}\rangle^k \langle \widebar{\Xi} \rangle^m c_0+b_+\Bigl( k \frac{\widebar{X}\cdot \widebar{\Xi}^{(2)}}{\langle \widebar{X}\rangle^2}+M \frac{(\widebar{X}-\widebar{X}_1)\cdot \widebar{\Xi}^{(2)}}{\langle \widebar{X}-\widebar{X}_1\rangle^3} \Bigr).
    \end{split}
    \end{equation}
Now, by \eqref{eq.usneg0}, we also get
\begin{equation}\label{eq.usneg}
k \frac{\widebar{X}\cdot \widebar{\Xi}^{(2)}}{\langle \widebar{X}\rangle^2}\geq \nu k|\xi|^2\langle x \rangle^{-1}-C|\xi|^2\langle x \rangle^{-2},
\end{equation}
for some constants $\nu,C>0$.
Therefore, from \eqref{eq.usneg0}, $|\xi|\approx |\widebar{\Xi}(x,\xi)|, \la \widebar{X}-\widebar{X}_1\ra \gtrsim \la x\ra $ and \eqref{eq.usneg}, by choosing $M>0$ sufficiently large we obtain, with $c'\geq 0$ universal constant, 
\begin{equation}\label{eq.boundc}
    c_+=\lbrace a_3,b_+ \rbrace \gtrsim b_+|\xi|^2 \Big(\nu k\langle x \rangle^{-1}+(c'M-C)\la x\ra^{-2}\Big) \geq 0.
\end{equation}
We now define the sets,
     \[
     \mathscr{E}_+(x_1,\xi_1):=\supp(b_+), \ \ \mathscr{E}_{+}^0(x_1,\xi_1):=\Bigl\lbrace (x,\xi) \in \R^{2n}; \  b_+(x,\xi) \geq c\langle x \rangle ^k \langle \xi \rangle^m \Bigr\rbrace,
     \]
     where $c>0$ is a suitable small constant. 
     They satisfy \eqref{lemE.eq1}, due to the support properties of $b_1$. By definition, we trivially have 
     \[
     b_+(x,\xi) \gtrsim \la x \ra^k \la \xi \ra^m, \quad \forall (x,\xi) \in \mathscr{E}^0_{+}(x_1,\xi_1),
     \]
     and by \eqref{eq.boundc}
     \[
     c_+(x,\xi) \gtrsim \langle x \rangle^{k-1} \langle \xi \rangle^{m+2}+C_M\langle x \rangle^{k-2} \langle \xi \rangle^{m+2}, \quad \forall (x,\xi) \in \mathscr{E}^0_{+}(x_1,\xi_1),
     \]
     since by Remark \ref{rmk.suppb0} $\mathrm{supp}(b_+)\subseteq \{(x,\xi)\in\R^{2n}; |\xi|\geq c'\}$ for some positive $c'$, where $C_M$ can be taken arbitrarily large by choosing $M>0$ large.
     
      We finally prove that 
     \begin{equation}\label{eq.posamkb}
     \lbrace a_3,\la x \ra^{-k}\la \xi \ra^{-m}b_+ \rbrace \geq 0. 
     \end{equation}
     For this purpose, note that 
     \[
     \lbrace a_3, \langle x \rangle^{-k} \rbrace= -k \langle x \rangle^{-k-2} \nabla_\xi a_3 \cdot x, \quad  \lbrace a_3, \langle \xi \rangle^{-m} \rbrace= m \langle \xi \rangle^{-m-2} \nabla_x a_3 \cdot \xi, 
     \]
    and so by using also the expression \eqref{eq.clowbound}
    we obtain
\begin{align}
&\{a_3,\langle x\rangle^{-k}\langle\xi\rangle^{-m}b_+\}
\notag\\
&=
e^{-M/\langle\widebar{X}-\widebar{X}_1\rangle}
\frac{\langle\widebar{X}\rangle^k\langle\widebar{\Xi}\rangle^m}
     {\langle x\rangle^k\langle\xi\rangle^m}
\Bigg(
c_0+b_0\Bigg(
M\frac{(\widebar{X}-\widebar{X}_1)\cdot\widebar{\Xi}^{(2)}}
        {\langle\widebar{X}-\widebar{X}_1\rangle^3}
\notag\\
&\hspace{48mm}
+k\frac{\widebar{X}\cdot\widebar{\Xi}^{(2)}}{\langle\widebar{X}\rangle^2}
-k\frac{x\cdot\nabla_\xi a_3}{\langle x\rangle^2}
+m\frac{\xi\cdot\nabla_x a_3}{\langle\xi\rangle^2}
\Bigg)
\Bigg).
\label{eq:normalized-bracket-exp-correction}
\end{align}
Since $b_0,c_0\geq0$, to get \eqref{eq.posamkb}, it suffices to prove that the term multiplied by $b_0$ is positive on $\mathrm{supp}(b_0)$.
To prove this desired property we recall that, by definition of scattering variables
\[
\widebar{X}-x=O(\langle x\rangle^{-1}),
\qquad
\widebar{\Xi}^{(2)}-\nabla_\xi a_3
=O(|\xi|^2\langle x\rangle^{-2}),
\]
and so
\begin{align*}
\widebar{X}\cdot\widebar{\Xi}^{(2)}-x\cdot\nabla_\xi a_3
&=(\widebar{X}-x)\cdot\widebar{\Xi}^{(2)}
  +x\cdot(\widebar{\Xi}^{(2)}-\nabla_\xi a_3)\\
&=O(|\xi|^2\langle x\rangle^{-1}).
\end{align*}
Moreover, using $\langle\widebar{X}\rangle\approx\langle x\rangle$, we have
\[
\langle\widebar{X}\rangle^{-2}-\langle x\rangle^{-2}
=
\frac{|x|^2-|\widebar{X}|^2}
     {\langle\widebar{X}\rangle^2\langle x\rangle^2}
=O(\langle x\rangle^{-4}).
\]
It follows that
\[
\left|
\frac{\widebar{X}\cdot\widebar{\Xi}^{(2)}}{\langle\widebar{X}\rangle^2}
-\frac{x\cdot\nabla_\xi a_3}{\langle x\rangle^2}
\right|
\lesssim |\xi|^2\langle x\rangle^{-3}.
\]
The decay assumption on the coefficients also yields
\[
\left|
m\frac{\xi\cdot\nabla_x a_3}{\langle\xi\rangle^2}
\right|
\lesssim \langle x\rangle^{-3}
       \frac{|\xi|^4}{\langle\xi\rangle^2}
\lesssim |\xi|^2\langle x\rangle^{-3}.
\]
On the other hand, the support properties of $b_0$ 
\[
\frac{(\widebar{X}-\widebar{X}_1)\cdot\widebar{\Xi}^{(2)}}
       {\langle\widebar{X}-\widebar{X}_1\rangle^3}
\geq
c_1|\xi|^2\langle\widebar{X}-\widebar{X}_1\rangle^{-2}
\]
for some $c_1>0$. Thus, since $\widebar{X}_1$ is fixed, and $\langle\widebar{X}-\widebar{X}_1\rangle\approx\langle x\rangle$ (recall $\widebar{X}_1$ is fixed) we have that the expression in parentheses in
\eqref{eq:normalized-bracket-exp-correction} is bounded from below by
\[
c_1 M|\xi|^2\langle x\rangle^{-2}
-c_2|\xi|^2\langle x\rangle^{-3},
\]
for some universal constants $c_1,c_2>0$. Then, by choosing $M>0$ sufficiently large, say $M>M_0$ for a certain $M_0$ large, we obtain \eqref{eq.posamkb} which yields \eqref{property-abmk}. Finally, this will give \eqref{property-abmk} at each point of $\mathscr{E}_+^0(x_0,\xi_0)$.

To conclude that the symbols $b_+,c_+$ have the desired order we observe that
\begin{align}
  & e^{-M/\langle \widebar{X}-\widebar{X_1}\rangle}\in S(0,0;0,1),\quad \langle \widebar{\Xi}(x,\xi) \rangle^m   \in S(m,0;0,1),\quad \langle \widebar{X}(x,\xi) \rangle^k   \in S(0,k;0,1),
\end{align}
therefore, by using the Leibniz rule in the expressions of $b_+$ and $c_+$ (see also \eqref{propertyComposition}), we get that $b_+\in S(m,k;\rho,\delta)$ and $c_+\in S(m+2,k-1;\rho,\delta)$.

When $k=0$, if $k'>1$ we set $\sigma=\min \lbrace 2,k'-1 \rbrace$, we define
 \[
    b_+(x,\xi):=e^{-M/\langle \widebar{X}-\widebar{X_1}\rangle^\sigma} \langle \widebar{\Xi}(x,\xi) \rangle^m b_0(x,\xi), \quad c_+(x,\xi):=\lbrace a_3,b_+ \rbrace(x,\xi), 
    \]
and we argue as before.

\end{proof}

At this point, to handle the case in which $(x_0,\xi_0) \notin \mathscr{S}_{+}$ the following additional lemmas are required.
\begin{lemma}\label{lem.x1xi1}
    If $(x_0,\xi_0) \in \R^{2n}\setminus 0$, 
    then there exists $s_1>0$ such that $(x_1,\xi_1)=\varphi(s_1;x_0,\xi_0) \in \mathscr{S}_{+}$.
\end{lemma}
\begin{proof}
    Since $(x_0,\xi_0)$ is forward nontrapped there exists $s_1>0$ such that $|X(s_1;x_0,\xi_0)|>R_0$  and $d/ds(|X(s_1;x_0,\xi_0)|^2)>0$. Hence, due to the identity 
    \[
    \frac{d}{ds}|X(s)|^2=2 X(s)\cdot \nabla_\Xi a_3(X(s),\Xi(s)),
    \]
    we obtain that $(x_1,\xi_1):=\varphi(s_1;x_0,\xi_0) \in \mathscr{S}_+$.
\end{proof}

We are now ready to prove Theorem \ref{thm.constrE-}.
\begin{proof}[Proof of Theorem \ref{thm.constrE-}]
We prove the theorem in the case of $k>0$. The case $k=0$ follows analogously (cf. \cite{CKS}).

Let $(x_0,\xi_0) \in \R^{2n}\setminus 0$. 
If $(x_0,\xi_0) \in \mathscr{S}_{+}$, the proof follows directly by Lemma \ref{lem.E.1} and Lemma \ref{lem.E.2} after the introduction of an auxiliary symbol $b_\ast$ (see the construction below) which covers the region of the bicharacteristic for $s$ between $s^\ast$ and $0$. Hence we consider the case $(x_0,\xi_0) \notin \mathscr{S}_{+}$. 

In this case, we choose $V(x_0,\xi_0)$ a conic neighborhood of $(x_0,\xi_0)$ and $\chi_0 \in C^\infty(\R^{2n}\setminus 0;[0,1])$ such that $\mathrm{supp}(\chi_0) \subset V(x_0,\xi_0)$ and
$\pi_x(\operatorname{supp}\chi_0)$ is compact and $\chi_0(x,\lambda\xi)=\chi_0(x,\xi)$
for all $\lambda>0$.
Moreover, we  require $\chi_0=1$ on a smaller conic neighborhood $V'(x_0,\xi_0)$
of $(x_0,\xi_0)$. Such a function can be obtained as the product of a smooth,
compactly supported function of $x$, equal to $1$ near $x_0$,
and a smooth function of $\xi/|\xi|$, equal to $1$ near
$\xi_0/|\xi_0|$, with sufficiently small supports.

 Hence, we  define 
\[
\begin{split}
&\tilde{b}_\ast(x,\xi):=\int_{-\infty}^0 \chi_0(\varphi(s;x,\xi))(-s)ds, \quad (x,\xi) \in \R^{2n}, \\
&\tilde{c}_\ast(x,\xi):=\lbrace a_3,\tilde{b}_\ast \rbrace (x,\xi), \quad (x,\xi) \in \R^{2n}.
\end{split}
\]
Observe that $\tilde{b}_\ast$ is supported in a conic neighborhood of $\lbrace \varphi(s;x_0,\xi_0); \ s>0 \rbrace$. Indeed, for each $s>0$
\begin{equation}\label{eq.posb2}
\tilde{b}_\ast(\varphi(s;x_0,\xi_0))=\int_{-\infty}^0 \chi_0(\varphi(s+\sigma;x_0,\xi_0))(-\sigma) d\sigma>0.
\end{equation}
Since $\pi_x(\operatorname{supp}\chi_0)$ is compact and the flow is nontrapping, the integrand vanishes outside a finite interval for each $(x,\xi)$. On compact sets of initial positions and unit frequencies this interval can be chosen uniformly, by the uniform escape estimate for $X_1$.

As for the symbol $\tilde c_*$, integration by parts gives
\[
\begin{split}
\tilde c_*(x,\xi)=\{a_3,\tilde b_*\}(x,\xi)
&=\int_{-\infty}^0(-s)\frac{d}{ds}\chi_0(\varphi(s;x,\xi))\,ds\\
&=\int_{-\infty}^0\chi_0(\varphi(s;x,\xi))\,ds\geq0.
\end{split}
\]
Hence, $\tilde{b}_\ast$ and $\tilde{c}_\ast$ satisfy \eqref{propertyPB} and \eqref{propertyb-c}.

Next, in order to deal with symbols of proper order we then define
\[
b_\ast(x,\xi):=|\widebar{\Xi}^{(2)}(x,\xi)|^{\frac{m}{2}+2}\psi(|\widebar{\Xi}^{(2)}(x,\xi)|)\tilde{b}_\ast(x,\xi), \quad (x,\xi) \in \R^{2n},
\]
where $\psi\in C^\infty(\R)$ is a smooth cutoff function supported in $[1/2,+\infty)$ such that $\psi'\geq 0$ and $\psi(s)=1$ for $s>1$.

Let now $(x_1,\xi_1)=\varphi(s_1;x_0,\xi_0) \in \mathscr{S}_+$, for $s_1>0$, be the point defined in Lemma \ref{lem.x1xi1} and let $(b_+,c_+)$ be the pair of associated symbols constructed in Lemma \ref{lem.E.2}. Let also $\chi\in C_c^\infty(\R^{n})$ be  such that $\chi\geq 0$, $\mathrm{supp}(\chi)\subseteq\{|x|\leq R+1\}$, with $R>0$ such that $\pi_x(\mathrm{supp}(b_+))\subseteq \{|x|\geq R\}$, $\chi\equiv1$ for $|x|\leq R$, and 
\[
\mathrm{supp}(\{a_3,\chi\}b_*)\subseteq \lbrace (x,\xi); \; b_+(x,\xi) \gtrsim \langle x \rangle^k \langle \xi \rangle^m \rbrace.
\] 
Now $\chi \tilde{b}_\ast\in S(-4,k;\rho,0)$ for all $\rho>0$ and $k\geq 0$. This is due to the fact that $\tilde{b}_\ast$ is homogeneous degree $-4$ in the $\xi$ variable. Indeed, by Lemma \ref{lem.hom.prop} and the property $\chi_0(x,\lambda\xi)=\chi_0(x,\xi)$, for $\lambda>0$, we have 
\[
\tilde{b}_\ast(x,\lambda\xi)=\int_{-\infty}^0 (-s) \chi_0(\varphi(\lambda^2s;x,\xi))ds=\lambda^{-4}\int_{-\infty}^0 (-\tau) \chi_0(\varphi(\tau;x,\xi))d\tau=\lambda^{-4}\tilde{b}_\ast(x,\xi).
\]
Then, rewriting $\tilde{b}_*(x,\xi)=\tilde{b}_*(x, |\xi|\frac{\xi}{|\xi|})$, we get $|\tilde{b}_*(x,\xi)|\lesssim |\xi|^{-4}\lesssim\langle\xi\rangle^{-4}\langle x\rangle^k$ for all $k\geq 0$ by the properties of the support.
To see that also $\chi \partial^\alpha_x\partial^\beta_\xi \tilde{b}_*$ satisfies the suitable bounds, it suffices to show that they are also homogeneous of the suitable order with respect to $\xi$. To see this, we proceed as above by using the homogeneity properties of the derivatives of the flow in Lemma \ref{lem.hom.prop}, the fact that $(\partial^\gamma_\xi\chi_0) (x,\lambda \xi)= \lambda^{-|\gamma|} \partial^\gamma_\xi\chi_0 (x, \xi)$, and the boundedness of the domain of integration. This shows that $\chi \tilde{b}_*\in S(-4,k,\rho,0)$. Finally, since $a_3\in S(3,0;0,1)$, by the properties of the calculus we get $\chi \tilde{c}_*=\{a_3,\tilde{b}_*\}\in S(-2,k-1;\rho,0)$.

Then, we define 
\[
b(x,\xi):=\chi(x)e^{Mx_1}b_\ast(x,\xi)+C_0b_+(x,\xi), \quad (x,\xi) \in \R^{2n},
\]
where $M,C_0>0$ are sufficiently large constants to be chosen later, and 
\[
c(x,\xi):=\lbrace a_3,b \rbrace(x,\xi), \quad (x,\xi) \in \R^{2n}.
\]
Note that, since $\chi$ has compact support, then $\chi b_*\in S(m,k;\rho,\delta)$ for all $0\leq \delta\leq 1$ and all $k\geq 0$.  Moreover, due to the support of $\psi$, $b_*(x,\xi)\gtrsim \la \xi\ra^{m}$ when $|\xi|^2\approx |\widebar{\Xi}^{(2)}|>1/2$.
Therefore, due to the inclusion $b_+\in S(m,k;\rho,1-\rho)$ and to the properties of the calculus, we get that $b\in S(m,k;\rho,\delta)$ and $c\in S(m+2,k-1,\rho,\delta)$, for $\delta=1-\rho$.

We now define the sets $\mathscr{E}_+^0(x_0,\xi_0) \subseteq \mathscr{E}_+(x_0,\xi_0)$ by
\[
\mathscr{E}_+(x_0,\xi_0):=\mathrm{supp}(b), \qquad
\mathscr{E}_+^0(x_0,\xi_0):= \lbrace (x,\xi) \in \R^{2n}\setminus 0; \ b(x,\xi) \geq \langle x \rangle^k\langle \xi \rangle^m \rbrace.
\]
Note that, due to the supports of $b_+$ and of $b_*$, we have $(x,\xi)\in\mathscr{E}_+^0$ if $|\xi|\geq c'$ for some $c'>0$ that can be chosen by using multiplicative constants in $b_+,b_*$. 

Let us now discuss the validity of \eqref{eq.asympset} and \eqref{eq.asympset0} for $\mathscr{E}_+(x_0,\xi_0)$ and $\mathscr{E}_+^0(x_0,\xi_0)$. 
They satisfy \eqref{eq.asympset} and \eqref{eq.asympset0} by the construction of $b_+$ (see Lemmas \ref{lem.E.1} and \ref{lem.E.2}) and, by definition, we trivially have $\mathrm{supp}(c) \subseteq \mathrm{supp}(b) \subseteq \mathscr{E}_+(x_0,\xi_0)$.
The fact that $\mathscr{E}_+^0(x_0,\xi_0) \subseteq \mathscr{E}_+(x_0,\xi_0)$ are two conic neighborhoods at high frequencies of $\lbrace \varphi(s;x_0,\xi_0); \ s>0 \rbrace$ follows from \eqref{eq.posb2} together with the fact that, by construction, the symbol $b_+$ defined in Lemma \ref{lem.E.2} is positive along the bicharacteristic curve  $\lbrace \varphi(s;x_1,\xi_1); \ \ s>0 \rbrace$, with $(x_1,\xi_1) \in \mathscr{S}_{+}$ defined above.

Finally, note that \eqref{property-abmk} follows directly by Lemma \ref{lem.E.2}, since such property is required on the scattering set $\mathscr{S}_+$ and since $b$, up to a constant, coincides with $b_+$ on $\lbrace |x |>R+1 \rbrace$. The validity of \eqref{property-abmk} can be proved similarly for $\{a_3,\la x\ra^{-k}\la \xi\ra^{-m}\chi e^{Mx_1}b_*\}$, using the spatial support of $\chi$ and choosing $M$ large enough as in \eqref{eq:normalized-bracket-exp-correction}. Therefore, we obtain \eqref{property-abmk} on the whole $\mathscr{E}_+^0(x_0,\xi_0)$.

The next goal is to prove the positivity of $b$ and $c$. 
The positivity of $b$ follows directly by its definition, while to establish the positivity of $c$, we write
\[
c=\lbrace a_3, b \rbrace=\chi e^{Mx_1}\lbrace a_3, b_\ast \rbrace+b_\ast e^{Mx_1}\lbrace a_3, \chi \rbrace +M \chi e^{Mx_1}b_\ast\partial_{\xi_1}a_3 +C_0 \lbrace a_3,b_+ \rbrace,
\]
and analyze the different terms separately.
We first note that, since 
\[
\ds \left( |\widebar{\Xi}^{(2)}(\varphi(s;x,\xi))|^{\frac{m}{2}+2}\psi(|\widebar{\Xi}^{(2)}(\varphi(s;x,\xi))|) \right)\Bigl|_{s=0}=0,
\]
we get
\begin{equation}\label{eq.posci}
\begin{split}
\lbrace a_3,b_\ast \rbrace&= |\widebar{\Xi}^{(2)}(x,\xi)|^{\frac{m}{2}+2}\psi(|\widebar{\Xi}^{(2)}(x,\xi)|)\lbrace a_3,\tilde{b}_\ast \rbrace + \tilde{b}_\ast \lbrace a_3, |\widebar{\Xi}^{(2)}(x,\xi)|^{\frac{m}{2}+2}\psi(|\widebar{\Xi}^{(2)}(x,\xi)|) \rbrace \\
&=|\widebar{\Xi}^{(2)}(x,\xi)|^{\frac{m}{2}+2}\psi(|\widebar{\Xi}^{(2)}(x,\xi)|)\lbrace a_3,\tilde{b}_\ast \rbrace \geq 0.
\end{split}
\end{equation}
and so 
\begin{equation}\label{eq.cstar}
   \chi e^{Mx_1}\{a_3,b_*\}\geq 0.
\end{equation}
Moreover, since 
\[
\mathrm{supp}(\{a_3,\chi\}b_*)\subseteq \lbrace (x,\xi); \; b_+(x,\xi) \gtrsim \langle x \rangle^k \langle \xi \rangle^m \rbrace, \quad \forall j=1,\dots,n,
\]
by choosing $C_0=C_0(M,R)>0$ sufficiently large and using the properties of $c_+:=\{a_3,b_+\}$ in Lemma \ref{lem.E.2}, 
we have that
\begin{equation}\label{eq.poscii}
b_\ast e^{Mx_1}\lbrace a_3, \chi \rbrace +C_0 \lbrace a_3,b_+ \rbrace\geq 0,
\end{equation}
and finally, 
\[
M\chi e^{Mx_1}b_\ast\partial_{\xi_1}a_3=M\chi e^{Mx_1} b_\ast(2\xi_1^2+a_2)\geq 0.
\]
Hence, by \eqref{eq.posci}, \eqref{eq.poscii} and \eqref{eq.cstar} we obtain 
\[
c= \chi e^{Mx_1}\lbrace a_3, b_\ast \rbrace+b_\ast e^{Mx_1}\lbrace a_3, \chi \rbrace +M\chi e^{Mx_1} b_\ast\partial_{\xi_1}a_3 +C_0 \lbrace a_3,b_+ \rbrace\geq 0,
\]
and, in particular, by Lemma \ref{lem.E.2} and for $C_0$ sufficiently large, we obtain 
\begin{equation}
    c(x,\xi)\gtrsim \la x\ra^{k-1}\la \xi\ra^{m+2} +\frac{M}{2} \la x\ra^{k-2}\la \xi\ra^{m+2}, \quad \forall (x,\xi)\in \mathscr{E}^0_+(x_0,\xi_0).  
\end{equation}
\end{proof}

We finally iterate this argument to construct a nested sequence of subneighborhoods $(\mathscr{E}_{+}^\ell(x_0,\xi_0))_{\ell \in \mathbb{N}_0}$ corresponding to a sequence of symbol pairs $((b_\ell,c_\ell))_{\ell\in \mathbb{N}_0}$ satisfying the aforementioned conditions. This construction is performed through the following theorem. 
\begin{theorem}\label{theo.nested}
Let $(x_0,\xi_0) \in \R^{2n} \setminus 0$. For all $J \in \mathbb{N}$ and for all $m_\ell, \in \mathbb{R}$, $k_\ell\geq 0$ with $\ell=1,\dots,J$, there exist
symbol pairs $(b_\ell,c_\ell)$, with $b_\ell \in S(m_\ell,k_\ell;\rho,\delta), c_\ell \in S(m_\ell+2,k_\ell-1;\rho,\delta)$, and a nested sequence of subneighborhoods 
\[
\mathscr{E}_{+}^J(x_0,\xi_0) \subseteq \mathscr{E}_+^{J-1}(x_0,\xi_0)\subseteq \dots \subseteq \mathscr{E}_+^{0}(x_0,\xi_0) 
\]
of the forward bicharacteristic $\lbrace \varphi(s;x_0,\xi_0); \; s>0 \rbrace$, defined by 
\[
\mathscr{E}_{+}^\ell(x_0,\xi_0)=\lbrace (x,\xi); \ b_\ell(x,\xi) \geq \langle x \rangle^{k_\ell} \langle \xi \rangle^{m_\ell}\rbrace , \quad \ell=1,\dots,J,
\]
such that, for each $\ell$, 
\[
\mathrm{supp}(c_\ell)\subseteq \mathrm{supp}(b_\ell)\subseteq \mathscr{E}_{+}^{\ell-1}(x_0,\xi_0),
\]
and they also satisfy the same conclusions of Theorem \ref{thm.constrE-} on the corresponding $\mathscr{E}_{+}^{\ell}(x_0,\xi_0)$.
Once again, if $k_\ell=0$, then the property \eqref{eq.propk>0bis} holds with $k'>1$.
\end{theorem}
\begin{proof}
We prove the theorem for $J=1$. The proof for a general $J \in \mathbb{N}$ follows by iterating the same argument. If $J=1$ and $m_1\in \mathbb{R}, k_1>0$ we have to prove that we can construct $b_1 \in S(m_1,k_1;\rho,\delta), \; c_1 \in S(m_1+2,k_1-1;\rho,\delta)$ and $\mathscr{E}_+^{1}(x_0,\xi_0)$ satisfying the conditions in the statement.

To this end, let us consider $(b,c)$ the symbol pairs constructed in Theorem \ref{thm.constrE-}. Without loss of generality, we may assume that $b \in S(0,0;\rho,\delta)$, for otherwise replace it by the expression $\langle x \rangle^{-k} \langle \xi \rangle^{-m} b$, which by the construction in Theorem \ref{thm.constrE-} (i.e. \eqref{property-abmk} on $\mathscr{E}_+^0(x_0,\xi_0)$) satisfies \eqref{propertyb-c} (cf. Theorem 3.5 in \cite{CKS}).  In particular, 
\[
\mathscr{E}^0_{+}(x_0,\xi_0)=\big\lbrace (x,\xi); b(x,\xi) \geq 1 \rbrace.
\]
Let us consider now 
$\psi \in C^\infty([0,+\infty))$ be such that $\psi(s)=0$ for $0 \leq s \leq 1$, $\psi(s)>0$ if $s>1$ and $\psi'(s)\geq 0$ for all $s \geq 0$, and set 
\[
\tilde{b}_1(x,\xi):=\psi(b(x,\xi)), \quad \tilde{c}_1(x,\xi)=\lbrace a_3,\tilde{b}_1 \rbrace(x,\xi).  
\]
From the properties of $\psi$ it follows that $\tilde{b}_1 \geq 0$, $\mathrm{supp}(\tilde{b}_1) \subseteq \mathscr{E}^0_{+}(x_0,\xi_0)$ and 
\[
\tilde{c}_1(x,\xi)=\psi'(b(x,\xi))\lbrace a_3,b \rbrace(x,\xi) \geq 0.  
\]
Furthermore, by reasoning as before, due to Faà di Bruno's formula (see Theorem \ref{theo.fdb}) we have $\tilde{b}_1 \in S(0,0;\rho,\delta)$ and $\tilde{c}_1 \in S(2,-1;\rho,\delta)$.
Finally, to get $b_1$ and $c_1$ of the required order and satisfying the required property, it is possible to argue as in Lemma \ref{lem.E.2}.

\end{proof}

\section{The microlocal smoothing effect}\label{sec.smoothing}
In this section we prove the smoothing results stated in the Introduction for Cauchy problems of the form 
\begin{equation}\label{eq.mainprob}
    \begin{cases}
        D_t u -Au=0 \\
        u(0,\cdot)=u_0 \in L^2(\R^n),
    \end{cases}
\end{equation}
where $A=D_{x_1}A_2$ satisfies the properties \ref{hp1}-\ref{hp3}. 

As a consequence of the smoothing effect for \eqref{eq.mainprob}, we get the same property for
\begin{equation}\label{eq.mainpro2b}
    \begin{cases}
        D_t v +Av=0 \\
        v(0,\cdot)=v_0 \in L^2(\R^n),
    \end{cases}
\end{equation}
since the solution to this problem is $v(t,x)=u(t,-x)$, where $u$ is the solution of \eqref{eq.mainprob} with 
\[ A=D_1\sum_{i,j}D_i a_{ij}(-x)D_j, \quad u_0(-x)=v_0(x).\]
Under this reflection, the forward regions for the reflected operator become backward regions for the original operator. Thus the microlocal statements for \eqref{eq.mainpro2b} use backward moment assumptions; the global weighted statements are unchanged.

\begin{remark}
Note that, when $A_2=-\Delta$ is the positive Laplacian in \eqref{eq.mainprob}, we obtain the (possibly multidimensional) linear ZK-type Cauchy problem 
\begin{equation}
    \begin{cases}
        \partial_t u +\partial_{x_1}\Delta u=0 \\
        u(0,\cdot)=u_0 \in L^2(\R^n).
    \end{cases}
\end{equation}
\end{remark}

Let us now return to problem \eqref{eq.mainprob}. The existence of a solution operator $\lbrace U(t) \rbrace_{t \in [0,T]}$, with  
\[
U(t):L^2(\R^n)\rightarrow L^2(\R^n), \quad u_0 \mapsto u(t,\cdot), \quad \quad t\in [0,T], 
\]
that sends the initial datum $u_0$ to the (unique) solution of \eqref{eq.mainprob} is established in \cite{Ch} (cf. \cite{FT}), since our assumptions \ref{hp1}-\ref{hp3} (see also Lemma \ref{lemma.nontrapping}) imply that the symbol $a$ of $A$ satisfies the conditions assumed in that work. Moreover, by Theorem 1.2 in \cite{Ch} we also have 
\begin{equation}\label{eq.boundL2}
\sup_{0 \leq t \leq T}\|u(t,\cdot)\|_{L^2} \leq C_T \|u_0\|_{L^2},
\end{equation}
for some $C_T>0$.

Now, let $b \in S(m,k;\rho,\delta)$ be the symbol defined in subsection \ref{subsec.bc} and let $B=\mathrm{Op}(b)$ be the corresponding pseudodifferential operator. If $u=u(t,x) \in C([0,T];L^2(\R^n))$ is a (sufficiently regular) solution of the problem \eqref{eq.mainprob}, we have 
\begin{align}\label{eq.partialtquadform0}
\partial_t  (u,Bu)_{L^2}&=(iAu,Bu)_{L^2}+(u,iBAu)+(u,(\p_tB)u)_{L^2}\\
&=-(u,i[A,B]u)_{L^2}+(u, i(A-A^*)Bu)+(u,(\p_tB)u)_{L^2},
\end{align}
which gives
\begin{equation}\label{eq.partialtquadform}
\partial_t \re (u,Bu)_{L^2}= -\re(u,i[A,B]u)_{L^2}+\re(u, i(A-A^*)Bu)+\re (u, (\partial_t B)u)_{L^2}.
\end{equation}
We now write 
\[
A=D_{x_1}A_2=D_{x_1}\sum_{i,j=1}^nD_{x_i}a_{ij}(x)D_{x_j}=A_3+A',
\]
where $A_3=\Op(a_3)\in \Psi^3(\R^n)$ is the quantization of the principal symbol $a_3$ defined in \eqref{eq.defa3}, namely
\[
A_3=\sum_{i,j=1}^n a_{ij}(x)D_{x_1}D_{x_i}D_{x_j},
\]
and the lower term $A'=\Op(a') \in \Psi^2(\R^n)$, $a'\in S(2,-2;\rho,\delta)$, is defined by 
\[
A'=\sum_{i,j=1}^n(D_{x_1}a_{ij})(x)D_{x_i}D_{x_j}+D_{x_1}\sum_{i,j=1}^n(D_{x_i}a_{ij})(x)D_{x_j}.
\]
This gives that
\begin{align}
    i[A,B]=\Op(\{a_3,b\}-i\{\{a_3,b\}\}+\{a',b\})+E_0+R_{L+1},
\end{align}
where
\begin{align}
    \{\{a_3,b\}\}&=\sum_{|\alpha|=2}\frac{1}{\alpha!}(\partial^{\alpha}_\xi a_3\partial^{\alpha}_xb-\partial^{\alpha}_x a_3\partial^{\alpha}_\xi b)
\end{align}
belongs to the class $S(m+1, k-2\min\{\delta,1-\rho \};\rho,\delta)$,  $R_{L+1}$, for $L$ sufficiently large, is a bounded operator on $L^2(\mathbb{R}^n)$, and $E_0=\mathrm{Op}(e_0)$, with $e_0\in S(m,k-2\min\{\delta,1-\rho\};\rho,\delta)$ given by
\begin{align}
    e_0&= \sum_{|\alpha|=2}\frac{1}{\alpha!}(-i)^{|\alpha|-1}(\partial^{\alpha}_\xi a'\partial^{\alpha}_xb-\partial^{\alpha}_x a'\partial^{\alpha}_\xi b)+\sum_{|\alpha|=3}^L \frac{1}{\alpha!} (-i)^{|\alpha|-1}(\partial^{\alpha}_\xi a\partial^{\alpha}_xb-\partial^{\alpha}_x a\partial^{\alpha}_\xi b).\label{eq.E0}
\end{align}
To check that the above symbols belong to the indicated symbol classes, it suffices to observe that, due to the space-decay properties of the coefficients $a_{ij}$, i.e. 
\[
|\p_x^\alpha a_{ij}(x)|\lesssim \la x\ra ^{-2-|\alpha|} \quad \text{ for} \ \ |\alpha|\geq 1,
\]
we have
\begin{align}
    \partial^{\alpha}_\xi a\partial^{\alpha}_xb\in S(m+3-|\alpha|, k-\delta|\alpha|;\rho,\delta), \quad \partial^{\alpha}_x a\partial^{\alpha}_\xi b\in S(m+3-|\alpha|,k-2-(1-\rho)|\alpha|;\rho,\delta),
\end{align}
therefore $\partial^{\alpha}_\xi a\partial^{\alpha}_xb-\partial^{\alpha}_x a\partial^{\alpha}_\xi b \in S(m+3-|\alpha|,k-\min\{\delta,1-\rho\}|\alpha|;\rho,\delta)$. The same argument applies when $a$ is replaced by $a_3$ or $a'$, and taking into account that $a'$ is a symbol of order two we conclude the class for $e_0$.

Now we observe that $i\{\{a_3,b\}\}$ and the principal part of $\{a',b\}$ are purely imaginary symbols. Consequently,
\[
\re(u,\Op(i\{\{a_3,b\}\}-\{a',b\})u)_{L^2}
=\re(u,E_1u)_{L^2}
\]
for an operator $E_1$ with symbol in $S(m,k-2\min\{\delta,1-\rho\};\rho,\delta)$.

On the other hand,
\begin{equation}
A-A^*=[D_1,A_2]
=-i\sum_{i,j}D_i(\partial_{x_1}a_{ij})D_j=\Op(\sum_{i,j}(-i(\partial_{x_1}a_{ij})\xi_i\xi_j-(\p^2_{i1}a_{ij})\xi_j),
\end{equation}
and, by the symbolic calculus, we have
\begin{align}
    i(A-A^*)B&=\Op(b\sum_{i,j=1}^n(\partial_{x_1}a_{ij})\xi_i\xi_j)\\
    &+\Op(-i\sum_{i,j=1}^n(\p^2_{i1}a_{ij})\,\xi_j\,b
-2i\sum_{i,j=1}^n(\p_1a_{ij})\,\xi_j\,\partial_{x_i}b)+R_m\\
&=:\Op(r_{m+2})+i\Op(r_{m+1})+ \Op(r_m) \label{eq.Rm.symbols}
\end{align}
where $r_m\in S(m, k-3-2\delta;\rho,\delta)$, $r_{m+1}\in S(m+1,k-3-\delta;\rho,\delta)$ and  $r_{m+2}\in S(m+2,k-3;\rho,\delta)$.

Therefore, from \eqref{eq.partialtquadform}, $\partial_tB=0$, and the fact that $ir_{m+1}$ is purely imaginary, we obtain (absorbing its real quadratic form into $E_1$)
\begin{align}
\partial_t\re(u,Bu)_{L^2}
&=-\re(u,\Op(\{a_3,b\})u)_{L^2}
+\re(u,\Op(r_{m+2})u)_{L^2}\nonumber\\
&\quad+\re(u,(-E_0+E_1+\Op(r_m))u)_{L^2}
-\re(u,R_{L+1}u)_{L^2}\nonumber\\
&=-\tfrac12\re(u,\Op(\{a_3,b\})u)_{L^2}
-\tfrac12\re(u,\Op(\{a_3,b\}-2r_{m+2})u)_{L^2}\nonumber\\
&\quad+\re(u,(-E_0+E_1+\Op(r_m))u)_{L^2}
-\re(u,R_{L+1}u)_{L^2}.\label{eq.commutator}
\end{align}

Hence, by \eqref{eq.commutator} we obtain
\begin{align}
    \partial_t \re(u, Bu)_{L^2}\leq-\frac 1 2\re(u,C u)_{L^2}-\frac 1 2\re(u,Du)_{L^2}+\re(u, E u)_{L^2}+C\|u\|_{L^2_x}^2, \label{eq.EnergyEstimate}
\end{align}
for $C=\Op(\{a_3,b\})$, $D=\Op(d)=\Op(\{a_3,b\}-2r_{m+2})$, and
$E=\Op(e)=-E_0+E_1+\Op(r_m)$, with $e\in S(m,k-2\min\{\delta,1-\rho\};\rho,\delta)$.

The next goal is to show some properties of the operator $D$ that we will need in the proof of our results. 

Recall that by \eqref{eq.Rm.symbols} we have
\[d=\lbrace a_3,b \rbrace-2r_{m+2}= \lbrace a_3,b \rbrace-2b\sum_{i,j=1}^n(\p_{x_1}a_{ij})\xi_i\xi_j,\]
so  $d\in S(m+2,k-1;\rho,\delta)$ and $\mathrm{supp}(d)\subseteq\mathrm{supp}(b) =\mathscr{E}_+(x_0,\xi_0).$
Moreover, $d\geq 0$ on its support and satisfies suitable bounds on $\mathscr{E}_+^0(x_0,\xi_0)$. In fact,
for $k>0$, by \eqref{eq.propk>0}, the decay properties of $a_{ij}$, and by choosing $M>0$ in \eqref{eq.propk>0} large enough, for all $(x,\xi)\in \mathscr{E}_+^0(x_0,\xi_0)$, we obtain 
\begin{align}
    d\gtrsim \la x\ra^{k-1}\la \xi\ra^{m+2}+\frac{M}{2}\la x\ra^{k-2}\la \xi\ra^{m+2}-C\la x\ra^{k-3}\la \xi\ra^{m+2}\gtrsim \la x\ra^{k-1}\la \xi\ra^{m+2}.
\end{align}
To see that $d\geq 0$ on its support  we denote by
\[q=\sum_{i,j}(\partial_1a_{ij})\xi_i\xi_j,
\]
and note that 
\[
\begin{aligned}
d=c-2qb\geq c-2|q|b
={}&e^{Mx_1}\chi
\big[c_*+(M\partial_{\xi_1}a_3-2|q|)b_*\big]\\
&+e^{Mx_1}(H_{a_3}\chi)b_*
+C_0(c_+-2|q|b_+).
\end{aligned}
\]
The term in the first line is nonnegative for $M$ sufficiently large, since, by \eqref{eq:normalized-bracket-exp-correction},
\[\partial_{\xi_1}a_3\ge c_0|\xi|^2,\qquad
|q|\le C_q\langle x\rangle^{-3}|\xi|^2,\quad c_*\geq 0.\]
The term 
\[c_+-2|q|b_+\geq (C_1M \la x\ra^{-2}|\xi|^{2}- C_2\la x\ra^{-3}|\xi|^{2})b_+\geq 
\frac{C_1M}{2}\langle x\rangle^{-2}|\xi|^2b_+, \quad \] 
by \eqref{eq:normalized-bracket-exp-correction} and choosing $M$ sufficiently large, since $C_1,C_2>0$. Finally, 
since $\chi$ is such that $b_*H_{a_3}\chi$ is supported where 
$c_+(x,\xi)\gtrsim \la x\ra^{k-1}\la \xi\ra^{m+2}+M  \la x\ra^{k-2}\la \xi\ra^{m+2}$, for $C_0=C_0(M)$ large enough, we get 
\[e^{Mx_1}(H_{a_3}\chi)b_*
+C_0(c_+-2|q|b_+) \geq 0,\]
which yields $d\geq0$. For $k=0$, the same argument uses the tail bound
$c_+\gtrsim M\sigma\langle x\rangle^{-1-\sigma}|\xi|^2b_+$,
where $0<\sigma=\min\{2,k'-1\}\leq2$; this again dominates $2|q|b_+$ for $M$ large enough.

Finally, observe that for $b_\ell,c_\ell$ as in Theorem \ref{theo.nested} we have a corresponding symbol
\begin{equation}\label{eq.dell}
    d_\ell:=\lbrace a_3,b_\ell \rbrace-2b_\ell\sum_{i,j=1}^n(\p_{1}a_{ij})\xi_i\xi_j=c_\ell-2b_\ell\sum_{i,j=1}^n(\p_{1}a_{ij})\xi_i\xi_j
\end{equation}
which inherits the same properties of $d$ but in terms of the suitable orders $m_\ell,k_\ell$ in Theorem \ref{theo.nested} with the suitable choice of the cutoff $\chi$ in $d_\ell$. Hence, $d_\ell\in S(m_\ell+2,k_\ell-1;\rho,\delta)$, it is nonnegative, has support satisfying
\begin{equation}
    \mathrm{supp}(d_\ell)\subseteq \mathrm{supp}(b_\ell)\subseteq \mathscr{E}_+^{\ell-1}(x_0,\xi_0),
\end{equation}
where the sets $\mathscr{E}_+^{\ell}(x_0,\xi_0)$ are those defined in Theorem \ref{theo.nested}, and satisfies positive lower on $\mathscr{E}_+^{\ell}(x_0,\xi_0)$.

The properties of $b_\ell, c_\ell, d_\ell$ and the energy estimate \eqref{eq.EnergyEstimate} are the core of the proofs of the forthcoming results. Throughout the following arguments, the nonnegative symbols to which the G{\aa}rding inequality (see Theorem \ref{thm.GardingInequality}) is applied have smooth square roots in the corresponding symbol
classes. This is due to the fact that the symbols in \ref{sec.smoothing} are constructed through cutoff functions that vanish at infinite order at any zeros (by a proper choice of the cut-off functions involved).

We consider $(x_0,\xi_0) \in \R^{2n}\setminus 0$ and $0 \leq \rho < \delta \leq 1$, $K \geq 0$. The set $\mathscr{E}_+(x_0,\xi_0)$ is a suitable neighborhood  of the forward bicharacteristic $\lbrace \varphi (s;x_0,\xi_0); s\geq 0 \rbrace$, and $b_0,c_0$ are two nonnegative symbols supported in $\mathscr{E}_+(x_0,\xi_0)$ such that 
\begin{equation}\label{b0c0}
0 \leq b_0 \in S(0,K;\rho,\delta), \quad 0 \leq c_0=\lbrace a_3,b_0 \rbrace \in S(2,K-1;\rho,\delta).
\end{equation}
Moreover, on denoting 
\begin{equation}\label{eq.E0.thm}
\mathscr{E}^0_{+}(x_0,\xi_0):=\big\lbrace (x,\xi); b_0(x,\xi) \geq \langle x \rangle^K \rbrace \subseteq \mathscr{E}_+(x_0,\xi_0), 
\end{equation}
we also have 
\[
0 \leq \lbrace a_3, \langle x \rangle^{-K}b_0 \rbrace (x,\xi), \quad \forall (x,\xi) \in \mathscr{E}^0_+(x_0,\xi_0)\cap \mathscr{S}^{+}.
\]
Finally, we define $s_0:=b_0^{1/2} \in S(0,K/2;\rho,\delta)$. Note that  due to its definition, $b_0$ has a smooth square root even if it is nonnegative.  

We are now ready to state the first smoothing result for solutions to \eqref{b0c0}.

\begin{theorem}\label{thm.nestedOperators}
    Let $u_0 \in L^2(\R^n)$ and, for $S_0=\Op(s_0)$ where $s_0:=\sqrt{b_0}$, assume 
    \begin{equation}\label{moment.condition}
    (u_0,S_0^\ast S_0u_0)_{L^2}<+\infty. 
     \end{equation}
    Let also $0 <r\leq (\delta-\rho)$ be such that $K/2r=J$ is an integer, where $K$ is such that \eqref{b0c0} holds. 
    Then, for all $\ell=1,\dots,J$ there exist
    symbol pairs $(b_\ell,c_\ell)$, with $b_\ell \in S(0,2\ell r;\rho,\delta), c_\ell \in S(2,2\ell r-1;\rho,\delta)$, and a nested sequence of subneighborhoods 
    \[
    \mathscr{E}_{+}^J(x_0,\xi_0) \subseteq \mathscr{E}_+^{J-1}(x_0,\xi_0)\subseteq \dots \subseteq \mathscr{E}_+^{0}(x_0,\xi_0)
    \]
    containing the forward bicharacteristic $\lbrace \varphi(s;x_0,\xi_0); s\geq 0 \rbrace$ and defined by 
    \[
    \mathscr{E}_{+}^\ell(x_0,\xi_0)=\lbrace (x,\xi); \ b_\ell(x,\xi) \geq \langle x \rangle^{2\ell r} \rbrace , \quad \ell=1,\dots,J,
    \]
    such that, for each $\ell$, 
    \[
    \mathrm{supp}(c_\ell)\subseteq \mathrm{supp}(b_\ell)\subseteq \mathscr{E}_{+}^{\ell-1}(x_0,\xi_0),
    \] 
    and, for all $T>0$, on denoting $B_\ell=\Op(b_\ell)$, $C_\ell=\Op(c_\ell)$ and $u(t,\cdot)$ the solution of \eqref{eq.mainprob}, we have 
    \begin{align}
        (i) \ \sup_{ 0 \leq t \leq T}\Bigl(\re(u(t,\cdot),B_\ell u(t,\cdot))_{L^2}\Bigr)<+\infty \\
        (ii) \ \re\int_0^T(u(t,\cdot),C_\ell u(t,\cdot))_{L^2}dt<+\infty.
    \end{align}
\end{theorem}
To prove Theorem \ref{thm.nestedOperators} we will need a technical lemma contained in \cite{CKS} that we recall here to keep the exposition self contained.

\begin{lemma}[Lemma 6.2 in \cite{CKS}]\label{lemmaCKS}
   Suppose there is a sequence of symbols $b_\ell\in S(0,2\ell r;\rho, \delta)$ for $1\leq \ell\leq L$ and with $r>0$, such that
   \begin{itemize}
       \item [(i)] $b_\ell(x,\xi)\geq 0$;
       \item [(ii)] $\mathrm{supp}(b_\ell)\subseteq \{ (x,\xi): b_{\ell-1}(x,\xi)\geq \la x\ra^{2(\ell-1)r}\}$.
   \end{itemize}
   Consider any symbol $p\in S(0,2Lr;\rho,\delta)$ with $\mathrm{supp}(p)\subseteq \{ (x,\xi): b_{L}(x,\xi)\geq \la x\ra^{2Lr}\}$. Then, there exist some constant and some constants  $M_0,\dots,M_L>0$ we have
   \begin{equation}
       |(\psi, \Op(p)\psi)_{L^2}|\leq \sum_{\ell=1}^L M_\ell \re(\psi,\Op(b_\ell)\psi)_{L^2}+M_0\|\psi\|_{L^2}^2, \quad \psi \in \mathscr{S}(\R^n).
   \end{equation}
\end{lemma}

\begin{remark}
    Note that in Theorem \ref{thm.nestedOperators} we require that $K=2Jr$ with $0<r\leq \delta-\rho$. These choices of $K$ and $r$ are made in order to be able to apply Lemma \ref{lemmaCKS} to the sequences of symbols we have in the proof, specifically the sequences $b_\ell$ and $e_\ell$ in the proof of Theorem \ref{thm.nestedOperators}. In particular, the choice $K=2Jr$ is made  to be in the position to apply the above lemma to sequences of symbols with supports contained in 
    \[\{(x,\xi): b_0(x,\xi)\geq \la x\ra^{2Jr} \},\]
    which indeed happens when $K=2Jr$ in \eqref{b0c0}. 
     The condition $0<r\leq \delta-\rho$  ensures that the sequences of symbols $b_\ell,e_\ell$ in the proof of Theorem \ref{thm.nestedOperators} have the $j$-th element in the class $S(0,2jr;\rho,\delta)$ as required in the hypotheses of Lemma \ref{lemmaCKS}.
\end{remark}

\begin{remark}\label{rmk.approx.argument}

The proofs of Theorems \ref{thm.nestedOperators} and \ref{thm.mainest} below are first carried out for solutions with Schwartz initial data. In this case, Sobolev well-posedness (see \cite{C}) and commutator estimates for polynomial weights ensure the regularity and decay required to justify all calculations. The conclusions for general initial data satisfying the distributional moment assumption follow by the approximation argument described below.

The key point is to obtain estimates that are uniform in the approximation parameter and depend on the initial data only through their \(L^2\) norm and initial moment. 
If $u_0 \in L^2(\R^n)$ we consider 
\[
u_{0,j}=\varphi(x/j)\Op(\psi(\xi/j))u_0,
\]
where $\varphi,\psi \in C_c^\infty(\R^n)$ equal $1$ near the origin. Then 
$u_{0,j} \in \mathscr S(\mathbb R^n)\) and \(u_{0,j}\to u_0\) in \(L^2\), as $j \rightarrow +\infty$.

We slightly shrink the microlocal neighborhoods and construct a new initial symbol \(b_0'\), with the same properties as \(b_0\), supported in the elliptic region of the original symbol \(b_0\) appearing in the hypothesis. A microlocal parametrix for \(\operatorname{Op}(\sqrt{b_0})\) then yields
\[
\bigl\|\operatorname{Op}(\sqrt{b_0'})u_{0,j}\bigr\|_{L^2}^{2}
\le C\left(
\bigl\|\operatorname{Op}(\sqrt{b_0})u_0\bigr\|_{L^2}^{2}
+\|u_0\|_{L^2}^{2}
\right),
\]
with \(C\) independent of \(j\). The resulting estimates are therefore uniform in \(j\) and pass to the limit by continuous dependence in \(L^2\) and weak lower semicontinuity of the nonnegative quadratic forms.

In the following proofs, we use \(b_0'\) as the starting symbol, denoting it again by \(b_0\).
\end{remark}

\proof[Proof of Theorem \ref{thm.nestedOperators}]
We assume that $u_0\in\mathscr{S}(\R^n)$, so that the corresponding solution $u$ has the regularity and decay properties needed to justify the computations below. The conclusion under the  assumptions of the theorem then follows by the approximation argument described in Remark \ref{rmk.approx.argument}.

We start by using Theorem \ref{theo.nested} to construct symbol pairs $(b_\ell,c_\ell)$ satisfying the same properties of $(b_0,c_0)$ in \eqref{b0c0}, that is,
\[0\leq b_\ell \in S(m_\ell,k_\ell;\rho,\delta), \quad   c_\ell=\{a_3,b_\ell\} \in S(m_\ell+2,k_\ell-1;\rho,\delta),\]
\[\mathrm{supp}(c_\ell)\subseteq \mathrm{supp}(b_\ell)\subseteq \mathscr{E}_+^{\ell-1}(x_0,\xi_0)\]
and 
\[
\mathscr{E}_+^\ell(x_0,\xi_0):=\lbrace (x,\xi); \; b_1(x,\xi) \geq \langle x \rangle^{k_\ell}\} \neq \emptyset \subseteq \lbrace (x,\xi); \; b_0(x,\xi) \geq \langle x \rangle^K \rbrace.
\]
Associated with the sequences $b_\ell,c_\ell$, we have the sequence of symbols $d_\ell$ as in \eqref{eq.dell}.

Then, we split the analysis into the cases $J=1$ and $J>1$.

When $J=1$ then $K=2r$ and we take $b_1$ from the above construction with $m_1=0$ and $k_1=K=2r$, so we have  
\[
\mathrm{supp}(b_1)\subseteq \lbrace (x,\xi); \; b_0(x,\xi) \geq \langle x \rangle^K \rbrace,
\]
and
\[
\mathscr{E}_+^1(x_0,\xi_0):=\lbrace (x,\xi); \; b_1(x,\xi) \geq \langle x \rangle^K\} \neq \emptyset.
\]
If $t \in \R$ and $u(t,\cdot)$ is the solution of \eqref{eq.mainprob}, then \eqref{eq.EnergyEstimate} holds with $B=B_1=\Op(b_1)$ and $C=C_1=\Op(c_1)$, so  after integration in time in  \eqref{eq.EnergyEstimate} we get  
\begin{equation}
\begin{split}
\re &(u(t,\cdot),B_1u(t,\cdot))_{L^2}+\frac 1 2\int_0^t \re (u(\tau,\cdot), C_1u(\tau,\cdot))d\tau +\frac{1}{2}\int_0^t \re (u(\tau,\cdot), D_1u(\tau,\cdot))d\tau \\
& \lesssim \re (u_0,B_1u_0)_{L^2}+\int_0^t \re (u(\tau,\cdot), Eu(\tau,\cdot))d\tau+ \|u_0\|_{L^2}^2,
\end{split}
\end{equation}
where in the last line we have used \eqref{eq.boundL2}.

Now, by the properties of the calculus, we have 
\[\re (u_0,B_1u_0)_{L^2}\lesssim (u_0,s_0^*s_0 u_0)_{L^2}+C\|u_0\|_{L^2}<+\infty,\]
and by \eqref{eq.E0} with $m=0$, we also have
\begin{equation}
\int_0^t\re (u(\tau,\cdot), Eu(\tau,\cdot))\lesssim \|u\|_{L^2}^2\lesssim \|u_0\|_{L^2}^2<+\infty,
\end{equation}
Therefore, putting all the bounds together, we obtain
\begin{equation}
    \label{eq.EnEst2}
    \re (u(t,\cdot),B_1u(t,\cdot))_{L^2}+\frac1 2\int_0^t \re (u(\tau,\cdot), C_1u(\tau,\cdot))d\tau+\frac 1 2\int_0^t \re (u(\tau,\cdot), D_1u(\tau,\cdot))d\tau<+\infty.
    \end{equation}
Note that $b_1, c_1,d_1$ are nonnegative and have a smooth square root due to their definition. This allows us us to use the Gårding inequality in Theorem \ref{thm.GardingInequality} and get

\begin{align}
\re(u(t,\cdot),B_1u(t,\cdot))_{L^2}&\geq \|S_0^{b_1}u(t,\cdot))\|^2_{L^2}-C\|u(t,\cdot)\|_{L^2}^2\gtrsim -C_T\|u_0\|^2_{L^2},\quad S_0^{b_1}:=\Op(\sqrt{b_1}) \label{eqb1}\\
\re(u(t,\cdot),D_1u(t,\cdot))_{L^2}&\geq \|S_0^{d_1}u(t,\cdot))\|^2_{L^2}-C\|u(t,\cdot)\|_{L^2}^2\gtrsim -C_T\|u_0\|^2_{L^2},\quad S_0^{d_1}:=\Op(\sqrt{d_1}) \label{eqd1}\\
\re(u(t,\cdot),C_1u(t,\cdot))_{L^2}&\geq \|S_0^{c_1}u(t,\cdot))\|^2_{L^2}-C\|u(t,\cdot)\|_{L^2}^2\gtrsim -C_T\|u_0\|^2_{L^2},\quad S_0^{c_1}:=\Op(\sqrt{c_1}).\label{eqc1}
\end{align}
Then, using these inequalities, at any fixed time we have
\begin{align}
 -C\|u_0\|_{L^2}^2&\leq   \re(u(t,\cdot),B_1u(t,\cdot))_{L^2}\lesssim\|u_0\|_{L^2}^2 +\|S_0 u_0\|^2_{L^2}\\
 -C\|u_0\|_{L^2}^2& \leq  \re \int_0^t (u(\tau,\cdot),C_1u(\tau,\cdot))_{L^2}\lesssim \|u_0\|_{L^2}^2+\|S_0 u_0\|^2_{L^2},
\end{align}
where the constants $C, C'$ may depend on time are uniformly bounded on $[0,T]$. 
Then, after the application of the supremum in time to the previous terms  we conclude the desired estimates in the case
 $\ell=1,J=1$, that is
\begin{align}
    \sup_{0\leq t\leq T}\Big(\re(u(t,\cdot),B_1u(t,\cdot))_{L^2}\Big)&<\infty,\\
    \re\int_0^T(u(t,\cdot),C_1u(t,\cdot))_{L^2}dt&<\infty.
\end{align}
 Note that $b_1\in S(0, 2r; \rho,\delta)$ and, consequently, $c_1=\{a_3,b_1\}\in S(2, 2r-1;\rho,\delta)$ and $d_1\in S(2,2r-1;\rho,\delta)$.

We are now left with proving the statement for every $\ell\leq J$ and  $J>1$. To do that, we proceed by induction on $\ell$ knowing, by the previous step, that the bounds are true for $\ell=1$. Then we assume that the bounds are true for all $0\leq j<\ell$, with $\ell\leq J$. Recall that the existence of the pairs $(b_j,c_j)$, $j\leq L$, is granted by Theorem \ref{theo.nested}. In particular, the sequence of symbols is chosen so that $b_j\in S(0, 2jr;\rho,\delta)$, $c_j\in S(0, 2jr-1;\rho,\delta)$ and there exist neighborhoods 
\[
    \mathscr{E}_{+}^\ell(x_0,\xi_0) \subseteq \mathscr{E}_+^{\ell-1}(x_0,\xi_0)\subseteq \dots \subseteq \mathscr{E}_+^{0}(x_0,\xi_0)
    \]
    containing the forward bicharacteristic $\lbrace \varphi(s;x_0,\xi_0); s\geq 0 \rbrace$. Moreover, 
    \[
    \mathscr{E}_{+}^j(x_0,\xi_0)=\lbrace (x,\xi); \ b_j(x,\xi) \geq \langle x \rangle^{2j r} \rbrace , j=0,\ldots,\ell,
    \]
    satisfy, for each $j\leq L$, the property
    \[
    \mathrm{supp}(c_j)\subseteq \mathrm{supp}(b_j)\subseteq \mathscr{E}_{+}^{j-1}(x_0,\xi_0).
    \] 
Since by construction $ \mathrm{supp}(b_j)\subseteq \mathscr{E}_{+}^{j-1}(x_0,\xi_0)\subseteq \lbrace (x,\xi); \ b_0(x,\xi) \geq \langle x \rangle^{2Jr} \rbrace$ for all $j<\ell$,
then Theorems \ref{thm.UpperBound} and \ref{thm.L2boundedness} yield
\begin{align}
\re( u_0, B_\ell u_0 )_{L^2}\lesssim ( u_0, s_0^*s_0 u_0)_{L^2}+C\| u_0\|_{L^2}^2.\label{eq.bddBell}
\end{align}
Now, we apply 
\eqref{eq.EnergyEstimate} with $B=B_\ell=\Op(b_\ell)$ and $C=C_\ell=\Op(c_\ell)$, which,  after integration in time in \eqref{eq.EnergyEstimate}, gives 
\begin{equation}\label{eq.EEbl2}
\begin{split}
\re &(u(t,\cdot),B_\ell u(t,\cdot))_{L^2}+\frac 1 2 \int_0^t \re (u(\tau,\cdot), C_\ell u(\tau,\cdot))_{L^2}d\tau+\frac 1 2 \int_0^t \re (u(\tau,\cdot), D_\ell u(\tau,\cdot))_{L^2}d\tau \\
& \lesssim \re (u_0,B_\ell u_0)_{L^2}+\int_0^t \re (u(\tau,\cdot), E_\ell u(\tau,\cdot))_{L^2}d\tau+ \|u_0\|_{L^2}^2
\end{split}
\end{equation}
where $E_\ell$ is obtained as in \eqref{eq.EnergyEstimate}.
The first term on the RHS of the inequality is bounded due to \eqref{eq.bddBell},
and the three terms on the LHS are bounded from below by the  Gårding inequality (our symbols have smooth square roots as in the previous step) and by Lemma \ref{lemmaCKS}, that is, 
for
\[S_0^{b_\ell }:=\Op(\sqrt{b_\ell }),\quad S_0^{c_\ell }:=\Op(\sqrt{c_\ell }), \quad S_0^{d_\ell }:=\Op(\sqrt{d_\ell }),\quad Q_0(u_0)^2:=\|u_0\|^2_{L^2}+\|S_0 u_0\|^2_{L^2}\]
\begin{align}
\re(u(t,\cdot),B_\ell u(t,\cdot))_{L^2}&\geq \|S_0^{b_\ell }u(t,\cdot))\|^2_{L^2}+\re(E_{b_\ell}u(t,\cdot),u(t,\cdot))_{L^2}-C\|u(t,\cdot)\|_{L^2}^2\gtrsim -C_{\ell,T}Q_0(u_0)^2,\label{eqbell}\\
\re(u(t,\cdot),D_\ell u(t,\cdot))_{L^2}&\geq \|S_0^{d_\ell }u(t,\cdot))\|^2_{L^2}+\re(E_{d_\ell}u(t,\cdot),u(t,\cdot))_{L^2}-C\|u(t,\cdot)\|_{L^2}^2\gtrsim -C_{\ell,T}Q_0(u_0)^2, \label{eqdell}\\
\re(u(t,\cdot),C_\ell  u(t,\cdot))_{L^2}&\geq \|S_0^{c_\ell }u(t,\cdot))\|^2_{L^2}+\re(E_{c_\ell}u(t,\cdot),u(t,\cdot))_{L^2}-C\|u(t,\cdot)\|_{L^2}^2\gtrsim -C_{\ell,T}Q_0(u_0)^2,\label{eqcell}
\end{align}
where in the final inequalities we have used Lemma \ref{lemmaCKS} which is possible since
$E_{b_\ell},E_{c_\ell}, E_{d_\ell}\in 
S(0,2(\ell-1)r;\rho,\delta)
$ and have support included in that of $b_{\ell-1}$. Since $\re(u,B_{j}u)$, for all $j\leq \ell-1$, is bounded by the inductive hypotheses, we have the above bounds.
 
Having the previous inequalities satisfied, to conclude the inductive step as before we only have to check that the integral term in the RHS is bounded.

Due to the form of $E_\ell=\Op(e_\ell)$, see \eqref{eq.EnergyEstimate}, we have $e_\ell\in S(0, 2r\ell-2\min\{\delta,1-\rho\};\rho,\delta)$ and $\mathrm{supp}(e_\ell)\subseteq \mathscr{E}_+^{\ell-1}$. Since $\min\{\delta,1-\rho\}\geq \delta-\rho\geq r$, we also have $e_\ell\in S(0, 2r(\ell-1);\rho,\delta)$.
Therefore, by the inductive hypotheses and the properties of $e_\ell$, we can apply Lemma \ref{lemmaCKS} and get  
\[\re(\psi,E_\ell \psi)\leq \sum_{i=1}^{L-1}c_i\re(\psi,B_i \psi)_{L^2}+C_0\|\psi\|^2_{L^2}<\infty.\]
Moreover, by the inductive hypothesis
\[ \sup_{0\leq t\leq T}\Big(\re(u(t,\cdot),B_ju(t,\cdot))_{L^2}\Big)<\infty\quad \forall j<L,\]
so, using \eqref{eq.boundL2}, we also have 
\[ \sup_{0\leq t\leq T}\Big(\re(u(t,\cdot),E_\ell u(t,\cdot))_{L^2}\Big)<\infty,\]
which yields
\[\int_0^T\re(u(\tau,\cdot),E_\ell u(\tau,\cdot))_{L^2}d\tau <\infty.\]
Applying this into \eqref{eq.EEbl2} we obtain the result for $j=\ell$, which concludes the proof by induction.
\endproof

We are now in a position to prove the main theorem of this work, that is, Theorem \ref{thm.mainest}, which is the building block of the proof of Theorem \ref{theo.WFWFs}.
To prove Theorem \ref{thm.mainest} we use the following lemma taken from \cite{CKS}.

\begin{lemma}[Lemma 6.5 in \cite{CKS}]\label{lemmaCKS2.2}
   Assume that there is a sequence of symbols $c_j(x,\xi) \in S ( j , k - j ;
\rho , \delta)$ for $j = 1 ,\ldots,\ell<k$ such that 
\begin{itemize}
    \item [(i)] $c_j(x,\xi)\geq 0$;
    \item[(ii)] $\mathrm{supp}(c_j)\subseteq \{(x,\xi):c_{j-1}(x,\xi)\geq \la x\ra^{k-j+1}\la \xi\ra^{j-1} \}$ for $j>1$.
\end{itemize}
Consider a symbol $e\in S(\ell,k-\ell;\rho,\delta)$ with $\mathrm{supp}(e)\subseteq \{(x,\xi):c_{\ell}(x,\xi)\geq \la x\ra^{k-\ell}\la \xi\ra^{\ell}\}$. 
Then there exist $e_0\in S(0,k;\rho,\delta)$, with $\mathrm{supp}(e_0)\subseteq\mathrm{supp}(c_1)$ and some constants  $M_0,\dots,M_\ell>0$ such that
\begin{equation}\label{eq.lemW0}
    \re(\psi,\Op(e)\psi)_{L^2}\leq \sum_{j=1}^\ell M_j \re(\psi,\Op(c_j)\psi)_{L^2}+\re(\psi,\Op(e_0)\psi)_{L^2}+M_0\|\psi\|_{L^2}^2, \quad \psi \in \mathscr{S}(\R^n).
\end{equation}
\end{lemma}
A natural generalization of Lemma \ref{lemmaCKS} is the following, which is the version we will use in the proof. 
\begin{lemma}\label{lemmaCKS2}
   Let $\ell \geq 1$ and assume that there is a sequence of symbols $c_j(x,\xi) \in S ( m_j , k_j ;
\rho , \delta)$ with $j = 1 ,\ldots,\ell$, $0 \leq m_1 \leq 2$,  $0<m_{j+1}-m_j\leq 2$, and $k_{j+1}\leq k_{j}$ for $j=1,\dots,\ell-1$, such that 
\begin{itemize}
    \item [(i)] $c_j(x,\xi)\geq 0$;
    \item[(ii)] $\mathrm{supp}(c_j)\subseteq \{(x,\xi):c_{j-1}(x,\xi)\geq \la x\ra^{k_{j-1}}\la \xi\ra^{m_{j-1}} \}$ for $j\geq 2$.
\end{itemize}
Consider a symbol $e\in S(m_\ell,k_\ell;\rho,\delta)$ with $\mathrm{supp}(e)\subseteq \{(x,\xi):c_{\ell}(x,\xi)\geq \la x\ra^{k_\ell}\la \xi\ra^{m_\ell}\}$. 
Then there is a symbol $e_0\in S(0,k_1-2(\delta-\rho);\rho,\delta)$, with $\mathrm{supp}(e_0)\subseteq\mathrm{supp}(c_1)$ and some constants $M_1,\dots,M_\ell>0$ such that
\begin{equation}\label{eq.lemW1}
    \re(\psi,\Op(e)\psi)_{L^2}\leq \sum_{j=1}^\ell M_j  \re(\psi,\Op(c_j)\psi)_{L^2}+\re(\psi,\Op(e_0)\psi)_{L^2}+M_0\|\psi\|_{L^2}^2, \quad  \psi \in \mathscr{S}(\R^n).
\end{equation}

\end{lemma}

The proof follows the inductive argument of
Lemma 6.5 in \cite{CKS}, using Theorems \ref{thm.GardingInequality} and \ref{thm.UpperBound}
and Remark \ref{rmk.theo2.7}.
At each level $j\geq 2$, the remainder belongs to
\[
S(m_j-2,k_j-2(\delta-\rho);\rho,\delta)
\subseteq S(m_{j-1},k_{j-1};\rho,\delta),
\]
and its support is contained in the elliptic region
of $c_{j-1}$ by assumption (ii).
At the final level, $m_1\leq 2$ yields a remainder in
$S(0,k_1-2(\delta-\rho);\rho,\delta)$ supported in
$\operatorname{supp}(c_1)$.

\begin{proof}[Proof of Theorem \ref{thm.mainest}]
  As in Theorem \ref{thm.nestedOperators} we assume $u_0\in\mathscr{S}(\R^n)$. As before, the conclusion under the  assumption $u_0 \in L^2(\R^n)$ follows by the density argument described in Remark \ref{rmk.approx.argument}.
  
    When $K=0$ the statement amounts to the existence of a pair $(b_0,c_0)$ with suitable sign and support properties and such that (i) and (ii) are satisfied. This pair is the one in \eqref{b0c0} whose existence  is given by the results in Subsection \ref{subsec.bc}, so we have that (i) is trivially satisfied since $b_0\in S(0,0;\rho,\delta)$, while (ii) amounts to
    \[\re\int_0^T(u(t,\cdot),C_0u(t,\cdot))_{L^2}dt,\]
     which can be proved using the arguments in the proof of Theorem \ref{thm.nestedOperators} as follows. First, we apply the energy estimate \eqref{eq.EnergyEstimate} with $B=B_0=\Op(b_0)$, $C=C_0=\Op(c_0)$ and $D$ correspondingly defined. Then we integrate in time the resulting inequality and obtain
     \begin{align}
\re &(u(t,\cdot),B_0u(t,\cdot))_{L^2}+\frac 1 2\int_0^t \re (u(\tau,\cdot), C_0u(\tau,\cdot))_{L^2}d\tau +\frac{1}{2}\int_0^t \re (u(\tau,\cdot), D_0u(\tau,\cdot))_{L^2}d\tau \\
& \lesssim \re (u_0,B_0u_0)_{L^2}+\int_0^t \re (u(\tau,\cdot), E_0u(\tau,\cdot))_{L^2}d\tau+ \|u_0\|_{L^2}^2,         
     \end{align}
 with $E_0$ suitably defined. Finally, by the $L^2$-boundedness of $E_0,B_0$, the LHS is bounded, while the RHS has a lower bound by the Gårding inequality in Theorem \ref{thm.GardingInequality}. Putting everything together, we conclude the desired bound. 
\medskip

We now consider the case where $K>0$ is an integer.

Let $b_j^0(x,\xi)$ be the sequence of symbols in Theorem \ref{thm.nestedOperators}, that is such that
$b_j^0\in S(0,2j r; \rho,\delta)$,
and let $\mathscr{E}_+^j$, $j\leq J$, be the corresponding sequence of nested subneghborhoods of $(x_0,\xi_0)$ in Theorem \ref{thm.nestedOperators} satisfying
 \[
    \mathscr{E}_{+}^J(x_0,\xi_0) \subseteq \mathscr{E}_+^{J-1}(x_0,\xi_0)\subseteq \dots \subseteq \mathscr{E}_+^{0}(x_0,\xi_0)
    \]
  where
    \[
    \mathscr{E}_{+}^\ell(x_0,\xi_0)=\lbrace (x,\xi); \ b_\ell(x,\xi) \geq \langle x \rangle^{2\ell r} \rbrace , \quad \ell=1,\dots,J,
    \]
and where
    \[
     \mathrm{supp}(b_j^0)\subseteq \mathscr{E}_{+}^{j-1}(x_0,\xi_0).
    \] 
The sequence $b_j^0(x,\xi)$, by Theorem \ref{thm.nestedOperators}, satisfies
\begin{equation}
    \sup_{0\leq t\leq T}\re(u(t,\cdot),\Op(b_j^0)u(t,\cdot))_{L^2}<\infty, \quad \forall j\leq J,
\end{equation}
for any fixed $T$.
Next, we use the sequence $\mathscr{E}_+^j$ (we shall omit the point $(x_0,\xi_0)$ from now on to make the notation shorter), to build the sequence $\mathscr{D}_+^j$. We set
\[\mathscr{D}_+^0\equiv \mathscr{E}_+^J\subseteq\mathscr{E}_+^0 \]
and proceed by induction to show the existence of the sequence
    \[\mathscr{D}_+^K\subseteq \mathscr{D}_+^{K-1}\subseteq \cdots \subseteq\mathscr{D}_+^0\equiv \mathscr{E}_+^J\]
    and of the symbols $(b_j,c_j)$ with the required property.

If $\ell=0, 0<k\leq K$, by Theorem \ref{thm.nestedOperators} there exist $\tilde{b}_0\in S(0,k;\rho,\delta), \tilde{c}_0\in S(2,k-1;\rho,\delta)$, with $\mathrm{supp}(\tilde{b}_0)\subseteq \mathscr{D}_+^0\equiv\mathscr{E}_+^J$. Therefore, an application of \eqref{eq.EnergyEstimate} with $B=\tilde B_0:=\Op(\tilde{b}_0), C=\tilde C_0:=\Op(\tilde{c}_0)$ gives
\begin{equation}\label{eq.EEbl0}
\begin{split}
\re &(u(t,\cdot), \tilde B_0 u(t,\cdot))_{L^2}+\frac 1 2 \int_0^t \re (u(\tau,\cdot),  \tilde C_0 u(\tau,\cdot))_{L^2}d\tau+\frac 1 2 \int_0^t \re (u(\tau,\cdot), \tau^\ell \tilde D_0 u(\tau,\cdot))_{L^2}d\tau \\
& \lesssim  \re (u_0,\tilde B_0 u_0)_{L^2}+\int_0^t  \re (u(\tau,\cdot), \tilde E_0 u(\tau,\cdot))_{L^2}d\tau+ \|u_0\|_{L^2}^2.
\end{split}
\end{equation}
Using that $\mathrm{supp}(\tilde{b}_0), \mathrm{supp}(e_0)\subseteq \mathscr{D}_+^0\equiv\mathscr{E}_+^J\subseteq \mathscr{E}_+^0$ 
we have, for $b_0$ as in Theorem \ref{thm.nestedOperators},
\begin{align}
    |\re (u_0,\tilde B_0 u_0)_{L^2}|&\lesssim \re(u_0,S_0^*S_0 u_0)_{L^2}+C\|u_0\|_{L^2}^2,\quad s_0=\sqrt{b_0}, \label{term1}\\
    \re (u(\tau,\cdot), \tilde E_0 u(\tau,\cdot))_{L^2} &\lesssim  \sum_{j=1}^J\re(u(\tau,\cdot), \Op(b_j^0)u(\tau,\cdot))_{L^2}+C\|u(t,\cdot)\|_{L^2}^2,\label{term2}
\end{align}
where the first inequality is due to Theorem \ref{thm.UpperBound} and the second to Lemma \ref{lemmaCKS}. 
In particular, \eqref{term1} is bounded by the moment condition \eqref{m.moment.condition}, while the time integral of the LHS of \eqref{term2} is bounded by Theorem \ref{thm.nestedOperators} applied on the sequence $b_j^0$ and the bound $\sup_{0\leq t\leq T}\|u(t,\cdot)\|_{L^2}\lesssim_T\|u_0\|_{L^2}$.
Hence, the RHS of \eqref{eq.EEbl0} is bounded. The LHS of \eqref{eq.EEbl0} is bounded from below by the Gårding, hence the the desired conclusion is true when $\ell=0$ and $0<k\leq K$.

We next proceed by induction in $\ell$, so we take $k\leq K$ and assume the inductive hypothesis that the result is true for all $0\leq j\leq \ell-1$, where $0< \ell \leq k\leq K$.
Hence, by the inductive hypothesis,  for each $j\leq \ell-1$ there exist a neighborhood
\[\mathscr{D}_+^j\subseteq \mathscr{D}_+^{j-1}\subseteq \cdots \subseteq\mathscr{D}_+^0\equiv \mathscr{E}_+^J\]
and symbols $b_j,c_j$ such that
\begin{align}
    0\leq b_j(x,\xi)&\in S(j,k-j;\rho,\delta)\\
    0\leq c_j(x,\xi)&\in S(j+2,k-j-1;\rho,\delta)
\end{align}
and
\[\mathrm{supp}(c_j)\subseteq \mathrm{supp}(b_j)\subseteq\mathscr{D}_+^{j-1}:=\{(x,\xi): c_{j-1}(x,\xi)\geq \la x \ra^{k-j} \la \xi\ra^{j+1}\}.\]
We now complete the sequence of symbols with a pair $(b_\ell,c_\ell)$
using the construction of Theorem \ref{theo.nested}. 
The sequences $b_j,c_j$, $j\leq \ell-1$, of the inductive hypothesis satisfy Theorem \ref{theo.nested}  with $m_j=j,k_j=k-j$, so we build $b_\ell,c_\ell$ from $b_{\ell-1},c_{\ell-1}$. Then, by construction,  $\mathrm{supp}(c_\ell)\subseteq \mathscr{D}_+^{\ell-1}$, $\mathscr{D}_+^\ell=\{(x,\xi): c_\ell(x,\xi)\geq \la x\ra^{k-\ell-1}\la \xi\ra^{\ell+2}\}$ when $\ell\neq k$, and  
\[\mathscr{D}^k_+=\{(x,\xi): c_k(x,\xi)\geq \la x\ra^{-k'}\la \xi\ra^{k+2}\}\neq \emptyset,\]
when $\ell=k$. This exception in the definition of $\mathscr{D}_+^k$ is due to the fact that, by Theorem \ref{theo.nested}, when $k_j=0$ 
\begin{equation}\label{eq,l=k}
    c_{k_j}(x,\xi)\gtrsim \la x\ra ^{-k'}\la \xi\ra^{m_j+2}, \quad \forall(x,\xi)\in\mathscr{E}_+^0(x_0,\xi_0),\end{equation}
for all $k'>1$.
Therefore, since for $k=\ell$ we have $k_\ell=0$ and the definition of $\mathscr{D}_+^K$ has to take a different form.
 
Now we apply the energy estimate \eqref{eq.EnergyEstimate} with $B=t^\ell B_\ell, C=t^\ell C_\ell,$ and $D=t^\ell D_\ell$, integrate in time, and get, since $\ell>0$,
\begin{equation}\label{eq.EEbl}
\begin{split}
\re &(u(t,\cdot),t^\ell B_\ell u(t,\cdot))_{L^2}+\frac 1 2 \int_0^t \re (u(\tau,\cdot), \tau^\ell C_\ell u(\tau,\cdot))_{L^2}d\tau+\frac 1 2 \int_0^t \re (u(\tau,\cdot), \tau^\ell D_\ell u(\tau,\cdot))_{L^2}d\tau \\
& \lesssim \int_0^t \ell \re (u(\tau,\cdot),\tau^{\ell-1}B_\ell u(\tau,\cdot))_{L^2}d\tau+\int_0^t \tau^\ell \re (u(\tau,\cdot), E_\ell u(\tau,\cdot))_{L^2}d\tau+ \|u_0\|_{L^2}^2.
\end{split}
\end{equation}
Our inductive hypothesis gives the validity of the bounds (i) and (ii) with $B_j, C_j$, for all $j<\ell$, respectively, so use  this fact in combination with Lemma \ref{lemmaCKS2} to prove (i) and (ii) for $B_j, C_j$ with $j=\ell$, which would conclude the proof by induction.

Note that $E_\ell=\Op(e_\ell)$ with $e_\ell\in S(\ell,k-\ell-2\delta;\rho,\delta)\subset S(\ell,k-(\ell-1)-1;\rho,\delta)\subset S(\ell+1,k-(\ell-1)-1;\rho,\delta)$ and $\mathrm{supp}(e_\ell)\subseteq \mathscr{D}_+^{\ell-1}=\{ (x,\xi): c_{\ell-1}(x,\xi)\gtrsim \la x\ra^{k-\ell}\la \xi\ra^{\ell+1}\}$. Hence, by Lemma \ref{lemmaCKS2} (with $c_j\in S(j+2,k-j-1)$),
\begin{align}
  \tau^\ell \re (u(\tau,\cdot), E_\ell u(\tau,\cdot))_{L^2}&\leq \tau^\ell\sum_{j=0}^{\ell-1} C_j\re (u(\tau,\cdot), C_j u(\tau,\cdot))_{L^2}\\
  &+\tau^\ell \re (u(\tau,\cdot), \Op(e_0) u(\tau,\cdot)), u(\tau,\cdot))_{L^2}+C\|u(t,\cdot)\|^2_{L^2},\label{eq.est.el} 
\end{align}
where the term containing $\Op(e_0)$, due to the condition $\mathrm{supp}(e_0)\subseteq \mathscr{D}_+^0\equiv\mathscr{E}_+^J$, by Lemma \ref{lemmaCKS} satisfies
\begin{equation}
 \tau^\ell \re (u(\tau,\cdot), \Op(e_0) u(\tau,\cdot))_{L^2}\leq T^\ell  \sum_{j=1}^J\re (u(\tau,\cdot), \Op(b_0^j)u(\tau,\cdot))_{L^2} +C_0\|u(\tau,\cdot))\|^2_{L^2}<+\infty. \label{eq.est.e0}
\end{equation}
Then, after integration in time in \eqref{eq.est.el}, using \eqref{eq.est.e0}, the inductive hypothesis for $j\leq \ell-1$, and the property $\sup_{0\leq t\leq T}\|u\|_{L^2}\leq C_T \|u_0\|_{L^2}$, we obtain
\begin{equation}
    \int_0^T\tau^\ell \re (u(\tau,\cdot), E_\ell u(\tau,\cdot))_{L^2}<+\infty.
\end{equation}
To conclude the boundedness of the RHS of \eqref{eq.EEbl}, we just need to  bound the integral containing $B_\ell$ in the RHS of \eqref{eq.EEbl}, which can be done by repeating the same steps we have just performed to get the boundedness of the integral containing $E_\ell$, so we do not rewrite. 
The bound from below can be obtained as in the proof of Theorem \ref{thm.nestedOperators}, more specifically using Theorem \ref{thm.GardingInequality} and the inductive hypothesis to bound the reminder terms (here we use the properties of support, see the case $K$ non integer below for more details). Hence, due to the upper and lower bounds, we get the result for $j=\ell$, which gives the result when $K>0$ is an integer.

We then focus on the case when $K$ is not an integer.
Here the delicate case is the first inductive step, that is when $\ell=0$. This amounts to construct symbols $\tilde b_0,\tilde c_0$ such that $\tilde b_0\in S(\kappa,K-\kappa;\rho,\delta)$ and $\tilde c_0\in S(2+\kappa,K-\kappa-1;\rho,\delta)$ satisfying (i), (ii). The construction is given by Theorem \ref{theo.nested}. However, for the bounds we cannot exploit the $L^2$-boundedness properties that one has when $\kappa=0$ (see Theorem \ref{thm.L2boundedness}) used in the case $K=0$ integer, so we proceed as follows.

We then focus on the case $\ell=0$.
Note that $\tilde{b}_0,\tilde{c}_0$ exist by Theorem \ref{theo.nested}, we need to check that (i) and (ii) hold.
Arguing as in the first part of the proof (specifically where we constructed $\tilde{b}_0\in S(0,k;\rho,\delta)$ and $\tilde{c}_0\in S(2,k;\rho,\delta)$ also used in Theorem \ref{thm.nestedOperators}, with $k\leq K$; note that this pair is the old one, not the one we are working with now) we can construct a new symbol pair $({b}'_0,{c}'_0)$ such that $b'_0\in S(0,K;\rho,\delta)$ and $c'_0\in S(2,K-1;\rho,\delta)$, satisfying (i) and (ii) respectively (this can be proved via the energy estimate, some a priori estimates as that in Lemma \ref{lemmaCKS} and \ref{lemmaCKS2} and support properties),  and such that
\[\mathrm{supp}(\tilde b_0)\subseteq \{(x,\xi): b'_0(x,\xi)\gtrsim \la x\ra^K\}\cap \{(x,\xi): c'_0(x,\xi)\gtrsim \la x\ra^{K-1} \la \xi\ra^2\}\subseteq \mathscr E_+^0,\]
where, recall, $b_0$ is as in \eqref{b0c0}.
In addition, we consider
\[
\operatorname{supp}c'_0
\subset\operatorname{supp}b'_0
\subset 
\mathscr{E}_+^J
=\{(x,\xi):b^0_J(x,\xi)\geq\langle x\rangle^K\}.
\]

Then we apply the energy estimate to $\tilde B_0,\tilde C_0$ and get
\begin{equation}
\begin{split}
\re &(u(t,\cdot),T^{\mu_0} \tilde B_0 u(t,\cdot))_{L^2}+\frac 1 2 \int_0^T \re (u(\tau,\cdot), \tau^{\mu_0} \tilde C_0 u(\tau,\cdot))_{L^2}d\tau+\frac 1 2 \int_0^T \re (u(\tau,\cdot), \tau^{\mu_0} \tilde D_0 u(\tau,\cdot))_{L^2}d\tau \\
& \lesssim \int_0^T  \re (u(\tau,\cdot),\tau^{\mu_0-1}\tilde B_0 u(\tau,\cdot))_{L^2}d\tau+\int_0^T \tau^{\mu_0} \re (u(\tau,\cdot), \tilde E_0 u(\tau,\cdot))_{L^2}d\tau+ \|u_0\|_{L^2}^2.\label{eqn.final}
\end{split}
\end{equation}
Note that each term on the LHS is bounded from below. This is proved, as above, by using the Gårding inequality and the fact that the additional reminder terms are given in terms of operators that can be bounded using Lemma \ref{lemmaCKS}, Lemma \ref{lemmaCKS2} (applied to $e$ and $-e$) and the moment condition. For instance,
\begin{align}
    \int_0^T \tau^{\mu_0}\re (u(\tau,\cdot), \tilde C_0 u(\tau,\cdot))_{L^2}\gtrsim -C_T\int_0^T|\re(\Op(e) u(\tau,\cdot),u(\tau,\cdot))|d\tau-C\|u(\tau,\cdot)\|^2_{L^2}
\end{align}
where, by Theorem \ref{thm.GardingInequality} $e\in S(\kappa,K-\kappa-1-2(\delta-\rho);\rho,\delta)$ with $\mathrm{supp}(e)\subseteq \mathrm{supp}(\tilde{c}_0)\subseteq \{ (x,\xi):c'_0(x,\xi)\gtrsim \la x\ra^{K-1}\la \xi\ra^2\}$. By applying Lemma \ref{lemmaCKS2} with $c_j\in S(3-j, K-j;\rho,\delta)$, $j=1$, so that $c_1=c_0'$,
we conclude
\begin{align}
   \int_0^T \tau^{\mu_0}\re (u(\tau,\cdot), \tilde C_0 u(\tau,\cdot))_{L^2}\gtrsim & -C_T\int_0^T\re(\Op(c_0') u(\tau,\cdot),u(\tau,\cdot))d\tau\\
   &- C_T\int_0^T\re(\Op(e_0') u(\tau,\cdot),u(\tau,\cdot))d\tau-C\|u(\tau,\cdot)\|^2_{L^2}\\
   \gtrsim &  - C_T(\|u_0\|^2+\|S_0u_0\|^2) >-\infty 
\end{align}
which is bounded by the construction of $c_0'$ and the properties of $e_0'$. In fact, $e_0'\in S(0,K-1;\rho,\delta)$ with $\mathrm{supp}(e_0')\subseteq \mathscr{E}_+^J$, so the term in the second line above can be bounded as in \eqref{eq.est.e0}. The lower bound for the term in $\tilde{B}_0$ is analogous but uses Lemma \ref{lemmaCKS}.

Next, to estimate the RHS from above we observe that
\begin{align}
|\tilde b_0|,\ |\tilde e_0|
&\lesssim
{b'}_0^{\,1-\kappa}{c'}_0^{\kappa},\\
t^{\mu_0}|\tilde e_0|
&\lesssim
b'_0+
t^{\mu_0/\kappa}c'_0,\\
t^{\mu_0-1}|\tilde b_0|
&\lesssim
t^{-\frac{1-\mu_0}{1-\kappa}}b'_0
+c'_0\longrightarrow 
C b'_0+C t^{\frac{1-\mu_0}{1-\kappa}} c'_0-t^{\kappa \frac{1-\mu_0}{1-\kappa}}\widetilde b_0=:h_t\geq 0,
\end{align}
where we take $1>\mu_0>\kappa$ to ensure the local integrability of $t^{-\frac{1-\mu_0}{1-\kappa}}$ and that $C b'_0+C t^{\frac{1-\mu_0}{1-\kappa}} c'_0-t^{\kappa \frac{1-\mu_0}{1-\kappa}}\widetilde b_0\in S(2,K; \rho,\delta)$ uniformly in time.
Therefore, we can apply Theorem \ref{thm.GardingInequality} (the square root of the considered symbol is well defined and belongs to $S(1,K/2;\rho,\delta)$ uniformly in time) and get

\begin{align}
    \int_0^T\re (u(\tau,\cdot),\tau^{\mu_0-1}\tilde B_0 u(\tau,\cdot))_{L^2}&\lesssim \int_0^T\re (u(\tau,\cdot),\tau^{\frac{\mu_0-1}{1-\kappa}}{B'}_0 u(\tau,\cdot))_{L^2}+\int_0^T
    \re (u(\tau,\cdot),{C'}_0 u(\tau,\cdot))_{L^2} d\tau\\
    &+ \int_0^T \tau^{\frac{\mu_0-1}{1-\kappa}}\re(E_{0} u(\tau,\cdot),u(\tau,\cdot)d\tau +C_T \|u_0\|_{L^2}^2,
\end{align}
where $E_{0}=\Op(e_0)$ with $e_0\in S(0, K-2(\delta-\rho);\rho,\delta)$ with $\mathrm{supp}(e_0)\subseteq \mathrm{supp}(h_t)\subseteq \mathrm{supp}(b_0')\subseteq \mathscr{E}_+^J$ ($e_0$ may depend on time but belongs to $S(0, K-2(\delta-\rho);\rho,\delta)$ uniformly in time). 
Now, the first two terms on the RHS are bounded by the choice of $\mu_0$ and the fact that $C',B'$ satisfy the suitable bounds by construction. On the other hand, due to $\mathrm{supp}(e_0)\subseteq \mathscr{E}_+^J$, then estimate \eqref{eq.est.e0} holds in this case as well,
which yields the boundedness of the first term in the of \eqref{eqn.final}. Since the second term on the RHS of of \eqref{eqn.final} is estimated similarly, we omit the details.

Finally, the second and last inductive step works exactly as in the previous case, that is exploiting the inductive hypothesis.  This completes the proof.
\end{proof} 

We now use this theorem to prove our main microlocal results (see Theorem \ref{theo.WFs} below), and to this end we first give the following definition.

\begin{definition}\label{def.finitemoments}
Let $(x_0,\xi_0) \in \R^{2n} \setminus 0$ and let $K\geq 0$. 
We say that  \textit{$u_0$ has finite moment of order $K$ along the forward bicharacteristic starting at $(x_0,\xi_0)$}
if there exist  two nonnegative symbols 
\[
0 \leq b_0 \in S(0,2K;0,1), \quad 0 \leq c_0=\lbrace a_3,b_0 \rbrace \in S(2,2K-1;0,1),
\]
with $\mathrm{supp}(c_0)\subseteq \mathrm{supp}(b_0)$, 
satisfying 
\[
0 \leq \lbrace a_3, \langle x \rangle^{-2K}b_0 \rbrace (x,\xi), \quad \forall (x,\xi) \in \mathscr{E}^0_+(x_0,\xi_0)\cap \mathscr{S}^{+}, \quad  b_0(x,\xi) \geq \langle x \rangle^{2K} \quad \forall (x,\xi)\in\mathscr{E}_+^0(x_0,\xi_0),
\]
where 
\[
\mathscr{E}^0_{+}(x_0,\xi_0):=\big\lbrace (x,\xi); b_0(x,\xi) \geq \langle x \rangle^{2K} \rbrace \subseteq \mathscr{E}_+(x_0,\xi_0):=\mathrm{supp}(b_0),
\]
are conic neighborhoods of  the forward bicharacteristic $\lbrace \varphi (s;x_0,\xi_0); s\geq 0 \rbrace$ starting from $(x_0,\xi_0)$. 
Moreover, $b_0$
has a smooth square root $s_{0,K}=\sqrt{b_0} \in S(0,K;0,1)$ for which the following holds
\[
\|S_{0,K}u_0\|_{L^2}<+\infty,\quad \text{with}\,\,S_{0,K}=\Op(s_{0,K}). 
\]
\end{definition}

We are in a position to prove Theorem \ref{theo.WFs}. 

\begin{proof}[Proof of Theorem \ref{theo.WFs}]
Once again we assume $u_0\in\mathscr{S}(\R^n)$.  The conclusion under the  assumption $u_0 \in L^2(\R^n)$ follows by the approximation argument described in Remark \ref{rmk.approx.argument}.

Let $t \in (0,T]$ and let $s \geq 0$ be fixed and prove that $(x_0,\xi_0) \notin WF^s(u(t,\cdot))$. By Theorem \ref{thm.mainest} we can construct a sequence of symbol pairs $(b_\ell,c_\ell)_{\ell=0,\dots,\lfloor 2s \rfloor}$, with $0 \leq b_\ell \in S(\ell+2s-\lfloor 2s \rfloor;\lfloor 2s \rfloor - \ell;0,1)$, such that $B_\ell=\Op(b_\ell)$ satisfies 
    \[
    \re (u(t,\cdot),B_\ell u(t,\cdot))_{L^2} <+\infty.
    \] 
    Moreover, Theorem \ref{thm.mainest} provides also the existence of neighborhoods of the bicharacteristic through $(x_0,\xi_0)$, denoted by $\mathscr{D}^\ell_+$, where $b_\ell$ is elliptic and admits an elliptic square root $s_\ell=\sqrt{b_\ell}$, that is such, $S_\ell=\Op(s_\ell)$ is elliptic at $(x_0,\xi_0)$ (recall Definition \ref{def.elliptic}).

   To prove the theorem we proceed by induction on $\ell$ and show that 
    for all $\ell \in \lbrace 0, \dots, \lfloor 2s \rfloor \rbrace$ we have that $S_\ell=\Op(s_\ell)$ is elliptic at $(x_0,\xi_0)$ and satisfies
    \begin{equation}\label{eq.SN}
    \|S_\ell^\ast u(t,\cdot)\|_{L^2}^2 <+\infty. 
    \end{equation}
    Indeed, if the latter holds for $\ell=\lfloor 2s \rfloor$, we get that $(x_0,\xi_0) \notin WF^s(u(t,\cdot))$.
    To establish the induction argument the purpose is to combine the estimates obtained in Theorem \ref{thm.mainest} with the G\aa rding inequality \ref{thm.GardingInequality}  as follows. 

    If $\ell=0$, by Theorem \ref{thm.mainest} there exists $b_0 \in S(2s-\lfloor 2s \rfloor,\lfloor 2s \rfloor;0,1)$, that admits a smooth square root whose corresponding operators $s_{0}:=\sqrt{b_{0}}$ $S_0=\Op(s_0)$ is elliptic at $(x_0,\xi_0)$, and  
    \begin{equation}\label{eq.B2}
    \re (u(t,\cdot), B_{0}u(t,\cdot))_{L^2}< +\infty.
    \end{equation}
    Moreover, $b_0$ can be taken as in the proof of Theorem \ref{thm.mainest} , that is with 
\[\mathrm{supp}(b_0)\subseteq \{(x,\xi): b'_0(x,\xi)\gtrsim \la x\ra^K\}\cap \{(x,\xi): c'_0(x,\xi)\gtrsim \la x\ra^{K-1} \la \xi\ra^2\}\subseteq \mathscr E_+^0,\]
where
\[
\operatorname{supp}c'_0
\subset\operatorname{supp}b'_0
\subset 
\mathscr{E}_+^J
=\{(x,\xi):b^0_J(x,\xi)\geq\langle x\rangle^K\},
\]
where $b_0'\in S(0,2s;0,1) $, $c_0'\in S(2,2s-1;0,1)$.

    Hence, by Theorem \ref{thm.GardingInequality}, there exists $e_{0} \in S(0,\lfloor 2s \rfloor-2; 0,1)$ with $ \mathrm{supp}(e_{0})\subseteq \mathrm{supp}(b_0)\subset \mathscr{E}_+^J$ such that, on denoting $E_{0}=\Op(e_0)$ we have 
    \[
     \|S_0^\ast u(t,\cdot)\|_{L^2}^2 \lesssim \re (u(t,\cdot), B_{0} u(t,\cdot))_{L^2}+|\re(u(t,\cdot),E_{0}u(t,\cdot))_{L^2}|+ \|u(t,\cdot)\|_{L^2}^2.
    \]
    Due to the inclusion $ \mathrm{supp}(e_{0})\subseteq \mathrm{supp}(b_0)\subset \mathscr{E}_+^J$, the term on the RHS containing $E_0$ can be bounded as in \eqref{eq.est.e0}. Then,
    using the boundedness of the latter and \eqref{eq.boundL2} into by \eqref{eq.B2}, we get  
    \begin{align}
    \|S_0^\ast u(t,\cdot)\|_{L^2}^2 \lesssim \re  (u(t,\cdot), B_{0}(u(t,\cdot))_{L^2}+\|u_0\|_{L^2}^2 + \|S_0u_0\|_{L^2}^2 < +\infty,
    \end{align}
    which concludes the first inductive step.
    
    For $\ell\geq 1$,  let us consider the sequence $(b_\ell,c_\ell)_{\ell=0,\dots, \lfloor 2s \rfloor}$ constructed in Theorem \ref{thm.mainest}. By the inductive hypothesis for $s_{\ell-1}=\sqrt{b_{\ell-1}}$, and $S_{\ell-1}=\Op(s_{\ell-1})$, we have
    \begin{equation}\label{eq.SN-1}
    \|S_{\ell-1}^\ast u(t,\cdot)\|_{L^2}^2 <+\infty.
    \end{equation}
    Furthermore, by construction,
    \[
    \mathrm{supp}(b_\ell) \subseteq \lbrace (x,\xi); b_{\ell-1} \geq \langle x \rangle^{\lfloor 2s \rfloor-\ell+1}\langle \xi \rangle^{2s-\lfloor 2s \rfloor +\ell-1}   \rbrace
    \]
    and
    \begin{equation}\label{eq.B2N}
    \re (u(t,\cdot),B_{\ell}u(t,\cdot))_{L^2}<+\infty.
    \end{equation}
    Next, an application of Theorem \ref{thm.GardingInequality} gives that there exists $e_{\ell} \in S(\ell+2s-\lfloor 2s \rfloor-2,\lfloor 2s \rfloor-\ell-2; 0,1)$ with $ \mathrm{supp}(e_{\ell})\subseteq \mathrm{supp}(b_{\ell})$ such that, on denoting $E_{\ell}=\Op(e_{\ell})$,
    \[
     \|S_\ell^\ast u(t,\cdot)\|_{L^2}^2 \lesssim \re (u(t,\cdot), B_\ell u(t,\cdot))_{L^2}+|\re(u(t,\cdot),E_{\ell}u(t,\cdot))_{L^2}|+ \|u(t,\cdot)\|_{L^2}^2.
    \]
    Since $\mathrm{supp}(e_\ell) \subseteq \lbrace (x,\xi); b_{\ell-1} \geq  \langle x \rangle^{\lfloor 2s \rfloor-\ell+1}\langle \xi \rangle^{2s-\lfloor 2s \rfloor +\ell-1}\rbrace$, by Theorem \ref{thm.UpperBound} (see Remark \ref{rmk.theo2.7}) we get that
    \[
    |\re(u(t,\cdot),E_{\ell}u(t,\cdot))_{L^2}| \lesssim \|S_{\ell-1}^\ast u(t,\cdot)\|_{L^2}^2+\|u(t,\cdot)\|^2_{L^2},
    \]
    and so 
    \[
    \|S_\ell^\ast u(t,\cdot)\|_{L^2}^2 \lesssim \re (u(t,\cdot), B_{\ell} u(t,\cdot))_{L^2}+\|S_{\ell-1}^\ast u(t,\cdot)\|_{L^2}^2+ \|u(t,\cdot)\|_{L^2}^2.
    \]
    Therefore, by \eqref{eq.B2N}, \eqref{eq.SN-1} and \eqref{eq.boundL2} we obtain 
    \[
    \|S_\ell^\ast u(t,\cdot)\|_{L^2}^2<+\infty.
    \]
\end{proof}
This theorem shows that if the initial datum $u_0$ has finite moments of order $s \in \R$ along the forward bicharacteristic starting from the point $(x_0,\xi_0)$, then the corresponding solution $u$ to \eqref{eq.mainprob} is microlocally $H^s$ at $(x_0,\xi_0)$, namely, $u$ is $H^s$ near $x_0$ in the direction $\xi_0$.


Finally, due to Theorem \ref{theo.singsupp}, we obtain the proof of Theorem \ref{theo.glob}, which is our smoothing result in the global setting. 
\begin{proof}[Proof of Theorem \ref{theo.glob}]
    The proof follows by combining Theorem \ref{theo.WFs} first with Theorem \ref{theo.Hs} and then \ref{theo.singsupp}. Indeed, if \eqref{eq.globalmomentK} holds, we have that the moment condition of Theorem \ref{theo.WFs} is satisfied at order $s$ at each $(x_0,\xi_0) \in \R^{2n}\setminus 0$, and so $WF^s(u(t,\cdot))=\emptyset$ for all $t \in (0,T]$, which implies $u \in H_{\mathrm{loc}}^s(\R^n)$ by Theorem \ref{theo.Hs}. 
    
    In particular, if \eqref{eq.globalmoment} holds (namely the moment condition at each order) then we have that $WF^s(u(t,\cdot))=\emptyset$, for each $s \in [0,+\infty)$, and so by Theorem \ref{theo.WFWFs} (see also Remark \ref{rmk.WFK}) we have that  $WF(u(t,\cdot))=\emptyset$ for all $t \in (0,T]$. Hence, by Theorem \ref{theo.singsupp} we have that $\mathrm{singsupp}(u(t,\cdot))=\emptyset$ for all $t \in (0,T]$ which in turn yields $u(t,\cdot) \in C^\infty(\R^n)$, for all $t \in (0,T]$.
\end{proof}

\section{Future perspectives}\label{sec.futper}
In this work, we have established smoothing results from different
perspectives for ZK-type \emph{linear} Cauchy problems of the form
\eqref{eq.mainprob}. A natural further question is whether these results
can be extended to the inhomogeneous or nonlinear setting
(namely, when the right-hand side of the first equation in
\eqref{eq.mainprob} contains a forcing term $f$ or a nonlinear
term $N(u)$), and whether they can be used to establish
well-posedness in the nonlinear case. We believe that these
questions can be addressed by combining our results with a
paradifferential calculus adapted to the problem, following
the approach developed for Schr\"odinger-type equations
in \cite{Sz}. We plan to investigate these questions in future work.

\sectionnotoc{AI declaration}
Most of this paper was written between September 2025 and May 2026, without the use of any AI tools. During the final stages of preparation, we used ChatGPT-6 Astra to revise the manuscript and simplify some arguments. We take full responsibility for the content of this work.

\appendix
\section{Auxiliary results}

\begin{proposition}\label{eq.usefineq} 
If $\frac 1 C|\xi|\leq |\Xi(s;x,\xi)|\leq C|\xi|$, for $C>0$ we have
\begin{align}
    |\p_{x_k}\Xi_i\Xi_j|\leq f(x,\xi)\quad \forall i,j=1,\ldots,n \Longrightarrow |\p_{x_k}\Xi_i|\leq f(x,\xi)/|\Xi|\lesssim f(x,\xi)/|\Xi| \quad \forall i=1,\dots,n.
\end{align}
\end{proposition}
\proof
Note that, when $i=j$, our hypothesis gives
\begin{align}
    |\p_{x_k}(\Xi_j\Xi_j)|=2|\Xi_j\p_{x_k}\Xi_j|\leq 2 f(x,\xi).
\end{align}
Therefore,
\begin{align}
    |\Xi|^2|\p_{x_k}\Xi_i|&=\sum_j |\Xi_j^2 \p_{x_k}\Xi_i|\\
    &=\sum_j |\Xi_j (\Xi_j\p_{x_k}\Xi_i)|\\
    &=\sum_j |\Xi_j \p_{x_k}(\Xi_j\Xi_i)-\Xi_i(\Xi_j \p_{x_k}\Xi_j)|\\
    &\leq \sum_j |\Xi_j| |\p_{x_k}(\Xi_j\Xi_i)|+\sum_j|\Xi_i| |\Xi_j \p_{x_k}\Xi_j|\\
    &\leq C(n) |\xi|f(x,\xi).
\end{align}
Hence,
$$ |\p_{x_k}\Xi_i|\lesssim f(x,\xi)/|\Xi|\lesssim f(x,\xi)/|\xi|,\quad \forall i=1,\ldots, n,$$
which proves the claim.
\endproof
\begin{proposition} \label{prop.App.2}
    Let $(x,\xi)\in\mathscr{S}_+$, then
    \begin{align}
        &|\widebar{\Xi}(x,\xi)|\approx |\xi|\\
        &\la \widebar{X}(x,\xi)\ra\approx \la x\ra.
    \end{align}
\end{proposition}

\proof
By Lemma \ref{lem.scatteringV} we have 
\begin{align}
    |\widebar{X}(x,\xi)-x|\leq C \la x\ra^{-2}
\end{align}
which yields
\begin{align}
   C_1^{-1}\la x\ra\leq|x|+1\leq |\widebar{X}(x,\xi)|+  C \la x\ra^{-2}+1\leq (C+1)C_1 \la \widebar{X}(x,\xi)\ra 
\end{align}
and
\begin{align}
C_1^{-1}\la \widebar{X}(x,\xi)\ra \leq    |\widebar{X}(x,\xi)|+1\leq  |x|+C\la x\ra^{-2}+1\leq (C+1)C_1 \la x\ra.
\end{align}
So we have
\begin{align}
\Big(\frac{C_1(C+1)}{C_1^{-1}}\Big)^{-1} \la \widebar{X}(x,\xi)\ra  \leq \la x\ra \leq \frac{C_1(C+1)}{C_1^{-1}}\la \widebar{X}(x,\xi)\ra,
\end{align}
proving one part of the statement.

To prove the other estimate we proceed similarly. By 
Lemma \ref{lem.scatteringV} we have 
\begin{align}
    |\widebar{\Xi}(x,\xi)-\Xi(s;x,\xi)|\leq C \la X(s;x,\xi)\ra^{-2}.
\end{align}
Moreover, by Lemma \ref{st.lem1} we have $|\Xi(s;x,\xi)|\approx |\xi|$ when $(x,\xi)\in\mathscr{S}_+$. Therefore, 
\begin{align}
-\frac{C}{\la X(s;x,\xi)\ra}+C_1^{-1}|\xi|\leq  -\frac{C}{\la X(s;x,\xi)\ra}+|\Xi(s;x,\xi)|\\
\leq |\widebar{\Xi}(x,\xi)|\leq \frac{C}{\la X(s;x,\xi)\ra}+|\Xi(s;x,\xi)|
\leq \frac{C}{\la X(s;x,\xi)\ra}+C_1|\xi|,
\end{align}
which, since $\la X(s)\ra\rightarrow +\infty$ as $s\rightarrow +\infty$, passing to the limit gives
\[ C_1^{-1}|\xi|\leq |\widebar{\Xi}(x,\xi)|\leq C_1 |\xi|.\]
This concludes the proof.
\endproof

We now recall the following version of Gronwall's lemma (see Theorem 1.10 in \cite{T}).
\begin{proposition}\label{prop.Gronwall}
Let $y:[0,+\infty]\rightarrow \R_+$ be a continuous nonnegative function such that 
\[
y(s) \leq c+\int_0^s k(\sigma)y(\sigma)d\sigma, \quad \forall s \geq 0,
\]
where $c\geq 0$ and $k:[0,+\infty]\rightarrow \R$ is a nonnegative continuous function that satisfies 
\[
\int_0^{+\infty}k(\sigma)d\sigma=C_0<+\infty.
\]
Then, for all $s \geq 0$, we have 
\[
y(s)\leq ce^{C_0}.
\]
\end{proposition}

We now recall the multivariate Fa\'a Di Bruno formula (see \cite{CS}). To this end we first introduce some preliminary notation. 

For $\mu,\nu \in \mathbb{N}_0^n$ we write $\mu \prec \nu$ if: $|\mu | < |\nu| $, or $|\mu|=|\nu|$ and $\mu_1 < \nu_1$, or $|\mu|=|\nu|$, $\mu_1=\nu_1, \dots, \mu_k=\nu_k$ and $\mu_{k+1}<\nu_{k+1}$ for some $k=1,\dots,n$. 
\begin{theorem}\label{theo.fdb}
    Let $g\in C^\infty(\R^n_x;\R^m_y)$ and $f\in C^\infty (\R^m_y;\mathbb{R}^N_z)$. Let
    \[ h(x_1,\dots,x_n)=f(g_1(x_1,\dots,x_n), \dots, g_m(x_1,\dots,x_n)).
    \]
    Then, for $x_0 \in \R^n$ and for $\alpha \in \mathbb{N}_0^n$, we have 
    \[
    \partial^\alpha_x h(x_0)=\sum_{1\leq |\gamma| \leq |\alpha|}\partial_y^\gamma f(g(x_0)) \sum_{(\lambda,\ell) \in p(\alpha,\gamma)} c_{\alpha,\lambda,\ell} \prod_{j=1}^{|\alpha|} [\partial^{\ell_j}_xg]^{\lambda_j}(x_0),
    \]
    for $c_{\alpha,\lambda,\ell}>0$ constant,  
    \[
    \begin{split}
    p(\alpha,\gamma)&= \Bigl\lbrace (\lambda,\ell)=(\lambda_1,\dots,\lambda_{|\alpha|}, \ell_1,\dots,\ell_{|\alpha|}); \ \lambda_j \in \mathbb{N}_0^m, \ \ell_j \in \mathbb{N}_0^n, \ \forall j=1,\dots, |\alpha|,\\
    &\quad \quad \text{for some $1\leq s \leq |\alpha|$ one has} \\
    & \quad \quad \lambda_j=\ell_j=0, \quad \forall j=1,\dots, |\alpha|-s, \\
    & \quad \quad |\lambda_j| >0, \quad   \forall j=|\alpha|-s+1 \leq 1\leq |\alpha|, \\
    & \quad \quad 0 \prec \ell_{|\alpha|-s+1} \prec \dots \prec \ell_{|\alpha|}, \\
    & \quad \quad \text{and} \\
    & \quad \quad \sum_{j=1}^{|\alpha|} \lambda_j=\gamma, \quad \sum_{j=1}^{|\alpha|} |\lambda_j|\ell_j= \alpha \Bigr\rbrace,
    \end{split}
    \]
    and 
    \[
    [\partial_x^{\ell_j}g]^{\lambda_j}(x_0)=\prod_{i=1}^m (\partial_x^{\ell_j}g_i(x_0))^{\lambda_{j,i}}.
    \]
\end{theorem}

\end{document}